 \documentclass{amsart}

 \usepackage{amscd, amsmath, amsthm, amssymb,multicol}
 \usepackage[curve,all]{xy}
 \xyoption{curve}
 \xyoption{line}

\theoremstyle{plain}

\newtheorem{theorem}{Theorem}[section]
\newtheorem{proposition}[theorem]{Proposition}
\newtheorem{corollary}{Corollary}[section]

\newtheorem{lemma}[theorem]{Lemma}

\theoremstyle{definition}

\newtheorem{example}[theorem]{Example}
\newtheorem{examples}[theorem]{Examples}

\theoremstyle{remark}
\newtheorem{remark}[theorem]{Remark}

\numberwithin{equation}{section}

\newcommand{\calE}{\mathcal{E}}

\newcommand{\calH}{\mathcal{H}}

\newcommand{\calU}{\mathcal{U}}
\newcommand{\calO}{\mathcal{O}}

\newcommand{\bbB}{\mathbb{B}}

\newcommand{\bbF}{\mathbb{F}}
\newcommand{\bbC}{\mathbb{C}}

\newcommand{\bbP}{\mathbb{P}}
\newcommand{\bbQ}{\mathbb{Q}}
\newcommand{\bbR}{\mathbb{R}}

\newcommand{\bbZ}{\mathbb{Z}}

\newcommand{\bfa}{\mathbf{a}}

\newcommand{\bfe}{\mathbf{e}}

\newcommand{\bfF}{\mathbf{F}}

\newcommand{\la}{\langle}
\newcommand{\ra}{\rangle}
\newcommand{\cha}{\textup{char}}

\newcommand{\SL}{\textrm{SL}}
\newcommand{\Aut}{\text{Aut}}
\newcommand{\Hom}{\text{Hom}}
\newcommand{\GL}{\text{GL}}

\newcommand{\PSL}{\text{PSL}}

\newcommand{\Or}{\textup{O}}
\newcommand{\PO}{\textup{PO}}
\newcommand{\rank}{\text{rank}}

\newcommand{\PGL}{\textup{PGL}}

\newcommand{\id}{\textup{id}}

\newcommand{\U}{\textup{U}}
\newcommand{\Lin}{\textup{Lin}}
\newcommand{\Ref}{\textup{Ref}}
\newcommand{\Tors}{\textup{Tors}}
\newcommand{\diag}{\textup{diag}}
\newcommand{\mult}{\textup{mult}}
\newcommand{\sign}{\textup{sign}}
\newcommand{\Jac}{\textup{Jac}}
\newcommand{\Sp}{\textup{Sp}}

\newcommand{\PU}{\textup{PU}}
\newcommand{\bbCH}{H_\bbC^n}
 
\begin{document}
 \title{Reflection Groups in Algebraic Geometry}

\author{Igor V. Dolgachev}

\address{Department of Mathematics, University of Michigan, 525 E. University Av., Ann Arbor, Mi, 49109}
\email{idolga@umich.edu}
\thanks{The  author was supported in part by NSF grant 0245203.}



\dedicatory{To Ernest Borisovich Vinberg}

\begin{abstract}
After  a brief exposition of the theory of discrete reflection groups in spherical, euclidean and hyperbolic geometry as well as their analogs in complex spaces we present a survey of appearances of these groups in various areas of algebraic geometry.
\end{abstract}

\maketitle

\section{Introduction}

The notion of a reflection in an euclidean space is one of the fundamental notions of symmetry of geometric figures and does not need an introduction. The theory of discrete groups of motions generated by reflections originates in the study of plane regular polygons and space polyhedra which goes back to ancient mathematics. Nowadays it is hard to find a mathematician  who has  not encountered reflection groups in his area of research. Thus a geometer sees them as examples of discrete groups of isometries of Riemannian spaces of constant curvature or examples of special  convex polytopes. An algebraist  finds them in group theory, especially in the theory of Coxeter groups,   invariant theory and representation theory. A combinatorialist may see them in the theory of arrangements of hyperplanes and combinatorics of permutation groups.  A number theorist meets them in arithmetic theory of quadratic forms and modular forms. For a topologist they turn up in the study of hyperbolic real and complex manifolds,  low-dimensional topology and singularity theory.  An analyst sees them in the theory of hypergeometric functions and automorphic forms, complex higher-dimensional dynamics and ordinary differential equations. All of the above  and much more appears in algebraic geometry. The goal of this survey is to explain some of ``much more''. 

One finds an extensive account of the history of the theory of reflection groups in euclidean and spherical spaces in Bourbaki's ``Groupes et Alg\`ebres de Lie'', Chapters IV-VI. According to this account the modern theory originates from the works of geometers A. Moebius and L. Schl\"afli in the middle of the 19th century, then extended  and applied to the theory of Lie algebras in the works of E. Cartan and W. Killing at the end of the same century,  and culminated in the works of H.S.M Coxeter \cite{Cox}. The first examples of reflection groups in hyperbolic plane  go back to  F. Klein and H. Poncar\'e at the end of the 19th century. 

It is a not a general knowledge  that reflection groups, finite and infinite, appear 
 in the works of S. Kantor in 1885-1895 \cite{Ka} on classification of subgroups of the Cremona group of birational transformations of the complex projective plane \cite{Ka}. He realized their reflection action   in the  cohomology space of rational algebraic surfaces obtained as blow-ups of the plane. In this way   all Weyl groups of type $A_2\times A_1,A_4,D_5,E_6,E_7,E_8$ appear naturally when the number of points blown-up is between $3$ and $8$. The last three groups appeared much earlier in algebraic geometry as the group of 27 lines on a cubic surface (type $E_6$),  the group of bitangents of a plane quartic (type $E_7$), and the group of tritangent planes of a space sextic of genus 4 (type $E_8$).  They  were widely known among algebraists since the appearance of C. Jordan's  ``Trait\'e  des substitutions'' in 1870. In 1910 P. Schoute discovered a convex polytope in six-dimensional space whose vertices are in a bijective correspondence with  27 lines on a cubic surface and the group of symmetries is isomorphic to the group of 27 lines \cite{Schoute}. A similar polytope in seven-dimensional space was found for the group of 28 bitangents of a plane quartic by Coxeter \cite{Coxeter2}.  The relationship between this six-dimensional space and the cohomology space of the blow-up of the plane at 6 points was explained by P. Du Val in the thirties of the last century. He also showed that all Kantor groups are reflection groups in euclidean, affine or hyperbolic spaces \cite{DV2}, \cite{DuVal2}. 
 
 The Coxeter diagram of the  Weyl group of type $E_8$ is of the form 

\medskip
\begin{figure}[ht]

\xy @={(0,0),(10,0),(20,0),(30,0),(40,0),(50,0),(60,0)}@@{*{\bullet}};
(0,0)*{};(60,0)**{}**\dir{-};
(-30,0)*{};
(20,0)*{};(20,-10)*{\bullet}**\dir{-};\endxy

\caption{}\label{d}
\end{figure}

It appears when the number of points blown up is equal to 8.  When the number of points is 9, we get a reflection group in affine space of dimension 8, the affine Weyl group of type $E_8$. Its Coxeter diagram is obtained by adding one point in the long arm of the diagram above. Starting from 10 points one gets Coxeter  groups in hyperbolic space with Coxeter diagram of type $E_n$ (extending the long arm of the  diagram on Fig. 1). In 1917, generalizing Kantor's work, Arthur Coble introduced  the notion of a regular Cremona transformation of a higher-dimensional projective space and considered more general Coxeter groups of type $W(2,p,q)$ \cite{Co}. Its Coxeter diagram is obtained from extending the two upper arms of the diagram. A modern account of Coble's theory is given in my book with D. Ortland \cite{DO}. Very recently, S. Mukai \cite{Mu} was able to extend Coble's construction to include all Coxeter groups of type $W(p,q,r)$ with Coxeter diagram obtained from the  above diagram by extending all its arms. 
 
Another class of algebraic surfaces where reflection groups in hyperbolic spaces arise naturally is the class of surfaces of type K3. An example of such a surface is a nonsingular quartic surface in complex three-dimensional projective space. F. Severi  gave the first example of a quartic  surface with explicitly  computed infinite group of birational automorphisms; the group turns out to be isomorphic to a plane reflection group \cite{Severi}.  In 1972   I.I. Pyatetsky-Shapiro and I. R. Shafarevich, answering a question of  A. Weil, proved that the complex structure of an algebraic (polarized) K3-surface is determined uniquely by the linear functional on its second cohomology space obtained by integrating  a nowhere vanishing holomorphic 2-form \cite{PS}, \cite{Sha}. They called  this result a Global Torelli Theorem. As a corollary of this result they showed that the group of automorphisms of a K3 surface is isomorphic, up to a finite group, to the quotient of the orthogonal group of the integral quadratic form defined on the group of integral algebraic 2-cycles modulo the subgroup $\Gamma$ generated by reflections in the cohomology classes of smooth rational curves  lying on the surface. Thus they reduced the  question of finiteness of the automorphism group to the question of finiteness of the volume of the fundamental polyhedron of $\Gamma$ in a real hyperbolic space. V. Nikulin \cite{Ni1} and E. Vinberg \cite{Vi1} determined which isomorphism types of integral quadratic forms so arise and have the property that the fundamental polyhedron in question has finite volume. This solves, in principle,  the problem of classification of  fields of algebraic dimension 2 over $\bbC$ whose group of automorphisms over $\bbC$ is infinite. 

In the thirties of the last century   Patrick Du Val \cite{DV} found the appearance of Coxeter diagrams in resolution of certain types of singularities on algebraic surfaces (nowadays going under the different names: simple singularities, ADE singularities, Du Val singularities, double rational points, Gorenstein quotient singularities, Klein singularities). However Du Val  did not find any reflection groups associated to these singularities. A conjectural relation to reflection groups and simple Lie algebras was suggested by A. Grothendieck in the sixties, and was confirmed by a construction of E. Brieskorn \cite{Br} (full details appeared in \cite{Sl2}).

The theory of finite complex reflection groups  was developed by G. C. Shephard and A. Todd in 1954  as a follow-up of the classical work on groups of projective transformations generated by homologies (see \cite{ST}).  Some examples of the arrangements of reflecting hyperplanes and the hypersurfaces defined by polynomial invariants of the groups  were known in classical geometry since the 19th century. 

Infinite  reflection groups of finite covolume in complex affine spaces were classified by V. Popov in 1982 \cite{Po}. They appear in the theory of compactification of versal deformation of simple elliptic singularities \cite{Lo1} and surface singularities with symmetries \cite{Go4}.

The most spectacular is the appearance of  reflection groups in complex hyperbolic spaces of dimension $> 1$. Extending the work of H. Terada \cite{Te}, P. Deligne and G. Mostow \cite{DM}, \cite{Mo} classified all hypergeometric functions whose monodromy groups $\Gamma$ are discrete reflection groups of finite covolume in a complex ball $\bbB^r$  (complex hyperbolic crystallographic groups, c.h.c. groups for short). The compactified orbit spaces $\overline{\bbB^r/\Gamma}$ turned out to be isomorphic to some  geometric invariant quotients  $\bbP^1(\bbC)^{r+3}/\!/\PGL(2,\bbC)$ for $r\le 9$.  No other  c.h.c.  groups  in $\bbB^r$ have been discovered until a few years ago (except one missed case in Deligne-Mostow's list noticed by W. Thurston \cite{Th}). The first new c.h.c. group in $\bbB^4$ appeared in a beautiful construction of D. Allcock, J. Carlson and D. Toledo of a complex ball uniformization of the moduli space of cubic surfaces \cite{ACT}. Later, using a similar uniformization construction for moduli spaces of other Del Pezzo surfaces,  new examples of c.h.c, groups were found in dimensions 6 and 8 \cite{HL}, \cite{Ko1},\cite{Ko2}. All these groups are commensurable with some of the Deligne-Mostow groups. A recent  work of Allcock, Carlson and Toledo \cite{ACT2} (see also \cite{Lo4}) on complex ball uniformization of the moduli space of cubic hypersurfaces in $\bbP^4$ produces a new complex reflection group in dimension 10. A generalization of the Deligne-Mostow theory due to W. Couwenberg, G. Heckaman and E. Looijenga \cite{CHL} gives another examples of complex reflection crystallographic groups.  A c.h.c. group in a record high dimension 13 was constructed by D. Allcock \cite{Al}. No geometrical interpretation so far is known for the corresponding ball quotients. 

The above discussion outlines the contents of the present paper.  As is the case with any survey paper, it is incomplete and the omitted material is either due to author's ignorance, or poor memory,  or limited  size of the paper.

\medskip\noindent
{\bf Acknowledgements.} This paper is dedicated  to Ernest Borisovich  Vinberg,  one of the  heroes in  the theory of reflection groups.  His  lectures  for high school children in Moscow were influential (without his knowledge) in my decision to become a mathematician. 

The paper is an expanded version of my colloquium lecture at the University di Roma Terzo in May 2006. I am thankful to Alessandro Verra for giving me an opportunity to give this talk, and hence to write the paper. 
I am very grateful to Daniel Allcock, Victor Goryunov  and the referee for numerous critical comments on earlier versions of the paper.

 \section{Real reflection groups} 
\subsection{Elementary introduction} The idea of a reflection transformation $r_H$ with respect to a mirror line $H$ is of course very familiar. 
 A picture on the plane is 

\vspace{14pt}

\begin{figure}[h] 

\xy (0,0)*{};(40,0)*{}**\dir{-};
(20,10)*{\bullet};(20,-10)*{\bullet}**\dir{-};
(22,10)*{p};(25,-10)*{r_H(p)};
(44,0)*{H};
(-40,0)*{};(-40,10)*{};(-40,-10)*{};\endxy
\caption{}  
\end{figure}

  Now suppose we have two mirror lines $H_1$ and $H_2$. Each line divides the plane into the disjoint union of two \emph{halfplanes} $H_i^\pm$.

  A choice of  halfplanes, say $H_1^+,H_2^+$ defines the \emph{angle} $ H_1^+\cap H_2^+$ with measure $\phi = \angle (H_1^+,H_2^+)$.  Here $\phi = 0$ if and only if $H_1^-\cap H_2^- =\emptyset.$

\vspace{14pt}
\begin{figure}[h] 
\xy (-20,-10)*{};(26,13)*{}**\dir{-};
(-20,0)*{};(30,)*{}**\dir{-};
(30,16)*{H_2};(35,0)*{H_1};(3,7)*{H_2^+};(5,-5)*{H_1^+};
(12,1)*\cir<15pt>{u^l};(12,3)*{\phi};
(-50,0)*{};
\endxy

\caption{}  
\end{figure}

Let $s_1 = r_{H_1}, s_2 = r_{H_2}$. The composition $s_2s_1$ is the counterclockwise rotation about the angle $2\phi$ if $\phi \ne 0$ and a translation if $\phi = 0$:

\vspace{14pt}
\begin{figure}[h]

\begin{multicols}{2}{
\xy (-20,-10)*{};(26,13)*{}**\dir{-};
(-20,0)*{};(30,)*{}**\dir{-};
(30,16)*{H_2};(35,0)*{H_1};(3,7)*{H_2^+};(5,-5)*{H_1^+};
(0,0)*{};(25,5)*{\bullet}**\dir{-};
(0,0)*{};(25,-5)*{\bullet}**\dir{-};
(0,0)*{};(20,18)*{\bullet}**\dir{-};
(25,22)*{s_2s_1(p)};(30,-5)*{p};(32,5)*{s_1(p)}
\endxy

\xy (40,5)*{};(80,5)*{}**\dir{-};
(40,-5)*{};(80,-5)*{}**\dir{-};
(60,-7.5)*{\bullet};(60,12.5)*{\bullet}*{};
(60,-2.5)*{\bullet};(66,-2.25)*{s_1(p)};
(67,12.5)*{s_2s_1(p)};(63,-7.25)*{p};

\endxy}
\end{multicols}\caption{}
\end{figure}

Let $G$ be the group generated by the two reflections $s_1,s_2$. We assume that $H_1\ne H_2$, i.e. $s_1\ne s_2$. The following two cases may occur:

\smallskip
\emph{Case 1}: The angle $\phi$ is of the form $n\pi/m$ for some rational number $r = n/m$. In the following  we assume that $m=\infty$ if $\phi = 0$.

In this case $s_2s_1$ is the rotation about the angle $2n\pi/m$ and hence 
$$(s_2s_1)^m = \textup{identity}.$$

The group $G$ is isomorphic to the  finite dihedral group $D_{2m}$ of order $2m$ (resp. \emph{infinite dihedral group} $D_\infty$ if $m = \infty$)  with presentation 
$$<s_1,s_2|s_1^2 = s_2^2 = (s_1s_2)^m = 1>.$$

 It acts as a discrete group of motions of the plane with fundamental domain  $H_1^-\cap H_2^-$.

Observe that the same group is generated by reflections with respect to the lines forming the angle obtained from the angle $\angle(H_1^-,H_2^-)$ subdividing it in $n$ equal parts. So, we may assume that $\phi = \pi/m$.

\smallskip
\emph{Case 2}: The angle $\phi$ is not of the form $r\pi$ for any rational $r$.

In this case $s_2s_1$ is of infinite order,  $G$  is isomorphic to  $D_\infty$ but 
it does not act discretely.

Now suppose we have a convex polygon given as the intersection of a finite set of halfplanes
$$P =  \bigcap_{i=1}^rH_i^-.$$

We assume that the interior $P^o$ is not empty and the set $\{H_1,\ldots,H_r\}$ is minimal in the sense that one cannot delete any of the halfplanes without changing $P$. 

More importantly, we assume that 
$$\angle(H_i^-,H_j^-) = \pi/m_{ij}$$
 for some positive integer $m_{ij}$ or equal to 0 ($m_{ij} = \infty$).
 
Let $G$ be the group generated by reflections with mirror lines $H_i$. It is a discrete group of motions of the plane. 
The polygon $P$ is a fundamental domain  of $G$ in the plane.

Conversely any discrete group of motions of the plane generated by reflections is obtained in this way.

Let $p_1,\ldots,p_r$ be the vertices of the polygon $P$. We may assume that  $p_i = H_i\cap H_{i+1},$ where $H_{r+1} = H_1$. Let $m_i = m_{ii+1}$. Since
$$\sum_{i=1}^r (\pi/m_{i}) = (r-2)\pi$$
we have
$$\sum_{i=1}^r\frac{1}{m_i} = r-2.$$

The only solutions for $(r;m_1,\ldots,m_r)$ are
$$(3;2,3,6), \quad (3;2,4,4), \quad (3;3,3,3) , \quad (4;2,2,2,2).$$

\begin{figure}[h]

\begin{multicols}{2}{

\xy @={(0,0),(15,0),(30,0),(0,15),(0,30),(30,0),(15,15),(15,30),(30,30),(30,15),(30,30)}@@{*{\bullet}};
(0,0)*{};(30,0)*{}**\dir{-};(0,0)*{};(0,30)*{}**\dir{-};(0,0)*{};(30,30)*{}**\dir{-};(0,30)*{};(30,0)*{}**\dir{-};
(0,30)*{};(30,30)*{}**\dir{-};(0,15)*{};(30,15)*{}**\dir{-};(15,0)*{};(15,30)*{}**\dir{-};(30,0)*{};(30,30)*{}**\dir{-};(0,7.5)*{};(30,7.5)*{}**\dir{-};(7.5,0)*{};(7.5,30)*{}**\dir{-};(0,22.5)*{};(30,22.5)*{}**\dir{-};(22.5,0)*{};(22.5,30)*{}**\dir{-};(0,15)*{};(15,30)*{}**\dir{-};(0,15)*{};(15,0)*{}**\dir{-};(15,0)*{};(0,15)*{}**\dir{-};(30,15)*{};(15,30)*{}**\dir{-};(15,0)*{};(30,15)*{}**\dir{-};
(15.25,14.75)*{};(15.25,7.5)*{}**\dir2{-};
(15.5,14.5)*{};(15.5,7.5)*{}**\dir2{-};
(15.75,14.25)*{};(15.75,7.5)*{}**\dir2{-};
(16,14)*{};(16,7.5)*{}**\dir2{-};
(16.25,13.75)*{};(16.25,7.5)*{}**\dir2{-};
(16.5,13.5)*{};(16.5,7.5)*{}**\dir2{-};
(16.75,13.25)*{};(16.75,7.5)*{}**\dir2{-};
(17,13)*{};(17,7.5)*{}**\dir2{-};
(17.25,12.75)*{};(17.25,7.5)*{}**\dir2{-};
(17.5,12.5)*{};(17.5,7.5)*{}**\dir2{-};
(17.75,12.25)*{};(17.75,7.5)*{}**\dir2{-};
(18,12)*{};(18,7.5)*{}**\dir2{-};
(18.25,11.75)*{};(18.25,7.5)*{}**\dir2{-};
(18.5,11.5)*{};(18.5,7.5)*{}**\dir2{-};
(18.75,11.25)*{};(18.75,7.5)*{}**\dir2{-};
(19,11)*{};(19,7.5)*{}**\dir2{-};
(19.25,10.75)*{};(19.25,7.5)*{}**\dir2{-};
(19.5,10.5)*{};(19.5,7.5)*{}**\dir2{-};
(19.75,10.25)*{};(19.75,7.5)*{}**\dir2{-};
(20,10)*{};(20,7.5)*{}**\dir2{-};
(20.25,9.75)*{};(20,7.5)*{}**\dir2{-};
(20.5,9.5)*{};(20.5,7.5)*{}**\dir2{-};
(20.75,9.25)*{};(20.75,7.5)*{}**\dir2{-};
(21,9)*{};(21,7.5)*{}**\dir2{-};
(21.25,8.75)*{};(21.25,7.5)*{}**\dir2{-};
(21.5,8.5)*{};(21.5,7.5)*{}**\dir2{-};
(21.75,8.25)*{};(21.75,7.5)*{}**\dir2{-};
(22,8)*{};(22,7.5)*{}**\dir2{-};
(22.25,7.75)*{};(22.25,7.5)*{}**\dir2{-};
(-20,0)*{};(15,-5)*{(3;2,4,4), G\cong \bbZ^2\rtimes D_8};

\endxy

\xy @={(-12.9,7.5),(-12.9,-7.5),(0,15),(0,-15),(12.9,7.5),(12.9,-7.5),(0,0)}@@{*{\bullet}};
(-17.4,0)*+{}\PATH ~={**\dir{-}}
'(-7.5,15)*+{}' (7.5,15)*+{}' (17.4,0)*+{}' (7.5,-15)*+{}' (-7.5,-15)*+{}' (-17.4,0)*+{}';
(-17.4,0)*{};(17.4,0)*{}**\dir{-};(0,15)*{};(0,-15)*{}**\dir{-};(-7.5,15)*{};(7.5,-15)*{}**\dir{-};
(7.5,15)*{};(-7.5,-15)*{}**\dir{-};
(-12.9,-7.5)*+{}\PATH~={**\dir{-}}
'(-12.9,7.5)*+{}',(0,15)*+{}',(12.9,7.5)*+{}',(12.9,-7.5)*+{}',(0,-15)*+{}',(-12.9,-7.5)*+{}';
(-12.9,-7.5)*{};(-12.9,7.5)*{}**\dir{-};
(12.9,-7.5)*{}; (-12.9,-7.5)*{}**\dir{-};(-12.9,-7.5)*{};(0,15)*{}**\dir{-}';(0,15)*{};(12.9,-7.5)*{}**\dir{-};
(-12.9,7.5)*{};
(0,-15)*{}**\dir{-};
(0,-15)*{};(12.9,7.5)*{}**\dir{-};
(12.9,7.5)*{};(-12.9,7.5)*{}**\dir{-};
(12.9,-7.5)*{};(-12.9,7.5)*{}**\dir{-};
(-12.9,-7.5)*{};(12.9,7.5)*{}**\dir{-};
(0,-20)*{(3;2,3,6),G\cong \bbZ^2\rtimes D_{12}};
(.25,-.5)*{};(.25,-7.5)*{}**\dir2{-};
(.5,-1)*{};(.5,-7.5)*{}**\dir2{-};
(.75,-1.5)*{};(.75,-7.5)*{}**\dir2{-}; 
(1,-2)*{};(1,-7.5)*{}**\dir2{-}; 
(1.25,-2.5)*{};(1.25,-7.5)*{}**\dir2{-}; 
(1.5,-3)*{};(1.5,-7.5)*{}**\dir2{-}; 
(1.75,-3.5)*{};(1.75,-7.5)*{}**\dir2{-}; 
(2,-4)*{};(2,-7.5)*{}**\dir2{-}; 
(2.25,-4.5)*{};(2.25,-7.5)*{}**\dir2{-}; 
(2.5,-5)*{};(2.5,-7.5)*{}**\dir2{-}; 
(2.75,-5.5)*{};(2.75,-7.5)*{}**\dir2{-}; 
(3,-6)*{};(3,-7.5)*{}**\dir2{-}; 
(3.25,-6.5)*{};(3.25,-7.5)*{}**\dir2{-}; 
\endxy}
\end{multicols}

\begin{multicols}{2}{

\xy 
@={(-12.9,7.5),(-12.9,-7.5),(0,15),(0,-15),(12.9,7.5),(12.9,-7.5),(0,0)}@@{*{\bullet}};
(-17.4,0)*+{}\PATH ~={**\dir{-}}
'(-7.5,15)*+{}' (7.5,15)*+{}' (17.4,0)*+{}' (7.5,-15)*+{}' (-7.5,-15)*+{}' (-17.4,0)*+{}';
(-17.4,0)*{};(17.4,0)*{}**\dir{-};
(7.5,15)*{};(-7.5,-15)*{}**\dir{-};
(12.9,-7.5)*{}; (-12.9,-7.5)*{}**\dir{-};(-12.9,-7.5)*{};(0,15)*{}**\dir{-}';(0,15)*{};(12.9,-7.5)*{}**\dir{-};
(-12.9,7.5)*{};
(0,-15)*{}**\dir{-};
(0,-15)*{};(12.9,7.5)*{}**\dir{-};
(12.9,7.5)*{};(-12.9,7.5)*{}**\dir{-};
(-35,0)*{}; (0,-20)*{(3;3,3,3), G\cong \bbZ^2\rtimes S_3};
(.25,-.5)*{};(.25,-7.5)*{}**\dir3{-};
(.5,-1)*{};(.5,-7.5)*{}**\dir3{-};
(.75,-1.5)*{};(.75,-7.5)*{}**\dir3{-}; 
(1,-2)*{};(1,-7.5)*{}**\dir3{-}; 
(1.25,-2.5)*{};(1.25,-7.5)*{}**\dir3{-}; 
(1.5,-3)*{};(1.5,-7.5)*{}**\dir3{-}; 
(1.75,-3.5)*{};(1.75,-7.5)*{}**\dir3{-}; 
(2,-4)*{};(2,-7.5)*{}**\dir3{-}; 
(2.25,-4.5)*{};(2.25,-7.5)*{}**\dir3{-}; 
(2.5,-5)*{};(2.5,-7.5)*{}**\dir3{-}; 
(2.75,-5.5)*{};(2.75,-7.5)*{}**\dir3{-}; 
(3,-6)*{};(3,-7.5)*{}**\dir3{-}; 
(3.25,-6.5)*{};(3.25,-7.5)*{}**\dir3{-}; 
(-.25,-.5)*{};(-.25,-7.5)*{}**\dir3{-};
(-.5,-1)*{};(-.5,-7.5)*{}**\dir3{-};
(-.75,-1.5)*{};(-.75,-7.5)*{}**\dir3{-}; 
(-1,-2)*{};(-1,-7.5)*{}**\dir3{-}; 
(-1.25,-2.5)*{};(-1.25,-7.5)*{}**\dir3{-}; 
(-1.5,-3)*{};(-1.5,-7.5)*{}**\dir3{-}; 
(-1.75,-3.5)*{};(-1.75,-7.5)*{}**\dir3{-}; 
(-2,-4)*{};(-2,-7.5)*{}**\dir3{-}; 
(-2.25,-4.5)*{};(-2.25,-7.5)*{}**\dir3{-}; 
(-2.5,-5)*{};(-2.5,-7.5)*{}**\dir3{-}; 
(-2.75,-5.5)*{};(-2.75,-7.5)*{}**\dir3{-}; 
(-3,-6)*{};(-3,-7.5)*{}**\dir3{-}; 
(-3.25,-6.5)*{};(-3.25,-7.5)*{}**\dir3{-}; 
\endxy

\xy @={(0,0),(15,0),(30,0),(0,15),(0,30),(30,0),(15,15),(15,30),(30,30),(30,15),(30,30)}@@{*{\bullet}};
(0,0)*{};(30,0)*{}**\dir{-};(0,0)*{};(0,30)*{}**\dir{-};
(0,30)*{};(30,30)*{}**\dir{-};(0,15)*{};(30,15)*{}**\dir{-};(15,0)*{};(15,30)*{}**\dir{-};(30,0)*{};(30,30)*{}**\dir{-};(0,7.5)*{};(30,7.5)*{}**\dir{-};(7.5,0)*{};(7.5,30)*{}**\dir{-};(0,22.5)*{};(30,22.5)*{}**\dir{-};(22.5,0)*{};(22.5,30)*{}**\dir{-};
(15,-5)*{(4;2,2,2,2), G\cong \bbZ^2\rtimes D_4};
(15.25,15)*{};(15.25,7.5)*{}**\dir3{-};
(15.5,15)*{};(15.5,7.5)*{}**\dir3{-};
(15.75,15)*{};(15.75,7.5)*{}**\dir3{-};
(16,15)*{};(16,7.5)*{}**\dir3{-};
(16.25,15)*{};(16.25,7.5)*{}**\dir3{-};
(16.5,15)*{};(16.5,7.5)*{}**\dir3{-};
(16.75,15)*{};(16.75,7.5)*{}**\dir3{-};
(17,15)*{};(17,7.5)*{}**\dir3{-};
(17.25,15)*{};(17.25,7.5)*{}**\dir3{-};
(17.5,15)*{};(17.5,7.5)*{}**\dir3{-};
(17.75,15)*{};(17.75,7.5)*{}**\dir3{-};
(18,15)*{};(18,7.5)*{}**\dir3{-};
(18.25,15)*{};(18.25,7.5)*{}**\dir3{-};
(18.5,15)*{};(18.5,7.5)*{}**\dir3{-};
(18.75,15)*{};(18.75,7.5)*{}**\dir3{-};
(19,15)*{};(19,7.5)*{}**\dir3{-};
(19.25,15)*{};(19.25,7.5)*{}**\dir3{-};
(19.5,15)*{};(19.5,7.5)*{}**\dir3{-};
(19.75,15)*{};(19.75,7.5)*{}**\dir3{-};
(20,15)*{};(20,7.5)*{}**\dir3{-};
(20.25,15)*{};(20,7.5)*{}**\dir3{-};
(20.5,15)*{};(20.5,7.5)*{}**\dir3{-};
(20.75,15)*{};(20.75,7.5)*{}**\dir3{-};
(21,15)*{};(21,7.5)*{}**\dir3{-};
(21.25,15)*{};(21.25,7.5)*{}**\dir3{-};
(21.5,15)*{};(21.5,7.5)*{}**\dir3{-};
(21.75,15)*{};(21.75,7.5)*{}**\dir3{-};
(22,15)*{};(22,7.5)*{}**\dir3{-};
(22.25,15)*{};(22.25,7.5)*{}**\dir3{-};

\endxy}
\end{multicols}
\caption{}
\end{figure}

\subsection{Spaces of constant curvature}
The usual euclidean plane is an example of a 2-dimensional space of zero constant curvature. Recall that a \emph{space of constant curvature} is a simply connected Riemannian homogeneous space $X$ such that  the isotropy subgroup of its group of isometries $\textup{Iso}(X)$ at each point coincides with the full orthogonal group of the tangent space. Up to isometry and rescaling the metric, there are three spaces of constant curvature of fixed dimension $n$.

\begin{itemize}
\item The \emph{euclidean space} $E^n$  with $\textup{Iso}(X)$ equal to the affine orthogonal group 
$\textup{AO}^n = \bbR^n\rtimes \Or(n)$, 
\item 
The \emph{$n$-dimensional sphere} 
$$S^n = \{(x_0,\ldots,x_n)\in \bbR^{n+1}:x_0^2+\ldots+x_n^2 = 1\}$$
with $\textup{Iso}(X)$ equal to the orthogonal group $\Or(n+1)$,
\item The \emph{hyperbolic} (or \emph{Lobachevsky}) space 
$$H^n = \{(x_0,\ldots,x_n)\in \bbR^{n+1}, -x_0^2+x_1^2+\ldots+x_n^2 = -1, x_0 > 0\}$$
with $\textup{Iso}(X)$ equal to the subgroup $\Or(n,1)^+$ of index 2 of the orthogonal group $\Or(n,1)$  which consists of transformations of spinor norm 1, that it is, transformations that can be written as a product of reflections in vectors with positive norm. The Riemannian metric is induced by the hyperbolic metric in $\bbR^{n+1}$
$$ds^2 = -dx_0^2+dx_1^2+\ldots+dx_n^2.$$
\end{itemize}

We will be  using a \emph{projective model} of $H^n$, considering $H^n$ as the image of the subset  
$$C = \{(x_0,\ldots,x_n)\in \bbR^{n+1}, -x_0^2+x_1^2+\ldots+x_n^2 < 0\}$$
in the projective space $\bbP^n(\bbR)$. The isometry group of the projective model is naturally identified with the group $\PO(n,1)$. By choosing a representative of a point from  $H^n$ with $x_0 = 1$, we can identify $H^n$ with the real ball 
$$K^n:\{(x_1,\ldots,x_n):|x|^2 = x_1^2+\ldots+x_n^2 < 1\}.$$ (the \emph{Klein model}). The metric is given by
$$ds^2 = \frac{1}{1-|x|^2}\sum_{i=1}^ndx_i^2+\frac{1}{(1-|x|^2)^2}(\sum_{i=1}^nx_idx_i)^2.$$
The closure of $H^n$ in $\bbP^n(\bbR)$ is equal to the image of the set 
$$\bar{C} = \{(x_0,\ldots,x_n)\in \bbR^{n+1}, -x_0^2+x_1^2+\ldots+x_n^2 \le 0\}$$
in $\bbP^n(\bbR)$. The boundary is called the \emph{absolute}.

One defines the notion of a \emph{hyperplane} in  a space of constant curvature.  If $X = E^n$, a hyperplane is an affine hyperplane. If $X= S^n$, a hyperplane is the intersection of $S^n$ with a linear hyperplane in $\bbR^{n+1}$ (a great circle when $n = 2$). If $X = H^n$, a hyperplane is the non-empty intersection of $H^n$ with a projective hyperplane in $\bbP^n(\bbR)$.

Each hyperplane $H$ in $E^n$ is a translate $a+L = \{x+a, x\in L\}$ of a unique  linear hyperplane $\tilde{H}$ in the corresponding standard euclidean space $V = \bbR^n$.  If $X^n = S^n$ or $H^n$, then a hyperplane $H$ is uniquely defined by a linear hyperplane $\tilde{H}$ in  $V = \bbR^{n+1}$ equipped with the standard symmetric bilinear form of Sylvester signature $(t_+,t_-) =  (n+1,0)$, or $(n,1)$.  

Any point $x\in V$ can be written uniquely in the form
$$x = h+v,$$
 where $h\in H$ and $v\in V$ is orthogonal to $\tilde{H}$.  We define  a\emph{reflection} with mirror hyperplane $H$ by the formula
 $$r_H(x) = h-v.$$
 
One can also give a uniform definition of a hyperplane in a space of constant curvature as a totally geodesic hypersurface and define a reflection in such a space as an isometric involution whose set of fixed points is a hyperplane. 

Let $H$ be a hyperplane in $X^n$. Its complement $X^n\setminus H$ consists of two connected components. The closure of a component  is  called a \emph{halfspace}. A reflection $r_H$ permutes the two halfspaces. One can distinguish the two halfspaces by a choice of one of the two unit vectors in $V$ orthogonal to the corresponding linear hyperplane $\tilde{H}$. We choose it so that it belongs to the corresponding halfspace. For any vector $v$ perpendicular to a hyperplane in $H^n$ we have $(v,v) > 0$ (otherwise the intersection of the hyperplane with $H^n$ is empty). 

Let $H_1^+,H_2^+$ be two halfspaces, and $e_1,e_2$ be the corresponding unit vectors. If $X\ne H^n$, the \emph{angle} $\phi = \angle(H_1^+,H_2^+) = \angle (H_1^-,H_2^-)$ is defined by 
$$\cos\phi = -(e_1,e_2), \quad 0\le \phi\le \pi.$$
If $X^n\ = H^n$ we use the same definition if $(e_1,e_2) \le 1$, otherwise we say that the angle is \emph{divergent} (in this case $(e_1,e_2) $ is equal to the hyperbolic cosine of the distance between the hyperplanes).

\begin{figure}[h]
\begin{multicols}{2}{
\xy
(-35,0)*{};(0,0)*\cir<40pt>{}; (-10,-10)*{};(0,25)*{}**\dir{-};
(10,-10)*{};(-5,25)*{}**\dir{-};
(0,-20)*{\text{divergent}};
\endxy

\xy
(0,25)*{};
(-15,0)*{};(0,0)*\cir<40pt>{}; (-8,-16)*{};(8,16)*{}**\dir{-};
(-8,16)*{};(8,-16)*{}**\dir{-};
(0,-20)*{\text{angle}};

\endxy}

\end{multicols}
\caption{}
\end{figure}

A \emph{convex polytope} in $X^n$ is a non-empty intersection of a locally finite\footnote{Locally finite means that each compact subset of $X^n$ is intersected by only finitely many hyperplanes.} set of halfspaces
$$P = \cap_{i\in I}H_i^-.
$$ 
The normal vectors $e_i$ defining the halfspaces $H_i^-$ all look inside the polytope. The hyperplanes $H_i$'s are called \emph{faces} of the polytope. 
In the case $X^n = H^n$ we also add to $P$ the points of the intersection lying on the absolute. We will always assume that no $e_i$ is a positive linear combination of others, or, equivalently, none of the halfspaces contains the intersection of the rest of half-spaces. In this case the set of bounding hyperplanes can be reconstructed  from $P$.  

A convex polytope has a finite volume if and only if  it is equal to a convex hull of finitely many points (vertices) from $X^n$ (or from the absolute if $X^n = H^n$). Such a polytope has finitely many faces. If $X^n = E^n$ or $S^n$, it is a compact polytope. A  polytope of finite volume in $H^n$ is compact only  if its vertices do not lie on the absolute.

\subsection{Reflection groups}
A \emph{reflection group} in a space of constant curvature is a discrete group of motions of $X^n$ generated by reflections.

\begin{theorem} Let $\Gamma$ be a reflection group in $X^n$.  There exists a convex polytope 
$P(\Gamma) = \cap_{i\in I}H_i^-$ such that
\begin{itemize} 
\item[(i)] 
$P  $ is a fundamental domain for the action of $\Gamma$ in $X^n$;
\item[(ii)] the angle between any two halfspaces  $H_i^-,H_j^-$ is equal to zero or $\pi/m_{ij}$ for some positive integer $m_{ij}$  unless the angle is divergent;
\item[(iii)] $\Gamma$ is generated by reflections $r_{H_i}, i\in I$.
\end{itemize}
Conversely, for every convex polytope $P$ satisfying property (ii)  the group $\Gamma(P)$ generated by the reflections into its facets is a reflection group and $P$ satisfies (i). 
\end{theorem}
 
\begin{proof} Consider the set $\calH$ of mirror hyperplanes of all reflections contained in $\Gamma$. For any mirror hyperplane $H$ and $g\in \Gamma$, the hyperplane $g(H)$ is the mirror hyperplane for the reflection $gr_Hg^{-1}$. Thus the set $\calH$ is invariant with respect to $\Gamma$. Let $K$ be a compact subset of $X^n$. For any hyperplane $H\in \calH$ meeting $K$, we have $r_H(K)\cap K\ne \emptyset$. Since
$\Gamma$ is a discrete group, the set $\{g\in G:g(K)\cap K\ne \emptyset\}$ is finite. This shows that the set $\calH$ is locally finite. The closure of a  connected component of 
$$X^n\setminus \bigcup_{H\in \calH} H.$$
is a convex polytope, called a \emph{cell} (or \emph{$\Gamma$-cell}, or a \emph{fundamental polyhedron}) of $\Gamma$. Its faces are called \emph{walls}. Two cells which share a common wall are called adjacent. The corresponding reflection switches the adjacent cells. This easily shows that the group $\Gamma$ permutes transitively the cells. Also it shows that any hyperplane from $\calH$ is the image of a wall of $P$ under an element of the group $\Gamma(P)$ generated by the reflections with respect to walls of $P$. Thus $\Gamma = \Gamma(P)$.
It is clear that the orbit of each point intersects a fixed  cell $P$. The proof  that no two interior points of $P$ belong to the same orbit follows from the last assertion of the theorem. Its proof is rather complicated, and we omit it (see \cite{VS}, Chapter V, Theorem 1.2).

Let   $H,H'$ be two hyperplanes bounding $P$ for which the angle $\angle(H^-,H'{}^-)$ is defined and is not zero. The corresponding unit vectors $e,e'$ span a plane in the vector space $V$ associated to $X^n$ and the restriction of the symmetric bilinear form to the plane is positive definite. The subgroup of $\Gamma$ generated by the reflections $r_{H}, r_{H'}$ defines a reflection subgroup in $\Pi$. Thus the angle must be of the form $r\pi$ for some rational number $r$. If $r$ is not of the form $1/m$ for some integer $m$, then $\Gamma$ contains a reflection with respect to a hyperplane intersecting the interior of $P$. By definition of $P$ this is impossible. This proves (i)-(iii).

\end{proof}

Define a \emph{Coxeter polytope} to be a convex polytope $P$ in which  any two faces are either divergent or form the angle equal to zero or $\pi/m$ for some positive integer $m$ (or $\infty$). Let $(e_i)_{\in I}$ be the set of unit vectors corresponding to the halfspaces $H_i^-$ defining $P$. The matrix 
$$G(P) = ((e_i,e_j))_{(i,j)\in I\times I}$$
is  the \emph{Gram matrix} of $P$ and its rank is the \emph{rank} of the Coxeter polytope. The polytope $P$ is called \emph{irreducible} if its Gram matrix is not equal to the non-trivial direct sum of matrices.

One can describe the matrix $G(P)$ via a certain labeled  graph, the \emph{Coxeter diagram} of $P$. Its vertices correspond to the walls of $P$. Two vertices corresponding to the hyperplanes with angle of the form $\pi/m, m\ge 3,$ are joined with an edge labeled with the number $m-2$ (dropped if $m = 3$) or joined with $m-2$ non-labeled edges. Two vertices corresponding to parallel hyperplanes (i.e. forming the zero angle) are joined by a thick edge, or is  labeled with $\infty$. Two vertices corresponding to divergent hyperplanes are joined by a dotted edge. 

Obviously an irreducible polytope is characterized by the condition that its Coxeter diagram is  a connected graph.

Let $\Gamma$ be a reflection group in $X^n$ and $P$ be  a $\Gamma$-cell. We apply the previous terminology concerning $P$ to $\Gamma$. Since  $\Gamma$-cells are transitively permuted by $\Gamma$, the isomorphism type of the Gram matrix does not depend on a choice of $P$. In particular, we can speak about \emph{irreducible reflection groups}. They correspond to Gram matrices which cannot be written as a non-trivial direct sum of its submatrices.

Suppose all elements of $\Gamma$ fix a point $x_0$ in $X^n$ (or on the absolute of $H^n$). Then all mirror hyperplanes  of reflections in $\Gamma$ contain $x_0$. Therefore each $\Gamma$-cell is a polyhedral cone with vertex at $x_0$. In the case when $X^n = E^n$, we can use $x_0$ to identify $E^n$ with its linear space $V$ and the group $\Gamma$ with a reflection group in $S^{n-1}$. The same is true if $X^n = S^n$. If $X^n = H^n$ and $x_0\in H^n$ (resp. $x_0$ lies on the absolute), then by considering the orthogonal subspace in $\bbR^{n+1}$ to the line defined by $x_0$ we find an isomorphism from $\Gamma$ to a reflection group in $S^{n-1}$ (resp. $E^{n-1}$).  

Any convex polytope of finite volume in $E^n$ or $H^n$ is nondegenerate in the sense that its faces do not have a common point and the unit norm vectors of the faces span the vector space $V$. A spherical convex polytope is nondegenerate if it does not contain opposite vertices.

\begin{theorem} Let $\Gamma$ be an irreducible reflection  group in $X^n$ with nondegenerate $\Gamma$-cell $P$. If $X^n = S^n$, then $\Gamma$ is finite and $P$ is equal to the intersection of $S^n$ with a simplicial cone in $\bbR^{n+1}$. If $X^n = E^n$, then $\Gamma$ is infinite and $P$ is a simplex in  $E^n$.
\end{theorem}

This follows from the following simple lemma (\cite{Bou}, Chapter V, \S3, Lemma 5):

\begin{lemma} Let $V$ be a real vector space with positive definite  symmetric bilinear form $(v,w)$ and 
$(v_i)_{i\in I}$ be vectors in $V$ with $(v_i,v_j) \le 0$ for $i\ne j$. Assume that the set $I$ cannot be nontrivially split into the union of two subsets $I_1$ and $I_2$ such that $(v_i,v_j) = 0$ for $i\in I_1, j\in I_2$. Then either the vectors $v_i$ are linearly independent or span a hyperplane and a  linear dependence can be chosen of the form $\sum a_iv_i$ with all $a_i$ positive.
\end{lemma}

The classification of irreducible non-degenerate Coxeter polytopes, and hence irreducible reflection groups $\Gamma$ in $S^n$ and $\bbR^n$ with nondegenerate  $\Gamma$-cell was given by Coxeter \cite{Cox}. The corresponding list of Coxeter diagrams is given in  Table \ref{spherical}.
\begin{table}[ht]
\caption{Spherical and euclidean real reflection groups}\label{spherical}
\xy (-30,0)*{};
(-10,0)*{A_n};@={(0,0),(5,0),(15,0),(20,0)}@@{*{\bullet}};
(0,0)*{};(7,0)*{}**\dir{-};(13,0)*{};(20,0)*{}**\dir{-};
(10,0)*{\ldots};
(40,0)*{\tilde{A}_n};
@={(50,0),(55,0),(65,0),(70,0),(60,-3)}@@{*{\bullet}};
(50,0)*{};(57,0)*{}**\dir{-};(63,0)*{};(70,0)*{}**\dir{-};(60,0)*{\ldots};(50,0)*{};(60,-3)*{}**\dir{-};
(70,0)*{};(60,-3)*{}**\dir{-};(60,0)*{\ldots};
(50,5)*{\bullet};(57,5)*{\bullet}**\dir{-};
(53,7.5)*{\infty};(40,5)*{\tilde{A}_1};
(-10,-10)*{B_n,C_n};@={(0,-10),(5,-10),(15,-10),(20,-10)}@@{*{\bullet}};
(0,-10)*{};(7,-10)*{}**\dir{-};(12,-10)*{};(15,-10)*{}**\dir{-};(15,-10)*{};(20,-10)*{}**\dir2{-};(10,-10)*{\ldots};(40,-10)*{\tilde{B}_n};
@={(50,-13),(50,-7),(55,-10),(60,-10),(70,-10),(75,-10),(80,-10)}@@{*{\bullet}};
(55,-10)*{};(62,-10)*{}**\dir{-};(68,-10)*{};(75,-10)*{}**\dir{-};(75,-10)*{};(80,-10)*{}**\dir2{-};(65,-10)*{\ldots};(50,-13)*{};(55,-10)*{}**\dir{-};
(50,-7)*{};(55,-10)*{}**\dir{-};
(40,-20)*{\tilde{C}_n};
@={(50,-20),(55,-20),(60,-20),(70,-20),(75,-20),(80,-20)}@@{*{\bullet}};
(50,-20)*{};(55,-20)*{}**\dir2{-};(55,-20)*{};(62,-20)*{}**\dir{-};(68,-20)*{};(75,-20)*{}**\dir{-};
(75,-20)*{},(80,-20)*{}**\dir2{-};(65,-20)*{\ldots};
(-10,-30)*{D_n};@={(0,-30),(5,-30),(15,-30),(20,-27),(20,-33)}@@{*{\bullet}};
(0,-30)*{};(7,-30)*{}**\dir{-};(12,-30)*{};(15,-30)*{}**\dir{-};(15,-30)*{};(20,-33)*{}**\dir{-};(10,-30)*{\ldots};
(15,-30)*{};(20,-27)*{}**\dir{-};
(40,-30)*{\tilde{D}_n};
@={(50,-33),(50,-27),(55,-30),(65,-30),(70,-33),(70,-27)}@@{*{\bullet}};
(55,-30)*{};(57,-30)*{}**\dir{-};(62,-30)*{};(65,-30)*{}**\dir{-};(65,-30)*{};(70,-33)*{}**\dir{-};(65,-30)*{};(70,-27)*{}**\dir{-};(60,-30)*{\ldots};(50,-27)*{};(55,-30)*{}**\dir{-};
(50,-33)*{};(55,-30)*{}**\dir{-};(60,-30)*{\ldots};
(-10,-40)*{E_6};@={(0,-40),(5,-40),(10,-40),(15,-40),(20,-40)}@@{*{\bullet}};
(0,-40)*{};(20,-40)*{}**\dir{-};(10,-40)*{};(10,-45)*{\bullet}**\dir{-};
(40,-40)*{\tilde{E}_6};
@={(50,-40),(55,-40),(60,-40),(60,-50),(65,-40),(70,-40)}@@{*{\bullet}};
(50,-40)*{};(70,-40)*{}**\dir{-};(60,-40)*{};(60,-45)*{\bullet}**\dir{-};(60,-40)*{};(60,-50)*{\bullet}**\dir{-};
(-10,-55)*{E_7};@={(0,-55),(5,-55),(10,-55),(15,-55),(20,-55),(25,-55),(10,-60)}@@{*{\bullet}};
(0,-55)*{};(25,-55)*{}**\dir{-};(10,-55)*{};(10,-60)*{\bullet}**\dir{-};
(40,-55)*{\tilde{E}_7};
@={(50,-55),(55,-55),(60,-55),(65,-55),(70,-55),(75,-55),(80,-55)}@@{*{\bullet}};
(50,-55)*{};(80,-55)*{}**\dir{-};(65,-55)*{};(65,-60)*{\bullet}**\dir{-};
(-10,-65)*{E_8};@={(0,-65),(5,-65),(10,-65),(15,-65),(20,-65),(25,-65),(30,-65),(10,-70)}@@{*{\bullet}};
(0,-65)*{};(30,-65)*{}**\dir{-};(10,-65)*{};(10,-70)*{\bullet}**\dir{-};
(40,-65)*{\tilde{E}_8};
@={(50,-65),(55,-65),(60,-65),(65,-65),(70,-65),(75,-65),(80,-65),(85,-65)}@@{*{\bullet}};
(50,-65)*{};(85,-65)*{}**\dir{-};(60,-65)*{};(60,-70)*{\bullet}**\dir{-};
(-10,-75)*{F_4};@={(0,-75),(5,-75),(10,-75),(15,-75)}@@{*{\bullet}};
(0,-75)*{};(5,-75)*{}**\dir{-};(5,-75)*{};(10,-75)*{}**\dir2{-};(10,-75)*{};(15,-75)*{}**\dir{-};
(40,-75)*{\tilde{F}_4};@={(50,-75),(55,-75),(60,-75),(65,-75),(70,-75)}@@{*{\bullet}};
(50,-75)*{};(55,-75)*{}**\dir{-};(55,-75)*{};(60,-75)*{}**\dir2{-};(60,-75)*{};(70,-75)*{}**\dir{-};
(-10,-85)*{I_2(m)};(0,-85)*{\bullet};(8,-85)*{\bullet}**\dir{-};(20,-85)*{m\ge 5};(4,-83)*{m};
(40,-85)*{\tilde{G}_2};@={(50,-85),(55,-85),(60,-85)}@@{*{\bullet}};(50,-85)*{}**\dir{-};(60,-85)*{}**\dir{-};(53,-82)*{6};
(-10,-95)*{H_3};@={(0,-95),(5,-95),(10,-95)}@@{*{\bullet}};
(0,-95)*{};(5,-95)*{}**\dir3{-};(5,-95)*{};(10,-95)*{}**\dir{-};
(-10,-105)*{H_4};@={(0,-105),(5,-105),(10,-105),(15,-105)}@@{*{\bullet}};(0,-105)*{};(5,-105)*{}**\dir3{-};(5,-105)*{};(10,-105)*{}**\dir{-};(10,-105)*{};(15,-105)*{}**\dir{-};
(-2,7)*{};(-2,-108)*{}**\dir{-};(33,7)*{};(33,-108)*{}**\dir2{-};(47,7)*{};(47,-108)*{}**\dir{-};
(10,-115)*{\text{Spherical groups}};(60,-115)*{\text{Euclidean groups}};(20,-115)*{};
\endxy
 \end{table}

 Here the number of nodes in the spherical (resp. euclidean) diagram is equal to the subscript  $n$ (resp. $n+1$) in the notation. The number $n$ is equal to the rank of the corresponding Coxeter polytope. We will refer to diagram from the first (resp. second) column as \emph{elliptic Coxeter diagrams}  (resp. \emph{parabolic Coxeter diagrams}) of rank $n$.

Our classification of plane reflection groups in section 1.1 fits in this classification
$$(2,4,4) \longleftrightarrow \tilde{C}_2,\ (2,3,6)  \longleftrightarrow \tilde{G}_2, \ 
(3,3,3)  \longleftrightarrow \tilde{A}_2, \  (2,2,2,2)  \longleftrightarrow \tilde{A}_1\times \tilde{A}_1. $$

The list of finite reflection groups not of type $H_3,H_4, I_2(m), m \ne 6$ ($I_2(6)$ is often denoted by $G_2$) coincides with the list of \emph{Weyl groups} of simple Lie algebras of the corresponding type $A_n, B_n$ or $C_n$, $G_2, F_4,E_6,E_7,E_8$.  The corresponding Coxeter diagrams coincide with the Dynkin diagrams only in the cases $A, D,$ and $E$. The second column corresponds to \emph{affine Weyl groups}.

The group of type $A_n$ is the symmetric group $\Sigma_{n+1}$. It acts   in the space
\begin{equation}\label{space}
V = \{(a_1,\ldots,a_{n+1})\in \bbR^{n+1}:a_1+\ldots a_{n+1} = 0\}\end{equation}
 with the standard inner product as the group generated by reflections  in vectors $e_i-e_{i+1}, i = 1,\ldots n$. 
 
 The group of type $B_n$ is isomorphic to the semi-direct product $2^n\rtimes \Sigma_n$. It acts in the euclidean space $\bbR^n$ as a group generated by reflections in vectors $e_i-e_{i+1}, i = 1,\ldots,n-1,$ and $e_n$. 
 
 The group of type $D_n$ is isomorphic to the semi-direct product $2^{n-1}\rtimes \Sigma_n$. It acts in the euclidean space $\bbR^n$ as a group generated by reflections in vectors  $e_i-e_{i+1}, i= 1,\ldots, n-1,$ and $e_{n-1}+e_n$.

A discrete group $\Gamma$ of motions in $X^n$ admitting  a fundamental domain of finite volume (resp. compact) is said to be \emph{of finite covolume} (resp. \emph{cocompact}). \footnote{If $\Gamma$ is realized as a discrete subgroup of a Lie group that acts properly and transitively on $X^n$, then this terminology agrees with the terminology of discrete subgroups of a Lie group.} Obviously, a simplex in $S^n$ or $E^n$ is compact. Thus  the previous list gives a classification of irreducible  reflection groups of finite covolume  in $E^n$ and $S^n$. They are are automatically cocompact.

The classification of  reflection groups of finite covolume in $H^n$ is known only for $n = 2$ (H. Poincar\'e) and $n = 3$ (E. Andreev \cite{An}).   It is known that they do not exist if $n\ge 996$ (\cite{Kh1}, \cite{Kh2}, \cite{Pr}) and even if $n > 300$ (see the announcement in \cite{Ni1}). There are no  cocompact reflection groups  in $H^n$ for $n\ge 30$ \cite{Vi6}.

One can give the following description of Coxeter diagrams defining  reflection groups of cofinite volume(see \cite{Vi1}).

\begin{proposition}\label{bugaenko} A reflection group $\Gamma$ in $H^n$ is of finite covolume if and only if any elliptic subdiagram of rank $n-1$  of its Coxeter diagram can be extended in exactly two ways to an elliptic subdiagram of rank $n$ or a parabolic subdiagram of rank $n-1$. Moreover, $\Gamma$ is cocompact  if the same is true but there are no parabolic subdiagrams of rank $n-1$. 
\end{proposition}

The geometric content of this proposition is as follows. The intersection of hyperplanes defining an elliptic subdiagram of rank $n-1$ is of  dimension 1 (a one-dimensional facet of the polytope). An elliptic subdiagram (resp. parabolic) of rank $n$ defines a proper (resp. improper) vertex of the polytope. So, the proposition says that $\Gamma$ is of finite covolume if and only if each one-dimensional facet joins precisely two vertices, proper or improper.

\subsection{Coxeter groups} Recall that a \emph{Coxeter group} is a group $W$ admitting an ordered  set of generators $S$ of order 2 with defining relations 
$$ (ss')^{m(s,s')} = 1, \ s,s'\in S,$$
where $m(s,s')$ is the order of the product $ss'$ (the symbol $\infty$ if the order is infinite). The pair $(W,S)$ is called a \emph{Coxeter system}. 

 The \emph{Coxeter graph} of $(W,S)$ is the graph whose vertices correspond to $S$ and two different vertices are connected by $m(s,s')-2$ edges or by an edge  labeled with $m(s,s')-2$ or  a thick edge  if $m(s,s') = \infty$.  We say that $(W,S)$ is \emph{irreducible} if the Coxeter graph is connected.

One proves the the following theorem (see \cite{Vi1}).

\begin{theorem} Let $P$ be a nondegenerate Coxeter polytope of finite volume in $X^n$ and $\Gamma(P)$ be the corresponding reflection group. The pair $(\Gamma(P),S)$, where $S$ is the set of reflections with respect to the set of faces of $P$ is a Coxeter system. Its Coxeter graph is equal to the Coxeter diagram of $P$ and the Gram matrix of $P$ is equal to the  matrix
\begin{equation}\label{cosine}
(-\cos \frac{\pi}{m(s,s')})_{(s,s')\in S\times S}.
\end{equation}
\end{theorem}

The converse is partially true. The following facts can be found in \cite{Bou}. Let $(W,S)$ be an irreducible  Coxeter system with no $m(s,s')$ equal to  $ \infty$. One considers the linear space  $V = \bbR^{S}$ and equips it with a symmetric bilinear form $B$ defined by 
\begin{equation}\label{cos}
B(e_s,e_{s'}) = -\cos \frac{\pi}{m(s,s')}.
\end{equation}

Assume that $B$ is positive definite ($(W,S)$ is \emph{elliptic}). Then $W$ is finite and isomorphic to a reflection group $\Gamma$ in the spherical space $S^n$, where $n+1 = \#S-1$. The corresponding $\Gamma$-cell can be taken  as the intersection of the sphere with the simplex in $\bbR^{n+1}$ with facets orthogonal to the vectors $e_s$. 

Assume that $B$ is degenerate and semipositive definite ($(W,S)$ is   \emph{parabolic}). Then its radical $V_0$ is one-dimensional and is spanned by a unique vector $v_0 = \sum a_se_s$ satisfying  $a_s > 0$ for all $s$ and $\sum a_s = 1$. The group $W$ acts naturally as a reflection group $\Gamma$ in the affine subspace $E^{n} = \{\phi\in V^*:\phi(v_0) = 1\}$ with the associated  the linear space 
$(V/V_0)^*$. The corresponding $\Gamma$-cell has $n+1$ facets orthogonal to the vectors $e_s$ and is a simplex in affine space.

Assume that $B$ is non-degenerate, indefinite and $W$ is of cofinite volume in the orthogonal group of  $B$ ($(W,S)$ is  \emph{hyperbolic}). In this case the signature of $B$ is equal to $(n-1,1)$ and  $C = \sum_{s\in S} \bbR_+e_s$ is contained in one of the two connected components of the set $\{x\in E: B(x,x) < 0\}$. Let $H^n$ be the hyperbolic space equal to the image of this component in $\bbP(E)$. Then the action of $W$ in $H^n$ is isomorphic to a reflection group $\Gamma$. A $\Gamma$-cell can be chosen to be the image of the closure $\bar{C}$ of $C$ in $\bbP(E)$. 

In all cases the Coxeter diagram of a $\Gamma$-cell coincides with the Coxeter graph of $(W,S)$. We will call the matrix given by \eqref{cos}, the Gram matrix of $(W,S)$.

\begin{remark} Irreducible reflection groups in $S^n$ and $E^n$ correspond to elliptic or parabolic irreducible Coxeter systems.  Hyperbolic Coxeter systems   define Coxeter simplices in $H^n$ of finite volume. Their Coxeter diagrams are called \emph{quasi-Lanner} \cite{VS} or \emph{hyperbolic} \cite{Bou}. They are characterized by the condition that each its proper subdiagram is either  elliptic  or parabolic. If no parabolic subdiagram is present then the simplex is compact and the diagram is \emph{Lanner} or  \emph{compact hyperbolic}. The complete  list of  hyperbolic Coxeter diagrams  can be found in \cite{VS} or  \cite{Hu}.
\end{remark}

\begin{example}\label{triangle} Let $W(p,q,r), 1\le p\le q\le r,$ be the Coxeter group with Coxeter graph of type $T_{p,q,r}$ given in Figure \ref{pqr}.

\begin{figure}[ht]
\xy
(-35,0)*{};
@={(0,0),(7,0),(20,0),(27,0),(34,0),(47,0),(54,0),(27,-7),(27,-20),(27,-27)}@@{*{\bullet}};
(0,0)*{};(10,0)*{}**\dir{-};(17,0)*{};(37,0)*{}**\dir{-};(44,0)*{};(54,0)*{}**\dir{-};
(27,0)*{};(27,-10)*{}**\dir{-};(27,-17)*{};(27,-27)*{}**\dir{-};(40.5,0)*{\ldots};(13.5,0)*{\ldots};(27,-12.5)*{\vdots}; 
(13.5,-12)*+++++++++++++{};(18,-12)*{}
**\frm{^\}};(14,5)*{\text{q}};
(40,-11)*++++++++++++{};(50,-11)*{}
**\frm{^\}};(40,5)*{\text{r}};
(37.5,-13)*++++++++++++{};(37.5,-18)*{}
**\frm{\{};(22,-13)*{\text{p}};
\endxy

\caption{}\label{pqr}
\end{figure}
Then 
\begin{itemize} 
\item $W(p,q,r)$ is of elliptic type if and only if  
$$(p,q,r) = (1,q,r), (2,2,r), (2,3,3), (2,3,4), (2,3,5);$$
 \item $W(p,q,r)$ is of parabolic  type if and only if 
 $$(p,q,r) =  (2,4,4), (2,3,6), (3,3,3);$$
 \item $W(p,q,r)$ is of hyperbolic type if and only if 
 $$(p,q,r) =  (3,4,4), (2,4,5), (2,3,7).$$
\end{itemize}

\end{example}

\section{Linear reflection groups}

\subsection{Pseudo-reflections} Let $E$ be a vector space over any field $K$. A \emph{pseudo-reflection} in $E$ is a linear invertible transformation $s:E\to E$ of finite order  greater than 1 which fixes pointwise a hyperplane.   A \emph{reflection} is a diagonalizable pseudo-reflection. A pseudo-reflection is a reflection if and only if its order  is coprime to the characteristic of $K$.

Let $v$ be an eigenvector of a reflection $s$ of order $d$. Its eigenvalue $\eta$ different from 1 is a $d$th root of unity $\eta\ne 1$. We can write $s$ in the form
\begin{equation}\label{refl}
s(x) = x-\ell(x)v
\end{equation}
for some linear function $\ell:E\to K$. Its zeroes define the hyperplane of fixed points of $s$, the \emph{reflection hyperplane}. Taking $x=v$ we obtain $\ell(v) = 1-\eta $. This determines $\ell$ uniquely when $v$ is fixed; we denote it by $\ell_v$.

A \emph{pseudo-reflection (reflection) subgroup} of $\GL(E)$ is a subgroup generated by pseudo-reflections (reflections).

Assume that we are given  an automorphism $\sigma$ of $K$ whose square is the identity. We denote its value on an element $\lambda\in K$  by $\bar{\lambda}$. Let $B:E\times E\to K$ be a  $\sigma$-hermitian form on $E$, i.e. $B$ is $K$-linear in the first variable and satisfies $B(x,y) = \overline{B(y,x)}.$  

Let $\U(E,B,\sigma)$ be the unitary group of $B$, i.e. the subgroup of $K$-linear transformations of $E$ which preserve $B$.  A pseudo-reflection subgroup in $\GL(E)$ is  a \emph{unitary pseudo-reflection group} if it is contained in a unitary group  $\U(E,B,\sigma)$ for some $\sigma,B$ and for any reflection \eqref{refl} one can choose a vector $v$  with $B(v,v) \ne 0$. 

The additional  condition implies that
\begin{equation}
\ell_v(x) = \frac{(1-\eta)B(x,v)}{B(v,v)}\ \text{for all}\ x\in V.
\end{equation}
In particular, the vector $v$ is orthogonal to the hyperplane $\ell_v^{-1}(0)$.

Finite reflection groups are characterized by the following property of its algebra of invariant polynomials (\cite{Bou},  Chapter V, \S 5, Theorem 4).

\begin{theorem}\label{chevalley} A finite subgroup $G$ of $\GL(E)$ of order prime to $\cha(K)$ is a reflection group if and only if the algebra $S(E)^G$ of invariants in the symmetric algebra of $E$ is isomorphic to a polynomial algebra.
\end{theorem}

In the case $K=\bbC$ this theorem was proven by Shephard and Todd \cite{ST}, and in the case of arbitrary characteristic, but for groups generated by reflections of order 2, it  was proven by C. Chevalley \cite{Ch}.

We will be concerned with the case where $K = \bbR$ and  $\sigma = \id_E$ or $K = \bbC$ and $\sigma$ is the complex conjugation. In the real case a reflection is necessarily of order 2. 

Let $G$ be a finite reflection  subgroup of $\GL(E)$. By taking some positive definite symmetric bilinear (resp. hermitian) form and averaging it, we see that $G$ is conjugate in $\GL(E)$ to a unitary reflection group. In fact, under isomorphism from $E$ to the standard euclidean (resp.unitary) space $\bbR^n$ (resp. $\bbC^n$), the group $G$ is isomorphic to a reflection subgroup of $\Or(n)$ (resp. $\U(n)$). 
In the real case $G$ becomes isomorphic to a reflection group in $S^{n-1}$, and hence is isomorphic to the product of irreducible reflection spherical groups.

An example of an infinite  real reflection group in $\bbR^{n+1}$ is a reflection group in the hyperbolic space $H^n$. It is an orthogonal reflection group with respect to a symmetric bilinear form of signature $(n,1)$.

\subsection{Finite complex linear reflection groups}
They were classified by Shephard and Todd \cite{ST}. Table \ref{complex} gives the list of irreducible finite linear reflection groups (in the order given by Shephard-Todd, see also the table in \cite{Kane}, p.166).

\begin{table}[ht]
\caption{Finite complex reflection groups}\label{complex}
\begin{center}
\begin{tabular}{|| l | r | r | r  |r||} \hline
Number&Name&Order&$\dim E$& Degrees\\ \hline\hline
1&$A_n = \Sigma_{n+1}$&$(n+1)!$&$n$&$2,3,\ldots,n+1$\\\hline
2&$G(m,p,n)$&$m^{n}n!/p$&$n$&$m,m+1,\ldots,(n-1)m,mn/p$\\ \hline
3&$[]^m$&$m$&1&$m$\\ \hline
4&3[3]3&24&2&4,6\\ \hline
5&3[4]3&72&2&6,12\\ \hline
6&3[6]2&48&2&4,12\\ \hline
7&$<3,3,3>_2$&144&2&12,12\\ \hline
8&4[3]4&96&2&8,12\\ \hline
9&4[6]2&192&2&8,24\\ \hline
10&4[4]3&288&2&12,24\\ \hline
11&$<4,3,2>_{12}$&576&2&24,24\\ \hline
12&$\GL(2,3)$&48&2&6,8\\ \hline
13&$<4,3,2>_2$&96&2&8,12\\ \hline
14&3[8]2&144&2&6,24\\ \hline
15&$<4,3,2>_6$&288&2&12,24\\ \hline
16&5[3]5&600&2&20,30\\ \hline
17&5]6]2&1200&2&20,60\\ \hline
18&5[4]3&1800&2&60,60\\ \hline
19&$<5,3,2>_{30}$&3600&2&60,60\\ \hline
20&3[5]3&360&2&12,30\\ \hline
21&3[10]2&720&2&12,60\\ \hline
22& $<5,3,2>_{2}$&240&2&12,20\\ \hline
23&$H_3$&120&3&2,6,10\\ \hline
24&$J_3(4)$&336&3&4,6,14\\ \hline
25&$L_3$&648&3&6,9,12\\ \hline
26&$M_3$&1296&3&6,12,18\\ \hline
27&$J_3(5)$&2160&3&6,12,30\\ \hline
28&$F_4$&1152&4&2,6,8,12\\ \hline
29&$N_4$&7680&4&4,8,12,20\\ \hline
30&$H_4$&14,440&4&2,12,20,30\\ \hline
31&$EN_4$&$64\cdot 6!$&4&8,12,20,24\\ \hline
32&$L_4$&$216\cdot 6!$&4&12,18,24,30\\ \hline
33&$K_5$&$72\cdot 6!$&5&4,6,10,12,18\\ \hline
34&$K_6$&$108\cdot 9!$&5&4,6,10,12,18\\ \hline
35&$E_6$&$72\cdot 6!$&6&2,5,6,8,9,12\\ \hline
36&$E_7$&$8\cdot 9!$&7&2,6,8,10,12,14,18\\ \hline
37&$E_8$&$192\cdot 9!$&8&2,8,12,14,18,20,24,30\\ \hline
\end{tabular}
\vspace{10pt}
\end{center}
\end{table}
The last column in the table gives the degrees of the generators of the algebra $S(E)^G$.

The group  $G(m,p,n)$ is equal to the semi-direct product $A(m,n,p)\rtimes \Sigma_n$, where $A(m,n,p)$ is a diagonal group of $n\times n$-matrices with $m$th roots of unity at the diagonal whose product is a $(m/p)$th root of unity. The semi-direct product is defined with respect to the action of $\Sigma_n$ by permuting the columns of the matrices. 

Groups No 3-22 are some extensions of binary polyhedral groups (i.e. finite subgroups of $\SL(2,\bbC)$). 

All real spherical irreducible groups are in the list. We have seen already the groups of types $A_n,E_6,E_7,E_8,H_3,H_4,F_4$. The groups of type $B_n$ are the groups $G(2,1,n)$. The groups of type $D_n$ are the groups $G(2,2,n)$. Finally, the groups of type $I_2(m)$ are the groups $G(m,m,2)$. 
These groups are distinguished from other groups by the property that one of invariant polynomials is of degree 2.

\subsection{Complex crystallographic  reflection groups} 
A complex analog of a space of constant curvature is a simply connected  complex K\"ahler manifold of constant holomorphic curvature (\emph{complex space form}).  There are three types of such spaces (see \cite{Igusa})
\begin{itemize}
\item $E_{\bbC}^n,$ the $n$-dimensional affine space equipped with the standard hermitian form
$|z|^2 = \sum |z_i|^2$. It is a homogeneous space $(\bbC^n\rtimes \U(n))/\U(n)$;
\item $\bbP^n(\bbC)$, the $n$-dimensional complex projective space equipped with the standard Fubini-Study metric. It  is a homogeneous space $\PU(n+1)/\U(n)$.
\item $\bbCH = \{z\in \bbC^{n}:|z| < 1\}$, the $n$-dimensional \emph{complex  hyperbolic space}. The hermitian metric on $H_\bbC^n$ is defined by 
$$\frac{1}{1-|z|^2}\bigl(\sum_{i=1}^n  z_id\bar{z}_i+\bar{z_i}dz_i)+(1-|z|^2)\sum_{i=1}^n dz_id\bar{z}_i\bigr).$$
\end{itemize}
They are  simply connected hermitian complex homogeneous manifolds of dimension $n$ with isotropy  subgroups equal to the unitary group $\U(n)$. 

The complex hyperbolic space has a model in complex projective space $\bbP^n(\bbC)$ equal to the image of
 the subset 
 $$C = \{(z_0,z_1,\ldots,z_n)\in \bbC^{n+1}:-|z_0|^2+|z_1|^2+\ldots+|z_n|^2 < 0\}.$$
The unitary group $\U(n+1)$ of the hermitian form $-|z_0|^2+|z_1|^2+\ldots+|z_n|^2$ of signature $(n,1)$ acts transitively on $\bbCH$ with isotropy subgroup $\U(n)$. It defines a transitive action of $\PU(n,1) = \U(n,1)/\U(1)$ with isotropy subgroups isomorphic to $\U(n)$. 

\bigskip
Let $X_\bbC^n$ be a $n$-dimensional complex space space forms. A reflection in $X_\bbC^n$ is a holomorphic isometry whose set of fixed points is a hypersurface. A reflection group is a discrete group of holomorphic automorphisms generated by reflections. A reflection group  $\Gamma$  of  $X_\bbC^n = \bbP^n(\bbC)$ (resp. $\bbCH$) can be centrally extended to a reflection  subgroup of $\U(n)$ (resp. $\U(n,1)$) such that over every reflection lies a linear reflection.

A reflection group of $E_{\bbC}^n$ is  a discrete subgroup of $\bbC^n\rtimes \U(n)$ which is generated by affine reflections.
It can be considered as a linear reflection group in a complex vector space $V$ of dimension $n+1$ equipped with a hermitian form of Sylvester signature $(t_+,t_-,t_0) = (n,0,1)$. We take for $E_\bbC^n$ the affine subspace of the dual linear space $V^*$ of linear functions $\phi\in V^*$ satisfying $\phi(v) = 1$, where $v$ is a fixed nonzero vector in $V^\perp$. The corresponding linear space is the hyperplane $\{\phi\in V^*:\phi(v) = 0\}$.

A reflection group $\Gamma$ of cofinite volume in $E_{\bbC}^n$ or $H_\bbC^n$ is called a \emph{complex crystallographic group} (affine, hyperbolic). If $X_\bbC^n = \bbP^n(\bbC)$ it is a finite group defines by a finite linear  complex reflection group in $\bbC^{n+1}$. If $X_\bbC^n = E_\bbC^n$, then $\Gamma$ leaves invariant  a lattice $\Lambda \subset \bbC^n$  of rank $2n$ (so that $\bbC^n/\Lambda$ is a compact complex $n$-torus), and $\Lin(\Gamma) $ is a finite subgroup of $\U(n)$. This implies that $\Gamma$ is also cocompact. 

If $X_\bbC^n = \bbP^n(\bbC)$, then $\Gamma$, being a discrete subgroup of a compact Lie group $\PU(n+1)$, is finite and  cocompact.


There is a classification of crystallographic reflection groups in $E_{\bbC}^n$ (due to V. Popov \cite{Po}\footnote{According to \cite{Go4} some groups are missing  in Popov's list}). 

First observe that any $a\in E_{\bbC}^n$ defines a surjective homomorphism $g\mapsto \bar{g}$ from the affine group to the linear group of the corresponding complex linear space $V$. We write any $x\in E_{\bbC}^n$ in the form
$x = a+v,$ for a unique $v\in V$, and get
$g(a+v) = g(a) +\bar{g}(v)$. This definition of $\bar{g}$ does not depend on the choice of $a$.  In particular, choosing $a$ on a reflecting affine hyperplane $H$, we see that $\bar{g}$ is a linear reflection which fixes $H-a$. This implies that the image of a crystallographic reflection subgroup $\Gamma$ of $\bbC^n\rtimes \U(n)$ is a finite reflection subgroup of $\U(n)$.

\begin{theorem} Let $G$ be a finite irreducible reflection group in $\U(n)$. Then the following properties are equivalent
\begin{itemize}
\item[(i)] there exists a complex reflection group $\Gamma$ in $E_{\bbC}^n$ with linear part $G$;
\item [(ii)]  there exists a $G$-invariant lattice $\Lambda\subset E_{\bbC}^n$ of rank $2n$;
\item [(iii)]  the number of the group $G$ in Table \ref{complex} is
$1, 2 (m=2,3,4,6),3 (m = 2,3,4,6), \\
4, 5, 8,12, 24-29, 31-37.$
\end{itemize}
If $G$ is not of type $G(4,2,n), n \ge 4$ (number 2), or $\GL(2,3)$ (type 12) or $EN_4$ (number 31)\footnote{also $G(6,6,n)$ as pointed out in \cite{Go4}.}, then 
$\Gamma$ is equal to the semi-direct product $\Lambda\rtimes G$. In the exceptional cases, $\Gamma$ is either the semi-direct product or some non-trivial extension of $G$ with normal subgroup $\Lambda$.
 \end{theorem}
 
A table in \cite{Po} describes all possible lattices and the extensions for each $G$ as above.

Recall from Theorem \ref{chevalley} that the algebra of invariant polynomials of a finite complex reflection group $\Gamma$ in $\bbC^n$ is a polynomial algebra. This can be restated as follows. One considers the induced action of $\Gamma$ in $\bbP^{n-1}(\bbC)$ and the orbit space $\bbP^{n-1}(\bbC)/\Gamma$ which exists as a projective algebraic variety. Now the theorem asserts that this variety is isomorphic to a weighted projective space $\bbP(q_1,\ldots,q_{n})$\footnote{This the quotient of $\bbC^n\setminus \{0\}$ by the action of $\bbC^*$ defined in coordinates by $(z_1,\ldots,z_n) \to (\lambda^{q_1}z_1,\ldots,\lambda^{q_n}z_n).$}, where the weights are equal to the degrees of free generators of the invariant algebra. The following  is an analog of Theorem  \ref{chevalley} for affine complex crystallographic groups due to Bernstein-Shwarzman \cite{BS} and Looijenga \cite{Lo2}. 

\begin{theorem}\label{BS} Assume that the linear part of a complex crystallographic group $\Gamma$ is a complexification of a real finite reflection group $W$. Then the orbit space $E_\bbC^n/\Gamma$ exists as an algebraic variety and is  isomorphic to a weighted projective space $\bbP(q_0,\ldots,q_n)$, where the weights are  explicitly determined by $W$. 
\end{theorem}

It is conjectured that the same is true without additional assumption on the linear part.

\begin{example} Let $G$ be a finite complex reflection groups arising from the complexification of a real reflection group $G_r$ . Any such group is realized as the linear part of a complex crystallographic group $\Gamma$ in affine space and $\Gamma$ is the semi-direct product of $G$ and a $G$-invariant lattice. Suppose $G_r$ is of $ADE$ type. Let $e_1,\ldots,e_n$ be the norm vectors of the Coxeter polytope in $\bbR^{n}$.   For any $\tau = a+bi, b > 0,$ consider the lattice $\Lambda_\tau$  in $\bbC^n$ spanned by the vectors $e_i$ and $\tau e_i$. This is a  $G$-invariant lattice and every $G$-invariant lattice is obtained in this way.  Moreover $\Lambda_\tau = \Lambda_{\tau'}$ if and only if $\tau$ and $\tau'$ belong to the same orbit of the modular group 
$\PSL(2,\bbZ)$ which acts on the upper half-plane $\{z = a+bi\in \bbC:b > 0\}$ by the Moebius transformations $z\mapsto (az+b)/(cz+d)$. The linear part $G$ is a finite group of automorphisms of the compex torus $(\bbC/\Lambda_\tau)^n$ with the orbit space  $(\bbC/\Lambda_\tau)^n/G$ isomorphic to a weighted projective space. In the case when $G_r$ is the Weyl group of a simple simply connected Lie group $H$,  this quotient  is naturally isomorphic to the moduli space of principal $H$-bundles on the elliptic curve $\bbC/\bbZ+\tau \bbZ$ (see \cite{FM}).
\end{example}

\begin{remark} If $\Gamma$ is a real crystallographic group in affine space $E^n$, then its complexification is a complex non-crystallographic reflection group in $E_{\bbC}^n$. Every complex non-crystallographic reflection group in $E_{\bbC}^n$ is obtained in this way (see \cite{Po}, 2.2).
\end{remark}

We will discuss later a construction of complex crystallographic reflection groups in $H_\bbC^n$ for $n\le 9$. The largest known dimension $n$ for which such group exists is 13 (\cite{All}). It is believed that these groups occur only in finitely many dimensions.

\section{Quadratic lattices and their reflection groups}
\subsection{Integral structure} Let $\Gamma$ be an orthogonal  linear reflection group in a real vector space $V$ of dimension $n$ equipped with a nondegenerate symmetric bilinear  form of signature $(n,0)$ or $(n-1,1)$. We assume that the intersection of its  reflection hyperplanes is the origin.  We say that $\Gamma$ admits an \emph{integral structure} if it leaves invariant  a free abelian subgroup $M\subset V$ of rank $n$ generating $V$. In other words,  there exists a basis $(e_1,\ldots,e_n)$ in $V$ such that 
$\sum_{i=1}^n\bbZ e_i$ is $\Gamma$-invariant.

A linear reflection group admitting an integral structure  is obviously a discrete subgroup of the orthogonal group $O(V)$ and hence acts discretely on the corresponding space of constant curvature $S^{n-1}$ or $H^{n-1}$ (because the isotropy subgroups are compact subgroups of $O(V))$. Thus $\Gamma$ is a reflection group of $S^{n-1}$ or $H^{n-1}$. By a theorem of Siegel \cite{Si}, the group $\Or(M) =  \{g\in O(V):g(M) = M\}$ is of finite covolume in the orthogonal group $\Or(V)$, hence $\Gamma$ is of finite covolume  if and only it is of finite index in $\Or(M)$.

Let $H \subset V$ be a reflection hyperplane in $\Gamma$ and  let $e_H$ be an orthogonal vector to $H$. In the hyperbolic case we assume that $H$ defines a hyperplane in $H^{n-1}$, hence $(e_H,e_H) > 0$. The reflection $r_H$ is defined by 
$$r_H(x) = x-\frac{2(x,\alpha_H)}{(\alpha_H,\alpha_H)}e_H$$
for some vector $\alpha_H$ proportional to $e_H$. Taking $x$ from $M$ we obtain that the vector 
$\frac{2(x,\alpha_H)}{(\alpha_H,\alpha_H)}\alpha_H$ belongs to $ M.$
 Replacing $\alpha_H$ by proportional vector, we may assume that $\alpha_H\in M$ and also that $\alpha_H$ is a primitive vector in $M$ (i.e. $M/\bbZ \alpha_H$ is torsion-free). We call such a vector a \emph{root vector} associated to $H$. The root vectors corresponding to the faces of a fundamental polytope  are called the \emph{fundamental root vectors}. 
 
A root is uniquely defined by $H$ up to multiplication by $-1$.  Using the primitivity property of 
root vectors it is easy to see that, for all $x\in M$,
\begin{equation}\label{ind}
2(x,\alpha)\in (\alpha,\alpha)\bbZ.\end{equation}
 In particular, if $\alpha,\beta$ are not perpendicular root vectors, then $2(\alpha,\beta)/(\beta,\beta)$ and $2(\alpha,\beta)/(\alpha,\alpha)$ are nonzero integers, so that the ratio $(\beta,\beta)/(\alpha,\alpha)$ is a rational number. If the Coxeter diagram is connected we can fix one of the roots $\alpha$ and multiply the quadratic form $(x,x)$ on $V$ by 
$(\alpha,\alpha)^{-1}$ to assume that $(\beta,\beta)\in \bbQ$ for all root vectors. This implies that 
the Gram matrix $(e_i,e_j)$ of a basis of $M$ has entries in $\bbQ$. Multiplying the quadratic form by an integer we may assume that it is an integral matrix, hence 
\begin{equation}\label{int}
(x,y)\in \bbZ, \ \text{for all}\ x,y\in M.
\end{equation}

\begin{example}\label{rootlat} Let $\Gamma$ be an irreducible finite real reflection group whose Dynkin diagram does not contain multiple edges (i.e. of types $A,D,E$). Let $(\alpha_1,\ldots,\alpha_n)$ be the unit norm vectors of its fundamental Coxeter polytope and let $M \subset V$ be the span of these vectors. If we multiply the inner product in $V$ by 2, we obtain  $(\alpha_i,\alpha_j) = -2\cos \frac{\pi}{m_{ij}} \in \{0,2,-1\}$. Hence  the reflections
$r_{\alpha_i}:x\mapsto x-(x,\alpha_i)\alpha_i$ leave $M$ invariant.  Thus $\Gamma$ admits an integral structure and the Gram matrix of its basis $(\alpha_1,\ldots,\alpha_n)$ is equal to the twice of  matrix \eqref{cosine}.  Note that $(\alpha_i,\alpha_i) = 2, i = 1,\ldots,n$.

Let $\Gamma$ be of type $B_n$ and $e_s,s=1,\ldots,n,$ be the unit normal vectors defined by a fundamental polytope. We assume that $m(n,n-1) = m(n-1,n) = 4$ and other $m(s,s')$  take value in $\{1,2,3\}$.  Let $\alpha_i = e_i$ if $i\ne n$ and $\alpha_n = \sqrt{2}e_n$. It is easy to see now that $\frac{2(x,\alpha_i)}{(\alpha_i,\alpha_i)}\in \bbZ$ for any $x$ in the span $M$ of the $\alpha_i$'s. This shows that $M$ defines an integral structure on $\Gamma$ and $(x,y) \in \bbZ$ for any $x,y\in M$. We have 
\begin{equation}\label{odd}
(\alpha_i,\alpha_i) = \begin{cases}
      1& \text{if}\  i\ne n, \\
      2& \text{otherwise}.
\end{cases}
\end{equation}

We leave to the reader to check that the groups of type $F_4$ and $G_2 = _2(6)$ also admit an integral structure. However the remaining groups  do not.
\end{example}

\subsection{Quadratic lattices} A (quadratic) \emph{lattice} is a free abelian group $M$ equipped with a symmetric bilinear form with values in $\bbZ$. The \emph{orthogonal group} $\Or(M)$ of a lattice is defined in the natural way as the subgroup of automorphisms of the abelian group preserving the symmetric bilinear form.  More generally one defines in an obvious way an \emph{isometry, or isomorphism} of lattices.

Let $V$ be a real vector space equipped with  a symmetric bilinear form $(x,y)$ and  $(e_i)_{i\in I}$ be a basis in $V$ such that   $(e_i,e_j)\in \bbZ$ for all $i,j\in I$. Then  the  $\bbZ$-span $M$ of the basis is equipped naturally with the structure of a quadratic lattice. We have  already seen this construction  in the beginning of the section. The orthogonal group $\Or(M)$ coincides with the group introduced there. Obviously every quadratic lattice is obtained in this way by taking $V = M_\bbR = M\otimes_\bbZ \bbR$ and extending the bilinear form by linearity.

Recall some terminology in the theory of integral quadratic forms stated in terms of lattices. The signature of a lattice $M$ is the Sylvester signature $(t_+,t_-,t_0)$ of the corresponding real quadratic form on $V= M_\bbR.$  We omit $t_0$ if it is equal to zero.  A lattice with $t_- = t_0 = 0$ (resp. $t_+ = t_0 = 0$)  is called \emph{positive definite} (resp. \emph{negative definite}). A lattice  $M$ of signature  with $(1,a)$ or $(a,1)$ where  $a\ne 0$, is called \emph{hyperbolic} (or \emph{Lorentzian}).

All latices are divided in two types: \emph{even} if the values of its quadratic form are even and \emph{odd} otherwise. 

Assume that the lattice $M$ is nondegenerate, that is $t_0 = 0$.  This ensures that the map
\begin{equation}\label{iota}
\iota_M:M\to M^* = \Hom_\bbZ(M,\bbZ), \quad m\mapsto (m,?)\end{equation}
is injective. Since $M^*$ is an abelian group of the same rank as $M$, the quotient group
$$D_M = M^*/\iota(M)$$
is a finite group (the \emph{discriminant group} of the lattice $M$). Its order $d_M$ is equal to the absolute value of the \emph{discriminant} of  $M$ defined as the determinant of  a Gram matrix of the symmetric bilinear form of $M$.  A lattice is called \emph{unimodular} if the map \eqref{iota} is bijective (equivalently, if its discriminant is equal to $\pm 1$).

\begin{example}\label{rootlat2} Let $M$ be the lattice defining an integral structure on a finite reflection group from Example \ref{rootlat}. It is an even positive definite lattice for the groups of types $A,D, E$ and odd positive definite lattice for groups of type $B_n,F_4,G_2$. These lattices are called finite root lattices of the corresponding type.
\end{example}

\begin{example}\label{latpqr} Let $\Gamma$ be an irreducible  linear reflection group in $V$ admitting an integral structure $M$. It follows from \eqref{int} that, after rescaling the inner product in $V$, we may assume that 
$M$ is a lattice in $V$ with $M_\bbR = V$.  For example, consider the group $\Gamma = W(p,q,r)$ from Example \ref{triangle} as a linear reflection group in $\bbR^n$, where $n = p+q+r-2$. The unit vectors $e_i$ of a fundamental Coxeter polytope  satisfy
$ (e_i,e_j) = -2\cos \frac{\pi}{m_{ij}}$, where  $m_{ij} \in \{1,2,3\}$. Thus, rescaling the quadratic form in $V$ by multiplying its values by 2, we find fundamental root vectors $\alpha_i$ such that $(\alpha_i,\alpha_j)\in \bbZ$. The lattice $M$ generated by these vectors defines an integral structure of $\Gamma$. The Gram matrix  $G$ of the set of fundamental root vectors has 2 at the diagonal and $2I_n-G$ is the incidence matrix of the Coxeter graph of type $T_{p,q,r}$ from Example \ref{triangle}. We denote the lattice $M$ by $E_{p,q,r}$. One computes directly the signature of $M$ to obtain that $E_{p,q,r}$ is nondegenerate and positive definite if and only if $\Gamma$ is a finite reflection group of type $A,D, E$ ($r = 1 (A_n)$, or $r=p = 2 (D_n)$, or $r=2,p=3, q= 3,4,5 (E_6,E_7,E_8))$. 

The lattice $E_{p,q,r}$ is degenerate  if and only if it corresponds to a parabolic reflection group of type $ \tilde{E}_6,\tilde{E}_7,\tilde{E}_8$. The lattice $E_{p,q,r}^\perp$ is of rank 1 and 
$E_{p,q,r}/E_{p,q,r}^\perp$ is isomorphic to the lattice $E_{p-1,q,r}, E_{p,q-1,r}, E_{p,q,r-1}$, respectively. 

In the remaining cases $E_{p,q,r}$ is a hyperbolic lattice of signature $(n-1,1)$. 

It is also easy to compute the determinant of the Gram matrix to obtain that the absolute value of the discriminant of a nondegenerate lattice $E_{p,q,r}$ is equal to $|pqr-pq-pr-qr|$. In particular, it is a unimodular lattice
 if and only if $(p,q,r) = (2,3,5)$ or $(2,3,7)$.

One can also compute  the discriminant group of  a nondegenerate lattice $E_{p,q,r}$ (see \cite{Bries}). 
\end{example}

Every subgroup of $M$ is considered as a lattice with respect to the restriction of the quadratic form (\emph{sublattice}).  The orthogonal complement of a subset $S$ of a lattice is defined to be the set of  vectors $x$ in $M$ such that $(x,s) = 0$ for all $s\in S$. It is a primitive sublattice of $M$ (i.e. a subgroup of $M$ such that the quotient group is torsion-free). Also one naturally defines  the orthogonal direct sum $M\perp N$ of two (and finitely many) lattices.

The lattice $M$ of rank 2 defined by the matrix 
$\left(\begin{smallmatrix}0&1\\1&0\end{smallmatrix}\right)$
is denoted by $U$ and is called the \emph{hyperbolic plane}.

For any lattice $M$ and an integer $k$ we denote by $M(k)$ the lattice obtained from $M$ by multiplying its quadratic form by $k$. For any integer $k$ we denote by  $\langle k\rangle $ the lattice of rank 1 generated by a vector $v$ with $(v,v) = k$.

The following theorem describes the structure of unimodular indefinite lattices (see \cite{Se}).

\begin{theorem}\label{serre} Let $M$ be a unimodular lattice of indefinite signature $(p,q)$. If $M$ is odd, then it is isometric to the lattice $I_{p,q} = \langle 1\rangle^p\perp \langle-1\rangle^q$. If $M$ is even, then $p-q \equiv 0\ \mod 8$ and $M$ or $M(-1)$ is isometric to the lattice $II_{p,q}, p < q,$ equal to the orthogonal sum $U^{p}\perp E_8^{\frac{q-p}{8}}$.
\end{theorem}

\subsection{Reflection group of a lattice} Recall that the orthogonal group of a non-degenerate symmetric bilinear form on a finite-dimensional vector space over a field of characteristic $\ne 2$ is always generated by reflections. This  does not apply to orthogonal groups of lattices.

Let $M$ be a nondegenerate quadratic lattice. A \emph{root vector} in $M$ is a primitive vector $\alpha$ with $(\alpha,\alpha) \ne  0$ satisfying \eqref{ind}. A root vector $\alpha$ defines a reflection 
$$r_\alpha:x\mapsto x-\frac{2(\alpha,x)}{(\alpha,\alpha)}\alpha$$
in $V = M_\bbR$ which  leaves $M$ invariant.  Obviously any vector $\alpha$ with $(\alpha,\alpha) = \pm 1$ or $\pm 2$ is a root vector. 
Suppose $(\alpha,\alpha) = 2k$. The linear function $M\to \bbZ, x\mapsto \frac{2(\alpha,x)}{(\alpha,\alpha)}$ defines a nontrivial  element from $M^*/M$  of order   $k$. Thus  $k$ must divide the order of 
the discriminant group. In particular, all root vectors of a unimodular even lattice satisfy $(\alpha,\alpha) = \pm 2$.

We will be  interested in positive definite lattice or  hyperbolic lattices $M$ of signature $(n,1)$.  For such a  lattice we define the \emph{reflection group} $\Ref(M)$ of $M$ as the subgroup of $\Or(M)$ generated by reflections $r_\alpha$, where $\alpha$ is a root vector with $ (\alpha,\alpha) > 0$. We denote by $\Ref_k(M)$ its subgroup generated by reflections in root vectors with $k = (\alpha,\alpha)$ (the \emph{$k$-reflection subgroup}).  We set
$$\Ref_{-k}(M(-1)) = \Ref_k(M).$$
Each group $\Ref_k(M)$ is a  reflection groups in corresponding hyperbolic  or spherical space.

 Suppose $M$ is of signature $(n,1)$. Let 
$$\Or(M)^+ = \Or(M)\cap \Or(n,1)^+, \quad \Or(M) = \Or(M)^+\times \{\pm 1\},$$
where $\Or(n,1)^+$ is the subgroup of index 2 of $\Or(n,1)$ defined in section 2.2. Note that every $\Ref_k(M)$ is a normal subgroup of $\Or(M)^+$. 

Let $P$ be a fundamental polyhedron of $\Ref_k(M)$  in $H^n$. Since $\Or(M)$ leaves invariant the set of root vectors $\alpha$ with fixed  $(\alpha,\alpha)$, it leaves invariant the set of reflecting hyperplanes of $\Ref_k(M)$. Hence, for any $g\in \Or(M)^+$, there exists $s\in \Ref_k(M)$ such that $g(P) = s(P)$. This shows that
\begin{equation}\label{semi}
\Or(M)^+ = \Ref_k(M)\rtimes S(P),
\end{equation}
where $S(P)$ is the subgroup of $\Or(M)$ which leaves $P$ invariant.

\begin{example} Let $M = E_{p,q,r}$ with finite $\Ref(M)$. Then $\Ref(M) = W(p,q,r)$ from Example \ref{triangle}, where $(p,q,r) = (1,1,n) (A_n), (2,2,n-2) (D_n), (2,3,3) (E_6), (2,3,4) \\(E_7),(2,3,5) (E_8)$. 
We have $S(P) = \bbZ/2\bbZ (A_n, E_6, D_n, n \ge 5), S(P) = S_3 (D_4)$ and $S(P)$ is trivial for $E_7,E_8$. 

The standard notations for the finite reflection groups $W(p,q,r)$ are $W(T)$, where $T = A_n,D_n,E_6,E_7,E_8$. The corresponding lattices $E_{p,q,r}$ are called \emph{finite root lattices}. Their reflection groups are the Weyl groups of the corresponding root systems. 

In general $\Ref(E_{p,q,r})$ is larger than the group $W(p,q,r)$ (see Example \ref{ex4}).
\end{example}

\begin{example}\label{leech} Let $M = II_{25,1}$ be an even  unimodular  hyperbolic lattice of rank 26. According to Theorem \ref{serre}
$$II_{25,1} \cong U\perp E_8^3$$
The lattice $II_{25,1}$ contains as a direct summand an even positive definite unimodular lattice $\Lambda$ of rank $24$ with  $(v,v)\ne 2$ for all $v\in \Lambda$. A lattice with such properties (which determine uniquely the isomorphism class) is  called a \emph{Leech lattice}. Thus $II_{25,1}$ can be also described as
\begin{equation}\label{leechd}
II_{25,1} = U\perp \Lambda.\end{equation}
The description of $\Ref(II_{25,1}) = \Ref_2(II_{25,1})$ was given by J. Conway \cite{Con}. Since $II_{25,1}$ is unimodular and even, all root vectors satisfy $(\alpha,\alpha) = 2$. The  group admits a fundamental polytope $P$ whose reflection hyperplanes  are orthogonal to the \emph{Leech roots}, i.e. root vectors of the form $(f-(1+\tfrac{(v,v)}{2})g,v)$, where $v\in \Lambda$ and $f,g$ is a basis of $U$ with Gram matrix $\left(\begin{smallmatrix}0&1\\
1&0\end{smallmatrix}\right)$. In other words, a choice of a decomposition \eqref{leechd} defines a fundamental polyhedron for the reflection group with fundamental roots equal to the Leech vectors. We have
$$\Or(II_{25,1})^+ = \Ref(II_{25,1})\rtimes S(P),$$
where $S(P) \cong \Lambda\rtimes \Or(\Lambda).$
\end{example}

Define a hyperbolic  lattice $M$ to be \emph{reflective} if its root vectors span $M$ and $\Ref(M)$ is of finite covolume (equivalently, its index in $\Or(M)$ is finite).\footnote{This definition is closely related but differs from the definition of reflective hyperbolic lattices used in the works of V. Gritsenko and V. Nikulin  on Lorentzian Kac-Moody algebras \cite{GN}.} 
 In the hyperbolic case the first condition follows from the second one. It is clear that the reflectivity property of $M$ is preserved when we scale $M$, i.e. replace $M$ with $M(k)$ for any positive integer $k$.
The following nice result is due to F.  Esselmann \cite{Es}.

\begin{theorem}\label{ess} Reflective lattices of signature $(n,1)$ exist only if $n \le 19$, or $n = 21$. 
\end{theorem}

The first example of a reflective lattice of rank 22 was given by Borcherds \cite{Bo1}. We will discuss this lattice later.

An important tool in the classification (yet unknown) of reflective lattices is the following lemma of Vinberg 
\cite{Vi1}.

\begin{lemma} Let $M$ be a hyperbolic reflective lattice. For any isotropic vector $v\in M$ the lattice
$v^\perp/\bbZ v$ is a definite reflective lattice.
\end{lemma}

Another useful result is the following (see \cite{Bu}).

\begin{theorem} Suppose a reflective hyperbolic lattice $M$ decomposes as an orthogonal sum of a lattice $M'$ and a definite lattice $K$. Then $M'$ is reflective.
\end{theorem}

\begin{examples} 
\noindent
1) A lattice of rank 1 is always reflective.

\smallskip\noindent
2) The lattice $U$ is reflective. The group $\Or(U)$ is finite.

\smallskip\noindent
3) All finite root lattices and their orthogonal  sums are reflective. 

\smallskip\noindent
4) An odd lattice $I_{n,1} = \la 1\ra^n\perp \la-1\ra$ is reflective if and only if $n\le 20$. The Coxeter diagrams of their reflection groups can be found in \cite{VS}, Chapter 6, \S2 ($n\le 18$) and in \cite{VK} ($n = 18,19$). Some of these diagrams are also discussed in \cite{CS}, Chapter 28. The lattices $I_{n,1}(2)$ are reflective for $n\le 19$ but 2-reflective only for $n\le 9$.   For example, the following is the Coxeter diagram of the reflection group of the lattice $I_{16,1}$

\begin{figure} [ht]
\xy
(-65,0)*{};
@={(-20,0),(-10,0),(10,0),(20,0),(0,20),(0,10),(0,-10),(0,-20),(14.15,14.15),(-14.15,-14.15),
(-14.15,14.15), (14.15,-14.15),(18,7.4),(-18,-7.4),(-18,7.4),(18,-7.4),(7.4,18),(-7.4,-18),(7.4,-18),(-7.4,18)}@@{*{\bullet}};
(-20,0)*{};(-18,7.4)*{}**\dir{-},(-18,7.4)*{};(-14.15,14.15)*{}**\dir{-};
(-14.15,14.15)*{};(-7.4,18)*{}**\dir{-},(-7.4,18)*{};(0,20)*{}**\dir{-};
(20,0)*{};(18,7.4)*{}**\dir{-},(18,7.4)*{};(14.15,14.15)*{}**\dir{-};
(14.15,14.15)*{};(7.4,18)*{}**\dir{-},(7.4,18)*{};(0,20)*{}**\dir{-};
(-20,0)*{};(-18,-7.4)*{}**\dir{-},(-18,-7.4)*{};(-14.15,-14.15)*{}**\dir{-};
(-14.15,-14.15)*{};(-7.4,-18)*{}**\dir{-},(-7.4,-18)*{};(0,-20)*{}**\dir{-};
(20,0)*{};(18,-7.4)*{}**\dir{-},(18,-7.4)*{};(14.15,-14.15)*{}**\dir{-};
(14.15,-14.15)*{};(7.4,-18)*{}**\dir{-},(7.4,-18)*{};(0,-20)*{}**\dir{-};
(-20,0)*{};(-10,0)*{}**\dir{-};(20,0)*{};(10,0)*{}**\dir{-};
(0,10)*{};(0,20)*{}**\dir2{-};(0,-10)*{};(0,-20)*{}**\dir2{-};
(0,-10)*{};(0,10)*{}**\dir2{-};(0,-10)*{};(0,10)*{}**\dir2{-};
(0,-10)*{};(0,10)*{}**\dir2{-};(0,-10)*{};(0,10)*{}**\dir{-};
(.05,-10)*{};(.05,10)*{}**\dir2{-};(-.05,-10)*{};(-.05,10)*{}**\dir2{-};

\endxy

\caption{}\label{fig3}
\end{figure}
It is easy to see that this is possible only if all vertices not connected by the thick line correspond to roots $\alpha$ with $(\alpha,\alpha) = 2$ and the remaining two vertices correspond to roots with 
$(\alpha,\alpha) = 1$.

Many other examples of Coxeter diagrams for 2-reflective lattices can be found in \cite{Ni2}.

\smallskip\noindent
5) Many examples (almost a classification) of reflective lattices of rank 3 and 4 can be found in \cite{SW}, \cite{Sh}.
\end{examples}

 All even hyperbolic lattices of rank $r > 2$ for which $\Ref_2(M)$ is of finite covolume (\emph{2-reflective lattices}) were found by V. Nikulin \cite{Ni2}  ($r\ne 4$) and E. Vinberg (unpublished) ($r = 4$) (a survey of Nikulin's results can be found in \cite{Dos}). They exist only in dimension $\le 19$. 

\begin{example}\label{ex4} A hyperbolic lattice $E_{2,3,r}$ is 2-reflective if and only if $7\le r \le 10$. A hyperbolic lattice $E_{2,4,r}$ is 2-reflective if and only if $r = 5,6,7$. A hyperbolic lattice $E_{3,3,r}$ is 2-reflective if and only if $r =4,5,6$. This easily follows  from Proposition \ref{bugaenko}.
The reflection groups of the  lattices $E_{2,3,7}, E_{3,3,4}$ and $E_{2,4,5}$ are quasi-Lanner and coincide with the Coxeter groups $W(2,3,7), W(3,3,4), W(2,4,5)$. The reflection groups of other lattices are larger than the corresponding groups $W(p,q,r)$. For example, the Coxeter diagram of $\Ref_2(E_{2,3,8}) = \Ref(E_{2,3,8})$ is the following one.

\begin{figure}[ht]
\xy (-30,0)*{};
@={(0,0),(7.5,0),(15,0),(22.5,0),(30,0),(37.5,0),(45,0),(52.5,0),(60,0),(67.5,0),(75,0),(82.5,0)}@@{*{\bullet}};
(0,0)*{};(82.5,0)**{}**\dir{-};
(15,0)*{};(15,-7.5)*{\bullet}**\dir{-};
(75,0)*{};(82.5,0)*{}**\dir{-};
(75,.1)*{};(82.5,.1)*{}**\dir2{-};
(75,-.1)*{};(82.5,-.1)*{}**\dir2{-};
(75,.05)*{};(82.5,-.05)*{}**\dir2{-};
(75,-.05)*{};(82.5,-.05)*{}**\dir2{-};

\endxy

\caption{}\label{fig5}

\end{figure}
To see this one first observes, using Proposition \ref{bugaenko},  that the Coxeter diagram defines a group of cofinite volume. The root vector corresponding to the extreme right vertex is equal to the vector $r+e$, where $r$ generates the kernel of the lattice $E_{2,3,6}$ embedded naturally in $E_{2,3,8}$ and the vector $e$ is the root vector corresponding to the extreme right vector in the subdiagram defining the sublattice isomorphic to $E_{2,3,7}$. 
\end{example}

\section{Automorphisms of algebraic surfaces}
\subsection{Quadratic lattices associated to an algebraic surface}
Let $X$ be a complex  projective algebraic surface. It has the underlying structure of a compact smooth oriented $4$-manifold. Thus the cohomology group $H^2(X,\bbZ)$ is a finitely generated abelian group  equipped with a unimodular symmetric bilinear form
$$H^2(X,\bbZ)\times H^2(X,\bbZ) \to \bbZ$$
defined by the cup-product. 

When we divide $H^2(X,\bbZ)$ by the torsion subgroup  we obtain a quadratic lattice $H_X$.  We will denote the value of the bilinear  form on $H_X$ induced by the cup product by $x\cdot y$ and write $x^2$ if $x=y$.

To compute its signature one uses the \emph{Hodge decomposition} (depending on the complex structure of $X$)
$$H^2(X,\bbC) = H^{2,0}(X)\oplus H^{1,1}(X,\bbC)\oplus H^{0,2}(X),$$
where $\dim H^{2,0} = \dim H^{0,2}$ and is equal to the dimension $p_g(X)$ of the space of holomorphic differential 2-forms on $X$. It is known that under the complex conjugation on $H^2(X,\bbC)$ the space $H^{1,1}(X)$ is invariant and the space $H^{2,0}(X)$ is mapped to $H^{0,2}(X)$, and vice versa. One can also compute the restriction of the cup-product on $H^2(X,\bbC)$ to each $H^{p,q}$ to conclude that the signature of the cup-product on $H^2(X,\bbR)$ is equal to $(b_2^+,b_2^-)$, where $b_2^+ = 2p_g+1$. 

The parity of the lattice $H_X$ depends on the property of its first Chern class $c_1(X)\in H^2(X,\bbZ)$. If is divisible by 2 in $H^2(X,\bbZ)$, then $H_X$ is an even lattice. 

The lattice $H_X$ contains an important primitive sublattice which depends on the complex structure of $X$. For any complex irreducible curve $C$ on $X$ its fundamental class $[C]$ defines a cohomology class in $H^2(X,\bbZ)$. The $\bbZ$-span of these classes defines a subgroup of $H^2(X,\bbZ)$ and its image $S_X$ in $H_X$ is called the \emph{Neron-Severi lattice} (or the \emph{Picard lattice}) of $X$.  It is a sublattice of $H_X$ of signature $(1,\rho-1)$, where $\rho = \rank~S_X$. An example of an element of positive norm is the class of a hyperplane section of $X$ in any projective embedding of $X$. Thus the lattice $S_X$ is hyperbolic in the sense of previous sections and we can apply the theory of reflection groups to $S_X$. 

It is known that the image in $H_X$ of the first Chern class $c_1(X)\in H^2(X,\bbZ)$ belongs to the Picard lattice. The negative $K_X = -c_1(X)$ is  the \emph{canonical class} of $X$.  It is equal  to  the image in $H^2(X,\bbZ)$ of a divisor of zeros and poles of a holomorphic differential 2-form on $X$. We denote the image of $K_X$ in $H_X$ by $k_X$. For any $x\in S_X$ we define
$$p_a(x) =  x^2+x\cdot k_X.$$ 
It is always an even integer.  If $x$ is the image in $H_X$ of the fundamental class of a nonsingular complex curve $C$ on $X$, then $p_a(x)$  is equal  to $2g-2$, where $g$ is the genus of the Riemann surface $C$  (the \emph{adjunction formula}). If $C$ is an irreducible complex curve with finitely many singular points, then this number is equal to $2g_0-2+2\delta$, where $g_0$ is the genus of a nonsingular model of $C$ (e.g. the normalization) and $\delta$ depends on the nature of singular points of $C$ (e.g. equal to their number if all singular points are ordinary nodes or cusps). This implies that $p_a(x)$ is always even. In particular, the sublattice 
\begin{equation}\label{canort}
S_X^0 = k_X^\perp = \{x\in S_X: x\cdot k_X = 0\}.
\end{equation}
is always even. Its signature is equal to 
\[
\sign(S_X^0) = \begin{cases}
 (1,\rho-1) & \text{if}\  k_X = 0, \\
    (0,\rho-2,1)  &\text{if}\ k_X^2 = 0, k_X\ne 0, \\
    (0,\rho-1) &\text{if}\ k_X^2 > 0,\\
    (1,\rho-2)&\text{if}\ k_X^2 < 0.
\end{cases}
\]
A \emph{$(-1)$-curve} (resp. \emph{$(-2)$-curve}) on $X$ is a nonsingular irreducible curve $C$ of genus $0$ (thus isomorphic to $\bbP^1(\bbC)$) with $[C^2] = -1$ (resp. $-2$). By the adjunction formula this is equivalent to $C\cong \bbP^1(\bbC)$ and $[C]\cdot K_X = -1$ (resp. $0$). A $(-1)$-curve curve appears as a fibre of a blow-up map  $f:X\to  Y$ of a point on a nonsingular algebraic surface $Y$. A $(-2)$-curve appears as a fibre of a resolution of an ordinary double point on a complex surface $Y$.

A surface $X$ which does not contain $(-1)$-curves  is called \emph{minimal}. The following results follow from the Enriques-Kodaira classification of complex algebraic surfaces (see \cite{BPV}).  

\begin{theorem}\label{clas} Let $X$ be a minimal complex algebraic surface. Then one of the following cases occurs. 
\begin{enumerate}
\item $X\cong \bbP^2$ or there exists a regular map $f:X\to B$ to some nonsingular curve $B$ whose fibres are isomorphic to $\bbP^1$. Moreover $S_X= H_X$ and we are in precisely one  of the following cases.
\begin{itemize}
\item [(i)] $X\cong \bbP^2$, $S_X \cong  \langle 1\rangle $, and $k_X = 3a$, where $a$ is a generator.
\item [(ii)] $S_X\cong U$, and $k_X = 2a+2b$, where $a,b$ are generators of $S_X$ with the Gram matrix $\left(\begin{smallmatrix}0&1\\
1&0\end{smallmatrix}\right)$ .
\item[(iii)] $S_X \cong I_{1,1} = \langle1\rangle\perp \langle-1\rangle$  and $k_X = 2a+2b$, where $a,b$ are generators of $S_X$ with the Gram matrix $\left(\begin{smallmatrix}1&0\\
0&-1\end{smallmatrix}\right)$.
\end{itemize}
\item $k_X = 0$;
\begin{itemize}
\item [(i)] $H_X \cong U^{3}\perp E_8(-1)^{2}$, $S_X$ is an even lattice of signature $(1,\rho-1)$, where $1\le \rho \le 20$;
\item [(ii)] $H_X = S_X \cong U \perp  E_8(-1)$;
\item [(iii)]  $H_X \cong U^{3}$, $S_X$ is an even lattice of signature $(1,\rho-1)$, where $1\le \rho \le 4$;
\item [(iv)] $H_X = S_X \cong U$.
\end{itemize}
\item $k_X\ne 0, k_X^2 = 0$,  $S^0/\bbZ k_X$ is a negative definite lattice.
\item [(4)] $k_X^2 > 0, K_X\cdot [C] \ge 0$ for any curve $C$ on $X$, 
$S_X^0$ is a negative definite lattice.
\end{enumerate}
\end{theorem}
 
 The four cases (1)-(4) correspond to the four possible values of the \emph{Kodaira dimension} $\kappa(X)$ of $X$ equal to $-\infty, 0, 1,2$, respectively. Recall that $\kappa(X)$ is equal to the maximal possible dimension of the image of $X$ under a rational map given by some multiple of the canonical linear system on $X$. The four subcases (i)-(iv) in (2) correspond to K3-surfaces, Enriques surfaces, abelian surfaces (complex algebraic tori), and hyperelliptic surfaces, respectively.

Let $\Aut(X)$ denote the group of automorphisms of $X$ (as an algebraic variety or as a complex  manifold). Any  $g\in \Aut(X)$  acts naturally on $H_X$ via the pull-backs of cohomology classes. Since the latter  is compatible with the cup-product, the action preserves the structure of a quadratic lattice on $H_X$ and also leaves invariant the sublattice $S_X$. This defines a homomorphism
\begin{equation}\label{aut}
\Aut(X) \to \Or(S_X), \ g\mapsto (g^*)^{-1}.
\end{equation}
Since $g^*(K_X) = K_X$, we see that the image of this homomorphism is contained in the stabilizer subgroup $\Or(S_X)_{k_X}$ of the vector $k_X$. In particular, it induces a homomorphism
\begin{equation}\label{auto}
a:\Aut(X) \to \Or(S_X^0).
\end{equation}
The group $\Aut(X)$ is a topological group whose connected component of the identity $\Aut(X)^0$ is a complex Lie group. One can show that $\Aut(X)^0$ acts identically on $S_X$ and the kernel of the induced map of the quotient group $\Aut(X)/\Aut(X)^0$ is finite (see \cite{Do2}). The group $\Aut(X)^0$ can be nontrivial only  for  surfaces of Kodaira dimension $-\infty$, or abelian surfaces, or surfaces of Kodaira dimension 1 isomorphic to some finite quotients of the products of two curves, one of which is of genus 1. It follows from Theorem \ref{clas} that $\Aut(X)$ is always finite for surfaces of Kodaira dimension 2. 

\subsection{Rational surfaces}\label{rat}
A rational surface $X$ is a nonsingular projective algebraic surface birationally isomorphic to $\bbP^2$. We will be interested only in  \emph{basic rational surfaces}, i.e. algebraic surfaces admitting a regular birational map $\pi:X\to \bbP^2$.\footnote{For experts:non basic rational surfaces are easy to describe, they are either minimal rational surfaces different from $\bbP^2$ or surfaces obtained from minimal ruled surfaces $\bfF_n, n\ge 2,$ by blowing up points on the exceptional section and their infinitely near points. The automorphism groups of non-basic rational surfaces are easy to describe and they are rather dull.}

It is known that any birational regular map of algebraic surfaces is equal to the composition of blowing-ups of  points (\cite{Ha}). Applying this to the map $\pi$, we obtain a factorization
\begin{equation}\label{decom}
\pi:X = X_N\overset{\pi_N}{\longrightarrow} X_{N-1}\overset{\pi_{N-1}}{\longrightarrow} \ldots \overset{\pi_2}{\longrightarrow} X_1\overset{\pi_1}{\longrightarrow} X_0 = \bbP^2,
\end{equation}
where $\pi_i:X_i\to X_{i-1}$ is the blow-up of a point $x_i\in X_{i-1}$. Let 
\begin{equation}\label{excl1}
E_i = \pi_i^{-1}(x_i),\quad \calE_i = (\pi_{i+1}\circ \ldots\pi_{N})^{-1}(E_i).
\end{equation}
For any rational surface $S_X = H_X$. Let $e_i$ denote the cohomology class  $[\calE_i]$ of the (possibly reducible) curve $\calE_i$. It satisfies $e_i^2 =  e_i\cdot k_X = -1$. One easily checks that $e_i\cdot e_j = 0$ if $i\ne j$. Let $e_0 = \pi^*([\ell])$, where $\ell$ is a line in $\bbP^2$. We have $e_0\cdot e_i = 0$ for all $i$. The classes $e_0,e_1,\ldots,e_N$ form a basis in $S_X$ which we call a \emph{geometric basis}. The Gram matrix of a geometric basis  is the diagonal matrix $\diag[1,-1,\ldots,-1]$. Thus the factorization \eqref{decom} defines an isomorphism of quadratic lattices
$$\phi_\pi:I_{1,N}\to S_X, \ \bfe_i\mapsto e_i,$$
where $\bfe_0,\ldots,\bfe_N$ is the standard basis of $I_{1,N}.$
It follows from the formula for the behavior of the canonical class under a blow-up that 
$k_X$ is equal to the image   of the vector 
$$k_N = -3\bfe_0+\bfe_1+\ldots+\bfe_N.$$
This implies that the lattice $S_X^0$ is isomorphic to the orthogonal complement 
$k_N^\perp$ in $I_{1,N}$. 

For $N\ge 3$, the vectors
\begin{equation}\label{alpha}
\bfa_1 = \bfe_0-\bfe_1-\bfe_2-\bfe_3, \quad \bfa_2 = \bfe_1-\bfe_2, \quad \ldots \quad, \bfa_{N} = \bfe_{N-1}-\bfe_N
\end{equation}
form a basis of $k_N^\perp$ with Gram matrix equal to $-2C$, where $C$ is the Gram matrix of the Coxeter group $W(E_N):= W(2,3,N-3)$ if $N\ge 4$ and $W(E_3) = W(A_2\times A_1)$ if $N = 3$.
This embeds the lattice 
\begin{equation}\label{en}
E_N = \begin{cases}
      A_2\perp A_1& \text{if}\ N= 3, \\
       E_{2,3,N-3}&\text{if}\  N\ge 4,
\end{cases}
\end{equation}
 in the lattice $I_{N,1}$ with orthogonal complement generated by $k_N$. The restriction of $\phi_\pi$ to $k_N^\perp$ defines   an isomorphism of lattices 
$$\phi_\pi:E_{N}(-1)\to S_X^0.$$
Let us identify the Coxeter group $W(E_N)$  with the subgroup of $\Or(I_{1,N})$ generated by the reflections in the vectors $\bfa_i$ from \eqref{alpha}. A choice of a geometric basis in $S_X$ defines an isomorphism from $W(E_N)$ to a subgroup of $\Or(S_X^0)$ generated by reflections in vectors 
$\alpha_i = \phi_\pi(\bfa_i)$.  It is contained in the reflection group $\Ref_{-2}(S_X^0)$.

\begin{theorem}\label{crem} The image $W_X$ of $W(E_N)$ in $\Or(S_X^0)$ does not depend on the choice of a geometric basis. The image of the homomorphism $a:\Aut(X)\to \Or(S_X^0)$ is contained in $W_X$.
\end{theorem}

\begin{proof} To prove the first assertion it suffices   to show that the transition matrix of two geometric bases defines an orthogonal transformation of $I_{1,N}$ which is the product of reflections in vectors $\bfa_i$. Let $(e_0,\ldots,e_N)$ and $(e_0',\ldots,e_N')$ be two geometric bases and 
$$e_0' = m_0e_0-m_1e_1-\ldots-m_Ne_N.$$ 
For any curve $C$ on $X$, we have $e_0'\cdot [C]) = \pi^*([\ell])\cdot [C]) = [\ell]\cdot [\pi'(C)]\ge 0$. Intersecting $e_0'$ with $e_i$ we obtain that $m_{0}> 0, m_{i} \ge 0, i > 0.$ We have
 \begin{eqnarray}\label{noether1}
1 &= &e_0^2 =  m_{0}^2-m_{1}^2-\ldots-m_{N}^2, \\ \notag
-3& = &e_0\cdot k_X = -3m_{0}+m_{1}+\ldots+m_{N}. 
\end{eqnarray}
Applying the reflections in vectors $\alpha_i = \phi_{\pi}(\bfa_i), i > 1$, we may assume that $m_1\ge m_2\ge\ldots \ge m_N$. Now we use the following inequality (\emph{Noether's inequality})
\begin{equation}\label{ne}
m_0 > m_1+m_2+m_3 \ \text{if}\  m_0 > 1.\end{equation}
To see this we multiply the second equality in \eqref{noether1} by $m_3$ and then subtract from the first one to get
$$m_1(m_1-m_3)+m_2(m_2-m_3)-\sum_{i\ge 4}m_i(m_3-m_i) = m_0^2-1-3m_3(m_0-1).$$
This gives
$$(m_0-1)(m_1+m_2+m_3-m_0-1) = (m_1-m_3)(m_0-1-m_1)+(m_2-m_3)(m_0-1-m_2)+$$
$$+\sum_{i\ge 4}m_i(m_3-m_i).$$
The first  inequality in \eqref{noether1} implies that $m_0^2-m_1^2> 0$, hence 
 $m_0-m_i\ge 1$. Thus the right-hand-side is non-negative, so the left-hand-side too.  This  proves the claim. Now consider the reflection $s = r_{\alpha_1}$. Applying it to $e_0'$ we get 
 $$s(e_0') = e_0'+((e_0-e_1-e_2-e_3)\cdot e_0')(e_0-e_1-e_2-e_3) = $$
 $$(2m_0-m_1-m_2-m_3)e_0-(m_0-m_2-m_3)e_1-(m_0-m_1-m_3)e_2-(m_0-m_1-m_2)e_3.$$
 Using \eqref{ne}, we obtain that the  matrix $S\cdot A$, where $S$ is the matrix of $s$, has the first column equal to $(m_0',-m_1,-\ldots,-m_N')$ with $m_0' < m_0$ and $m_i'\ge 0$ for $i > 0$. Since our transformations are isometries of $S_X$, the inequalities \eqref{noether1} hold for the vector 
$(m_0',-m_1,-\ldots,-m_N')$. 
So, we repeat the argument in order to decrease $m_0'$. After finitely many steps we get the transformation with first column vector equal to $(1,0,\ldots,0)$. 
Now the matrix being orthogonal matrix of the quadratic form $x_0^2-x_1^2-\ldots-x_N^2$ must have the first row equal to $(1,0,\ldots,0)$. Thus the remaining rows and columns define an orthogonal matrix of the quadratic form $x_1^2+\ldots+x_N^2$ with integer entries. 
This implies that, after reordering the columns and the rows we get a matrix with $\pm 1$ at the diagonal and zero elsewhere. It remains to use that the transformation leaves the vector $k_N$ invariant to conclude that the matrix is the identity. Thus the transition matrix is the product of the matrices corresponding to reflections in vectors $\alpha_i, i = 1,\ldots, N$.

Let us prove the  last  statement. Suppose $g^*$ is the identity on $S_X$. Then $g^*(e_N) = [g^{-1}(\calE_N)] = e_N$. Since $\calE_N^2< 0$, it is easy to see that $\calE_N$ is homologous to $g^{-1}(\calE_N)$ only if $g(\calE_N) = \calE_N$. This implies that $g$ descends  to the surface $X_{N-1}$. Replacing $X$ with $X_{N-1}$ and repeating the argument, we see that $g$ descends to $X_{N-2}$. Continuing in this way we obtain that $g^*$ is the identity and descends to a projective automorphism $g'$ of $\bbP^2$. If all curves $\calE_i$ are irreducible, their images on $\bbP^2$ form an ordered  set of $N$ distinct points which must be preserved under $g'$. Since a square matrix of size $3\times 3$ has at most 3 linear independent eigenvectors, we see that $g'$, and hence $g$ must be the identity. The general case requires a little more techniques to prove and we omit the proof.

\end{proof}

The main problem is to describe all possible subgroups of the Weyl group $W(E_N)$ which can be realized as the image of a group $G$ of automorphisms of a rational surface obtained by blowing-up $N$ points in the plane.

First of all we may restrict ourselves to \emph{minimal pairs} $(X,G\subset \Aut(X))$. Minimal means that any $G$-equivariant  birational regular map $f:X\to X'$ of rational surfaces  must be an isomorphism. A factorization \eqref{decom} of  $\pi':X'\to \bbP^2$ can be extended to a factorization  $\pi:X\to \bbP^2$, in such a way that the geometric basis $\phi_{\pi'}(\bfa)$ of $S_{X'}$  can be extended to a geometric basis $\phi_\pi(\bfa)$ of $S_X$. This  gives a natural inclusion $W_{X'}\subset W_X$ such that the image of $G$ in $W_{X'} $ coincides with the image of $G$ in $W_X$.

The  next result  goes back to the classical work of S. Kantor \cite{Ka} and now easily follows from an equivariant version of Mori's theory of minimal models. 

\begin{theorem} Let $G$ be a finite group of automorphisms of a rational surface $X$ making a minimal pair $(X,G)$.  Then either $X\cong \bbP^2$, or $X$ is a conic bundle with $(S_X)^G \cong \bbZ^2$ or $\bbZ$, or $X$ is a Del Pezzo surface with $(S_X)^G \cong \bbZ$. \end{theorem} 

Here a \emph{Del Pezzo  surface} is a rational surface $X$ with ample  $-K_X$.\footnote{This means that in some projective embedding of $X$ its positive multiple  is equal to the fundamental class of a hyperplane section.} Each Del Pezzo surface is isomorphic to either $\bbP^2$ or $\bbP^1\times \bbP^1$, or admits a factorization \eqref{decom} with $N \le 8$ and the images of the points $x_1,\ldots,x_N$ in $\bbP^2$ are all distinct  and satisfy 
\begin{itemize}
\item no three are on a line;
\item no six are on a conic;
\item not contained on a  plane cubic with one of them being its singular point $(N= 8)$.
\end{itemize}

\smallskip
A  \emph{conic bundle} is a rational surface which admits a regular map to a nonsingular curve with fibres isomorphic to a conic (nonsingular or the union of two distinct lines).

 A partial classification of  finite groups $G$ which can be realized as  groups of automorphisms of some  rational surface was given by  S. Kantor (for a complete classification and the history of the problem  see  \cite{DI}).

\begin{example} Let $X$ be a Del Pezzo surface with $N+1 = \rank~ S_X$. We refer to the number $9-N$ as the \emph{degree} of $X$. A Del Pezzo surface of degree $d > 2$ is isomorphic to a nonsingular surface of degree $d$ in $\bbP^d$. The most famous example is a cubic surface in $\bbP^3$ with  27 lines on it. The Weyl group $W(E_6)$ is isomorphic to the \emph{group of 27 lines on a cubic surface}, i.e. the subgroup of the permutation group $\Sigma_{27}$ which preserves the incidence relation between the lines.  Although the group of automorphisms of a general cubic surface is trivial, some special cubic surfaces admit nontrivial finite automorphism groups. All of them have been essentially classified in the 19th century.
\end{example}

The situation with infinite groups is more interesting and difficult. Since $E_N$ is negative definite for $N\le 8$, a basic rational surface $X$ with infinite automorphism group is obtained by blowing up $N \ge 9$ points. It is known that when the points are in general position, in some precisely defined sense, the group $\Aut(X)$ is trivial  \cite{Hi}, \cite{Ko}. So  surfaces with nontrivial automorphisms are obtained by blowing up a set of points  in some special position.

\begin{example}\label{ex1}Let $X$ be obtained by blowing up $9$ points $x_1,\ldots,x_9$ contained in two distinct irreducible plane cubic curves $F,G$. The surface admits a fibration $X\to \bbP^1$ with general fibre an elliptic curve. The image of each fibre in $\bbP^2$ is a plane cubic from the pencil of cubics spanned by the curves $F,G$. The exceptional curves $\calE_1,\ldots,\calE_9$ are sections of this fibration. Fix one of them, say $\calE_1$ and equip each nonsingular fibre $X_t$ with the group law with the zero point $X_t\cap \calE_1$. Take a point $x\in X_t$ and consider the sum $x+p_i(t)$, where $p_i(t) = X_t\cap \calE_i$. This defines an automorphism on an open subset of $X$ which can be extended to an automorphism $g_i$ of $X$. When the points $x_1,\ldots,x_9$ are general enough, the automorphisms $g_2,\ldots,g_9$ generate a free abelian group of rank 8. In general case it is a finitely generated group of rank $\le 8$.  In the representation of $\Aut(X)$ in $W_X\cong W(E_9)$ the image of this group is the lattice subgroup of the euclidean reflection group of type $E_8$.

This example can be generalized by taking general points $x_1,\ldots,x_9$ with the property that there exists an irreducible curve of degree $3m$ such that each $x_i$ is its singular point of multiplicity $m$. In this case the image of $\Aut(X)$ in $W(E_9)$ is the subgroup of the lattice subgroup $\bbZ^8$ such that the quotient group is isomorphic to $(\bbZ/m\bbZ)^8$ (see \cite{DO}, \cite{Gi}).

\end{example} 

\begin{example}\label{ex2} Let $X$ be obtained by blowing up 10  points $x_1,\ldots,x_{10}$ with the property that there exists an irreducible curve of degree 6 with ordinary double points at each $x_i$ (a \emph{Coble surface}).
The image of $\Aut(X)$ in $W_X\cong W(E_{10})$ is the subgroup
\begin{equation}\label{cs}
W(E_{10})(2) = \{g\in W(E_{10}):g(x)-x\in 2E_{10}\ \text{for all}\ x\in E_{10}\}.\end{equation}
(see \cite{Coble}). This group is the smallest normal subgroup which contains  the involution of the lattice $E_{10} = U\perp E_8$ equal to $(\id_U,-\id_{E_8})$.
\end{example}

\begin{example}\label{ex3} Let $(W,S)$ be a Coxeter system with finite set $S$. A \emph{Coxeter element} of $(W,S)$ is the product of elements  of $S$ taken in some order. Its conjugacy class does not depend on the order if  the Coxeter diagram is a tree (\cite{Bou}, Chapter V, \S 6, Lemma 1).  The order of a Coxeter element is finite if and only if $W$ is finite. Let $h_N$ be a Coxeter element of $W= W(E_N)$. In a recent paper \cite{Mc} C. McMullen realizes $h_N$ by an automorphism of a rational surface. The corresponding surface $X$ is obtained by blowing up $N$ points in special position lying on a cuspidal plane cubic curve. One can check that  for $N = 9$ or $10$ a Coxeter element does not belong to the subgroup described in the previous two examples. It is not known whether a rational surface realizing a Coxeter element   is unique up to isomorphism for $N\ge 9$. It is known to be unique for $N\le 8$. For example, for $N = 6$ the surface is isomorphic to the cubic surface 
$$T_0^3+T_1^3+T_3^2T_1+T_2^2T_3 = 0.$$
The order of $h_6$ is equal to $12$.
\end{example}

Until very recently, all known examples of minimal pairs $(X,G)$ with infinite  $G$ satisfied  the following condition:
\begin{itemize}
\item There exists $m > 0$ such that the linear system $|-mK_X|$ is not empty, or, in another words, the cohomology class $mc_1(X)$ can be represented by an algebraic curve.
\end{itemize}

 Without the minimality condition the necessity of this condition was conjectured by M. Gizatullin but a counter-example was found by B. Harbourne \cite{Harb}. A recent preprint of Eric  Bedford and Kyounghee Kim \cite{BK} contains an example of a minimal surface with infinite automorphism group with $|-mK_X| =\emptyset$ for all $m > 0$.

\subsection{K3 surfaces} 
These are surfaces from case 2 (i) of Theorem \ref{clas}. They are characterized by the conditions
$$c_1(X) = 0, \quad H^1(X,\bbC) = 0.$$
In fact, all K3 surfaces are simply connected and belong to the same diffeomorphism type. 

\begin{examples} 
1) $X$ is a nonsingular surface of degree 4 in $\bbP^3$;

\smallskip\noindent
2) $X$ is the double cover of a rational surface $Y$ branched along a nonsingular curve $W$ whose cohomology class $[W]$ is equal  to $-2K_X$. For example, one may  
take $Y = \bbP^2$ and $W$ a nonsingular curve of degree 6. Or, one takes $Y = \bbP^1\times \bbP^1$ and $W$ a curve of bi-degree $(4,4)$.

 \smallskip\noindent
 3) Let $A$ be a compact complex torus which happens to be a projective algebraic variety. This is a surface from case 2 (iii). The involution $\tau:a\mapsto -a$ has 16 fixed points and the orbit space 
 $A/(\tau)$ acquires 16 ordinary double points. A minimal nonsingular surface   birationally equivalent to the quotient is a K3 surface, called \emph{Kummer surface} associated to $A$.
\end{examples}

Let $\Aut(X)$ be the group of biregular automorphisms of $X$. It is known that $X$ does not admit nonzero  holomorphic vector fields and hence the Lie algebra of the maximal Lie subgroup of $\Aut(X)$ is trivial. This shows that the group $\Aut(X)$ is a discrete topological group and the kernel of the natural representation \eqref{aut} of $\Aut(X)$ in $\Or(S_X)$ is a finite group.

\begin{remark} One can say more about the kernel $H$ of homomorphism \eqref{aut} (see \cite{Ni2}, \S 10). Let $\chi:\Aut(X)\to \bbC^*$ be the one-dimenensional representation of $\Aut(X)$ in the space $\Omega^2(X)$ of holomorphic 2-forms on $X$. The image of $\chi$ is a cyclic group of some order $n$. First Nikulin proves that the value of the Euler function $\phi(n)$ divides $22-\rank S_X$. Next he proves that  the restriction of $\chi$ to $H$ is injective. This implies that $H$ is a cyclic group of order dividing $n$. All possible values of $n$ which can occur are known (see \cite {KondoNew}). The largest one is equal to $ 66$ and can be  realized for a K3 surface birationally isomorphic to a surface in the weighted projective space $\bbP(1,6,22,33)$ given by the equation
$x^{66}+y^{11}+z^3+w^2 = 0.$
\end{remark}

It follows from the adjunction formula that any smooth rational curve on a K3 surface is a $(-2)$-curve. The class of this curve in $S_X$ defines a reflection. Let $W_X^+$ denote the subgroup of $\Or(S_X)$ generated by these reflections. 

\begin{proposition} 
$$\Ref_{-2}(S_X) = W_X^+.$$
\end{proposition}

\begin{proof} Let $C$ be an irreducible curve on $X$. By adjunction formula, $[C]^2 \ge -2$ and $C^2 = -2$ if and only if $C$ is a $(-2)$-curve. We call a divisor class \emph{effective} if it can be represented by  a (possibly reducible) algebraic curve on $X$. Using Riemann-Roch Theorem, one shows that any divisor class $x$ with $x^2\ge -2$ is either effective or its negative is effective. Let $D = \sum_{i\in I}C_i$ be any algebraic curve written as a sum of its irreducible components. We have $[D]\cdot [C] \ge 0$ for any irreducible curve $C$ unless $C$ is an irreducible component of $D$ with $C^2 = -2$. A divisor class $x$ is called \emph{nef} if $x\cdot d\ge 0$ for any effective divisor class $d$.  It is not difficult to show that  the $W_X^+$-orbit of any effective divisor $x$ with $x^2\ge 0$ contains a unique nef divisor (use that $x\cdot e < 0$ for some effective $e$ with $e^2 = -2$ implies $r_e(x)\cdot e > 0$).  Let $V_X^0 = \{x\in V_X: x^2 > 0\}$.  We take for the model of the hyperbolic space $H^n$ associated with $S_X$ the connected component of $V_X^0/\bbR_+ \subset \bbP(V_X)$ which contains the images of effective divisors from $V_X^0$. Then the image $P^+$ in $H^n$ of the convex hull $N$ of nef effective divisors from $V_X^0$  is a fundamental polytope for the reflection group $W_X^+$. Its fundamental roots are the classes of $(-2)$-curves. Let $r_\alpha$ be a reflection from $ \Ref_{-2}(S_X)$. Replacing $\alpha$ with $-\alpha$ we may assume that $\alpha$  is  effective. Suppose $\alpha$ is a fundamental root for a fundamental polytope $P$ of $\Ref_{-2}(S_X)$ which is contained in $P^+$. Since all vectors from $N$ satisfy $x\cdot \alpha \ge 0$, we see that $P^+\subset P$ and hence $P = P^+$. This shows that $W_X^+$ and $\Ref_{-2}(S_X)$ are defined by the same convex polytope and hence the  groups are equal.
\end{proof}

The following result follows from the fundamental \emph{Global Torelli Theorem for K3 surfaces} due to I. Shafarevich and I. Pjatetsky-Shapiro \cite{PS}.

\begin{theorem} Let $A_X$ be the image of $\Aut(X)$ in $\Or(S_X)$. Let $G$ be the subgroup of $\Or(S_X)$ generated by $W_X^+$ and $A_X$. Then $G = W_X^+\rtimes A_X$ and its index in $\Or(S_X)$ is finite.
\end{theorem}

The  difficult part is the finiteness of the index. 

\begin{corollary} The following assertions are equivalent:
\begin{itemize}
\item $W_X^+$ is of finite index in $\Or(S_X)$;
\item $\Aut(X)$ is a finite group;
\item the $(-2)$-reflection group of $S_X$ admits a fundamental polytope with finitely many faces defined by the classes of smooth rational curves.
\item $S_X(-1)$ is a 2-reflective lattice.
\end{itemize}
\end{corollary}

The third property implies that the set of $(-2)$-curves is finite if $\Aut(X)$ is finite but the converse is not true.
 
  The following  are  even 2-reflective lattices of  rank $\ge 17$, the Coxeter diagrams of the reflection groups $\Ref_{-2}(S_X)$, and the corresponding K3 surfaces.
\begin{figure}[h]
\begin{multicols}{2}{
\xy
(-20,0)*{};
@={(0,0),(5,0),(10,0),(15,0),(20,0),(25,0),(30,0),(35,0),(40,0),(45,0),(50,0),(0,5),(0,-5),(0,-10),(0,-15),
(50,5),(50,-5),(50,-10),(50,-15)}@@{*{\bullet}};
(0,0)*{};(50,0)*{}**\dir{-};(0,5)*{};(0,-15)*{}**\dir{-};(50,5)*{};(50,-15)*{}**\dir{-};
(0,-15)*{};(50,-15)*{}**\dir2{-};
(0,-15.05)*{};(50,-15.05)*{}**\dir2{-};
(0,-15.1)*{};(50,-15.1)*{}**\dir{-};
(0,-14.95)*{};(50,-14.95)*{}**\dir2{-};
(0,-14.9)*{};(50,-14.9)*{}**\dir{-};
(25,-20)*{S_X= U\perp E_8(-1)\perp E_7(-1)};
\endxy
\xy
(-20,0)*{};
(10,7.2)*\cir<20pt>{};
(10,7.2)*\cir<23pt>{};
(10,7.2)*!<-.2pt,.2pt>\cir<20pt>{};
(10,7.2)*!<-.4pt,.4pt>\cir<20pt>{};
(10,7.2)*!<-.6pt,.6pt>\cir<20pt>{};
(10,7.2)*!<-.8pt,.8pt>\cir<20pt>{};
(0,0)*{};(40,0)*{}**\dir{-};
(0,10)*{};(0,-10)*{}
**\crv{(5,4)&(20,2)&(25,1)&(28,.5)&(35,0)&(28,-.5)&(25,-1)&(20,-2)&(5,-4)}
\endxy}
\end{multicols}
\caption{}
\end{figure}

\begin{figure}[h]
\begin{multicols}{2}{
\xy 
@={(-5,0),(0,0),(5,0),(10,0),(15,0),(20,0),(25,0),(30,0),(35,0),(40,0),(45,0),(50,0),(55,0),(60,0),
(65,0),(70,0),(75,0)}@@{*{\bullet}};
(-5,0)*{};(75,0)*{}**\dir{-};
(5,0)*{};(5,-5)*{\bullet}**\dir{-};(65,0)*{};(65,-5)*{\bullet}**\dir{-};
(30,-15)*{S_X=U\perp E_8(-1)^2};
\endxy
\xy
(-25,0)*{};(50,-10)*{};
(.5,15)*{};(20,-10)*{}**\dir2{-};
(.5,15)*{};(20,-10)*{}**\dir{-};
(.5,15.2)*{};(20,-9.8)*{}**\dir{-};
(.5,14.8)*{};(20,-10.2)*{}**\dir{-};
(.5,-15)*{};(20,10)*{}**\dir2{-};
(.5,-15)*{};(20,10)*{}**\dir{-};
(.5,-14.8)*{};(20,10.2)*{}**\dir{-};
(.5,-15.2)*{};(20,9.8)*{}**\dir{-};
(0,0)*{};(35,0)*{}**\dir{-};
(0,18)*{};(0,-18)*{}
**\crv{(5,1)&(20,2)&(30,.2)&(34,0)&(30,-.2)&(20,-2)&(5,-1)}
\endxy}
\end{multicols}
\caption{}
\end{figure}
\begin{figure}[h]
\begin{multicols}{2}{
\xy @={(0,0),(5,0),(10,0),(15,0),(20,0),(25,0),(30,0),(35,0),(40,0),(45,0),(50,0),(55,0),(60,0),
(0,5),(0,-5),(60,5),(60,-5),(30,-5),(30,-10)}@@{*{\bullet}};
(0,0)*{};(60,0)*{}**\dir{-};
(0,5)*{};(0,-10)*{\bullet}**\dir{-};(60,5)*{};(60,-10)*{\bullet}**\dir{-};
(30,0)*{};(30,-15)*{\bullet}**\dir{-};
(30.05,-5)*{};(30.05,-15)*{\bullet}**\dir2{-};
(29.95,-5)*{};(29.95,-15)*{\bullet}**\dir2{-};
(30.1,-5)*{};(30.1,-15)*{\bullet}**\dir{-};
(29.9,-5)*{};(29.9,-15)*{\bullet}**\dir{-};
(.5,-10)*{};(30,-15)*{\bullet}**\dir2{-};(60,-10)*{};(30,-15)*{\bullet}**\dir2{-};
(.5,-10.1)*{};(30,-15.1)*{\bullet}**\dir{-};(60,-10.1)*{};(30,-15.1)*{\bullet}**\dir{-};
(.5,-9.9)*{};(30,-14.9)*{\bullet}**\dir{-};(60,-9.9)*{};(30,-14.9)*{\bullet}**\dir{-};
(30,-20)*{S_X = U\perp E_8^{2}\perp A_1};
(-15,0)*{};
\endxy

\xy 
(-20,0)*{};
(0,13)*{};(21,-10)*{}**\dir2{-};
(0,13)*{};(21,-10)*{}**\dir{-};
(0,13.2)*{};(21,-9.8)*{}**\dir{-};
(0,12.8)*{};(21,-10.2)*{}**\dir{-};
(.5,-15)*{};(20,10)*{}**\dir1{-};
(.5,-15)*{};(20,10)*{}**\dir2{-};
(.5,-14.8)*{};(20,10.2)*{}**\dir{-};
(.5,-15.2)*{};(20,9.8)*{}**\dir{-};
(0,0)*{};(35,0)*{}**\dir{-};
(0,18)*{};(0,-18)*{}
**\crv{(5,2)&(20,2)&(25,1)&(28,.5)&(30,.2)&(33,0)&(30,-.1)&(28,-.5)&(25,-1)&(20,-2)&(5,-1)}
\endxy}
\end{multicols}
\caption{}
\end{figure}

The K3 surface is birationally isomorphic the surface obtained as the double cover of $\bbP^2$ branched along the  curve of degree 2 drawn thick, followed by the double cover along the proper inverse transform of the  remaining curve of degree 4 (the cuspidal cubic and its cuspidal tangent line).

 It follows from the classification of 2-reflective lattices (see section 4) that each of them is isomorphic to the lattice $S_X$ for some K3-surface $X$. 

Similar assertion for any even reflective lattice is not true for the following trivial reason. Replacing $M$ by $M(2k)$ for some $k$ we  obtain a reflective lattice with discriminant group whose  minimal number of generators $s\ge \rank M$. Suppose $\rank M \ge 12$ and $M(k)$ is primitively embedded in $L_{K3}$. Its orthogonal complement is a lattice of rank $22-12\le 10$ with isomorphic discriminant group generated by $\le 10$ elements. This is a contradiction.

A more serious reason is the existence of an even  reflective hyperbolic lattice $M$ of rank 22. One can take $M = U\perp D_{20}$ realized as the sublattice of $I_{19,1}$ of vectors with even $v^2$.
The reflectivity of $M$ was proven by R. Borcherds \cite{Bo1}. Since $\rank~S_X\le 20$ for any complex K3 surface $X$, the lattice $M$ cannot be isomorphic to $S_X$. However,  the lattice $M$ is realized as the Picard lattice of a K3 surface over an algebraically closed field of characteristic 2 isomorphic to a quartic surface in $\bbP^3$ with equation
$$T_0^4+T_1^4+T_2^4+T_3^4+T_0^2T_1^2+T_0^2T_2^2+T_1^2T_2^2+T_0T_1T_2(T_0+T_1+T_2)=0.$$
(see \cite{DK}). Note that over a field of positive characteristic the Hodge structure and  the inequality $\rho = \rank S_X\le 20$ does not hold. However, one can show that for any K3 surface over an algebraically closed field of positive characteristic
\begin{equation}\label{ineq1}
\rho \le 22, \quad \rho \ne 21\end{equation}
K3 surfaces with $\rho = 22$ are called \emph{supersingular} (in sense of Shioda).
Observe the striking analogy of inequalities \eqref{ineq1}  with the inequalities from Theorem \ref{ess}.

Besides  scaling, one can consider the following operations over nondegenerate quadratic lattices which preserve the reflectivity property (see \cite{SW}). The first operation replaces a lattice $M$ with $p^{-1}(M\cap p^2M^*)+M$ for any $p$ dividing the discriminant of $M$. This allows one to replace $M$ with a lattice such that the exponent of the  discriminant groups is square free. The second operation replaces $M$ with $N(p)$, where $N = M^*\cap p^{-1}M$. This allows one to replace $M$ with a lattice such that the largest power $a$ of $p$ dividing the discriminant of $M$ satisfies  $a\le \frac{1}{2}\rank M$. 

I conjecture that up to scaling and the above two  operations any even hyperbolic reflective lattice is isomorphic to the lattice $S_X(-1)$ for some K3 surface defined over an algebraically closed field of characteristic $p \ge 0$.
 
It is known that the lattice $S_X(-1)$ for a supersingular K3-surface $X$ over a field of characteristic $p >0$  is always of rank 22 and its discriminant group is isomorphic to a $p$-elementary group  $(\bbZ/p\bbZ)^{2\sigma}, \sigma \le 10,$ (see \cite{RS}). No two such lattices are equivalent in the sense of the operations on lattices described in above. There is only one such reflective lattice, namely $U\perp D_{20}(-1)$.

There are only a few cases where one can compute explicitly the automorphism group of a K3 surface when it is infinite and the rank of $S_X$ is large. This requires to construct explicitly a Coxeter polytope of  $\Ref_2(S_X(-1))$ which is of infinite volume. As far as I know this has been accomplished only in the following cases
\begin{itemize}
\item $X$ is the Kummer surface of the Jacobian variety of a general curve of genus 2 (\cite{Kondo});
\item $X$ is the Kummer surface of the product of two non-isogeneous elliptic curves (\cite {KK});
\item $X$ is birationally isomorphic to the  Hessian surface of a general cubic surface (\cite{DoKe});
\item $S_X(-1)= U\perp E_8^{ 2}\perp A_2$ (\cite{Vi4});
\item $S_X(-1) = U\perp E_8^{2}\perp A_1^{ 2}$ (\cite{Vi4});
\item $S_X$ is of rank $20$ with discriminant 7 (\cite{Bo3});
\item $S_X(-1) = U\perp D_{20}$ (\cite {DK}) (characteristic 2).
\end{itemize}

What is common about these examples is that the lattice $S_X(-1)$ can be primitively embedded in the lattice $II_{25,1}$ from Example \ref{leech} as an orthogonal sublattice to a finite root sublattice of $II_{25,1}$. We refer to \cite{Bo3} for the most general method to derive from this an  information about $\Aut(X)$. 

\subsection{Enriques surfaces} These are the surfaces from Case 2 (ii). They satisfy
$$2c_1(X) = 0, \quad  c_1(X) \ne 0, \quad H_1(X,\bbZ) = \Tors(H^2(X,\bbZ)) =\bbZ/2\bbZ.$$
The cover $Y\to X$ corresponding the generator of $H_1(X,\bbZ)$ is a K3 surface. So, Enriques surfaces correspond to pairs $(Y,\tau)$, where $Y$ is a K3 surface and $\tau$ is its fixed-point-free involution. 

We have
$$S_X(-1) = U\perp E_8 \cong E_{10}.$$
The reflection group of  $E_{10}$ is the group $W(2,3,7)$. Since its fundamental polytope does not have non-trivial symmetries we obtain from \eqref{semi}
\begin{equation}
\Or(S_X)^+ = \Ref(S_X(-1)) \cong W(2,3,7).\end{equation}
Let $W_X^+$ be the subgroup of $\Or(S_X)$ generated by reflections $r_\alpha$, where $\alpha$ is the image in $S_X$ of the cohomology class  of a $(-2)$-curve on $X$. 
The following theorem follows (but non-trivially) from  the Global Torelli Theorem for K3 surfaces (\cite \cite{Na}, {Ni2}).

\begin{theorem} Let $A_X$ be the image of $\Aut(X)$ in $\Or(S_X)$. Then $A_X\subset \Or(S_X)^+$, its intersection with $W_X^+$ is trivial and 
$W_X^+\rtimes A_X$ is of finite index in $\Or(S_X)$.
\end{theorem}

This gives as a corollary that $\Aut(X)$ is finite if and only $W_X^+$ is of finite index in $\Or(S_X)$. A general Enriques surface (in some precise meaning) does not contain $(-2)$-curves so $A(X)$  is isomorphic to a subgroup of finite index of the orthogonal group $\Or(E_{10})$. In fact, more precisely,  the group $\Aut(X)$ is isomorphic  to  the 2-level congruence subgroup of $W(E_{10})$ defined in \eqref{cs} (see \cite{BP}, \cite{Ni2}). The fact that the automorphism groups of a general Enriques surface and a general Coble surface are isomorphic  is not a coincidence, but I am not going to explain it here (see \cite{Do2}).  

All Enriques surfaces $X$ with finite $\Aut(X)$ were classified and explicitly constructed  by S. Kondo \cite{Kond} and (not constructively) by V. Nikulin \cite{Ni3}. In \cite{Do} I gave an example of an Enriques surface with finite automorphism group, believing that it was the first example of this kind.  After the paper had been published I  found that the existence of another example  was claimed much earlier  by G. Fano \cite{Fano}. However his arguments are very obscure and impossible to follow.   The Coxeter diagram of the reflection group $W_X^+$ in my example is given in Figure \ref{fig6}.
\vspace{15pt}

\begin{figure} [ht]
\xy
(-65,0)*{};
@={(-20,0),(20,0),(0,20),(0,-20),(14.15,14.15),(-14.15,-14.15),
(-14.15,14.15), (14.15,-14.15),(12,0),(-12,0),(4,0),(-4,0)}@@{*{\bullet}};
(-20,0)*{};(-18,7.4)*{}**\dir{-},(-18,7.4)*{};(-14.15,14.15)*{}**\dir{-};
(-14.15,14.15)*{};(-7.4,18)*{}**\dir{-},(-7.4,18)*{};(0,20)*{}**\dir{-};
(20,0)*{};(18,7.4)*{}**\dir{-},(18,7.4)*{};(14.15,14.15)*{}**\dir{-};
(14.15,14.15)*{};(7.4,18)*{}**\dir{-},(7.4,18)*{};(0,20)*{}**\dir{-};
(-20,0)*{};(-18,-7.4)*{}**\dir{-},(-18,-7.4)*{};(-14.15,-14.15)*{}**\dir{-};
(-14.15,-14.15)*{};(-7.4,-18)*{}**\dir{-},(-7.4,-18)*{};(0,-20)*{}**\dir{-};
(20,0)*{};(18,-7.4)*{}**\dir{-},(18,-7.4)*{};(14.15,-14.15)*{}**\dir{-};
(14.15,-14.15)*{};(7.4,-18)*{}**\dir{-},(7.4,-18)*{};(0,-20)*{}**\dir{-};
(-20,0)*{};(20,0)*{}**\dir{-};
(20,0)*{};(10,0)*{}**\dir{-};
(-12,.05)*{};(12,.05)*{}**\dir{-};
(-12,-.05)*{};(12,-.05)*{}**\dir{-};
(-12,.1)*{};(12,.1)*{}**\dir{-};
(-12,-.1)*{};(12,-.1)*{}**\dir{-};
(-12,0)*{};(12,0)*{}**\dir2{-};

\endxy

\caption{}\label{fig6}

\end{figure}

 \section{Cremona transformations}
\subsection{Plane Cremona transformations} A \emph{Cremona transformation} of projective space $\bbP^n$ (as always over complex numbers)  is a birational map  of algebraic varieties. It can be given in projective coordinates by $n+1$ homogeneous polynomials of the same degree $d$:
\begin{equation}\label{cr}
T: (t_0,\ldots,t_n) \mapsto (P_0(t_0,\ldots,t_n),\ldots, P_n(t_0,\ldots,t_n)).
\end{equation}
Dividing by a common multiple of the polynomials we may assume that the map is not defined on a closed subset of codimension $\ge 2$, the set of common zeros of the polynomials $P_0,\ldots,P_n$. 
Let $U = \textup{dom}(T)$ be the largest open subset where $T$ is defined and let $X$ be the Zariski closure of the graph of $T:U\to \bbP^n$ in $\bbP^n\times \bbP^n$.  By considering the two projections of $X$, we get a commutative diagram  of birational maps 
\begin{equation}\label{diag}
\xymatrix{&X{}\ar[ld]_\pi\ar[rd]^\sigma&\\
\bbP^2\ar^T@{-->}[rr]&&\bbP^2}
\end{equation}

By taking a resolution of singularities $X'\to X$, we may assume that such diagram exists with $X$ nonsingular. We  call it a \emph{resolution of indeterminacy points} of $T$. 

We denote the projective space of homogeneous polynomials  of degree $d$ in $n+1$ variables by $|\calO_{\bbP^n}(d)|$. The projective subspace of dimension $n$ spanned by the polynomials $P_0,\ldots,P_n$ defining the map $T$ is denoted by $L(T)$ and is called the \emph{linear system} defining $T$. It  depends only on $T$ (recall that two birational maps are equal if they coincide on an Zariski open subset). Any $n$-dimensional projective subspace $L$ of $|\calO_{\bbP^n}(d)|$ (linear system) which consists of polynomials without common factor defines a rational map $\bbP^n-\to \bbP^n$ by simply choosing a linear independent ordered set of $n+1$ polynomials from $L$.  When the map happens to be birational, the linear system is called \emph{homaloidal}. One obtains $L(T)$ as follows. First one considers the linear system $|\calO_{\bbP^n}(1)|$ of hyperplanes in the target $\bbP^n$. The pre-image of its member on $X$ under the map $\sigma$ is a hypersurface on $X$, we push it down by $\pi$ and get a hypersurface on the domain $\bbP^n$. The set of such hypersurfaces forms  a linear system $L(T)$.

Now let us assume that $n = 2$. As in section \ref{rat} we consider a factorization \eqref{decom} of $\pi$ 
\begin{equation}
\pi:X = X_N\overset{\pi_N}{\longrightarrow} X_{N-1}\overset{\pi_{N-1}}{\longrightarrow} \ldots \overset{\pi_2}{\longrightarrow} X_1\overset{\pi_1}{\longrightarrow} X_0 = \bbP^2,
\end{equation}
where $\pi_{k+1}$ is the blow-up of a point $x_k\in X_k, k = 0,\ldots,N-1$. Recall that the \emph{multiplicity} $\mult_x(D)$ of a hypersurface $D$ at a point $x$ on a nonsingular variety is the degree of the first nonzero homogeneous part in the Taylor expansion of its local equation at $x$. Define inductively the numbers $m_i$   as follows. Let $L = L(T)$ be the linear system of curves on $\bbP^2$ defining $T$. First we set
$$m_1 = \min_{D\in L}\mult_{x_1}D$$
The linear system $\pi_1^*(L)$ on $X_1$ which consists of the full pre-images of hypersurfaces from $L(T)$ on $X_1$ has the hypersurface $m_1E_1$ as a fixed component. Let
$$  L_1 = \pi_1^*(L)- m_1E_1.$$
This is a linear system on $X_1$ without fixed components.
Suppose $m_1,\ldots,m_i$ and $L_1,\ldots,L_i$  have been defined. Then we set
$$m_{i+1} = \min_{D\in L_i}\mult_{x_{i+1}}D,$$
$$  L_{i+1} = \pi_{i+1}^*(L_i)-m_iE_{i+1}.$$
It follows from the definition that 
$$L_N = \pi^*(L)-\sum_{i=1}^Nm_i\calE_i$$
has no   fixed components and is equal to the pre-image of $|\calO_{\bbP^2}(1)|$ under $\sigma$. The image of $L_N$ in $\bbP^2$ is equal to the linear system $L(T)$. It  is denoted by 
$$|\calO_{\bbP^2}(d)-m_1x_1-\ldots-m_Nx_N|.$$
The meaning of the notation is that $L(T)$ consists of plane curves of degree $d$ which pass through the points $x_i$ with multiplicities $\ge m_i$.

We have similar decomposition for the map $\sigma$ which defines the linear system $L(T^{-1})$.

As we have shown in section \ref{rat},  two factorizations of birational regular maps from $X$ to $\bbP^2$  in a sequence of blow-ups define two geometric  bases of the lattice  $S_X$. 
A Cremona transformation \eqref{cr} together with a choice of a diagram \eqref{diag} and the factorizations \eqref{decom} for $\pi$ and $\sigma$ is called a \emph{marked Cremona transformation}. It follows from the proof of Theorem \ref{crem}  that any marked Cremona transformation defines an element of the Coxeter group $W(E_N)$. The corresponding matrix is called the \emph{characteristic matrix} of the marked Cremona transformation.

\begin{example}\label{qc} Let $T$ be the \emph{standard quadratic Cremona transformation}  defined by 
$$T:(x_0,x_1,x_2)\mapsto (x_0^{-1},x_1^{-1},x_2^{-1}).$$
(to make sense of this one has to multiply all coordinates at the output by $x_0x_1x_2$).  It is not defined at points $p_1 = (1,0,0), p_2 = (0,1,0), p_3 = (0,0,1)$. Let $\pi:X= X_3\to X_2\to X_1 \to X_0= \bbP^2$ be the
composition of the blow-up of $p_1$, then the blow-up of the pre-image of $p_2$, and finally the pre-image of $p_3$. It is easy to see on the open subset $U$ where $T$ is defined the coordinate line $t_i = 0$ is mapped to the point $p_i$. This implies that the Zariski closures on $X$ of the pre-images of  intersections of these lines with $U$ are the curves $E_1',E_2',E_3'$ on $X$ which are blown down to the points $p_1,p_2,p_3$ under $T$. The factorization for $\sigma:X\to \bbP^2$ could be chosen in such a way that $E_i' = \calE_i'$ are the exceptional curves. These curves define a new geometric basis in $X$ with 
$$e_1'= e_0-e_2-e_3, \quad e_1'= e_0-e_1-e_3, \quad e_1'= e_0-e_0-e_1.$$
Since $k_X = -3e_0+e_1+e_2=e_3 = -3e_0'+e_1'+e_2'+e_3',$ we also get
$e_0' = 2e_0'-e_1'-e_2'-e_3'.$ The corresponding transformation is the reflection with respect to the vector $e_0-e_1-e_2-e_3$.
Given any set of 3 non-collinear points $q_1,q_2,q_3$, one can find a quadratic\footnote{i.e. defined by polynomials of degree 2} Cremona transformation $T'$ with indeterminacy points $q_1,q_2,q_3$. For this we choose a projective transformation $g$ which sends $q_i$ to $p_i$ and take $T' = T\circ g$. 
\end{example}

\begin{theorem}\label{charr} Let $A$ be the matrix representing an element from $W(E_N)$. Then there exists a marked Cremona transformation whose characteristic matrix is equal to $A$.
\end{theorem}

\begin{proof} Let $A$ be the matrix of $w\in W(E_N)$ with respect to the standard basis 
$\bfe_0,\ldots, \bfe_N$ of $I_{1,N}$. Its first column is a vector $(m_0,-m_1,\ldots,-m_N)$. Write $w$ as a word in reflections $r_{\bfa_i}$ and use induction on the length of $w$ to prove the Noether inequalities 
$m_i\ge 0, i\ge 1$ and $m_0 > 3\max\{m_i,i\ge 1\}$. Now, let us use induction on the length of $w$ to show that the linear system $|\calO_{\bbP^2}(m_0)-m_1p_1-\ldots-m_Np_N)|$ is homaloidal for some points $p_1,\ldots, p_N$ in general position. If $w$ is a simple reflection $r_{\bfa_i}$ we get $w(\bfe_0) = \bfe_0$ or $2\bfe_0-\bfe_1-\bfe_2-\bfe_3$. In the first case the linear system defines a projective transformation and in the second case it defines a standard quadratic  transformation if we choose three non-colinear points $p_1,p_2,p_3$. Now write $w = r_{\bfa_i}w'$ where the length of $w'$ is less than the length of $w$. By induction, $w'$ defines 
a homaloidal linear system $|\calO_{\bbP^2}(m_0')-m_1'p_1-\ldots-m_N'p_N)|$ for some points $p_1,\ldots,p_N$ in general position.  Let $\Phi':\bbP^2- \to \bbP^2$ be the corresponding Cremona transformation. If $i\ne 0$, the reflection permutes the $m_i', i > 0,$ and the linear system is still homaloidal and the set of points does not change. Since the points are in general position we may assume that $p_1,p_2,p_3$ are distinct non-colinear points. Composing $\Phi'$ with a projective transformation we may assume that  $p_1 = (1,0,0), p_2 = (0,1,0), p_3 = (0,0,1)$. The computation from the previous example shows that 
$r_{\bfa_0}w'(e_0) = m_0\bfe_0-\sum_{i\ge 1} m_i\bfe_i$, where
$m_i = 2m_i'-m_1-m_2-m_3, i \le 3, $ and $m_i = m_i', i > 3$. Let $\Phi =  \Phi'\circ T$, where $T$ is the standard quadratic transformation discussed in the previous example. Then it is easy to see that $\Phi$ is given by the linear system $|\calO_{\bbP^2}(m_0)-m_1p_1-\ldots-m_Np_N)|$.  Let $\tilde{w}\in W(E_N)$ correspond to the characteristic matrix of $\Phi$ (with respect to an appropriate marking). We have proved that $w(\bfe_0) = \tilde{w}(\bfe_0)$.  This implies that $\tilde{w}w^{-1}(\bfe_0) = \bfe_0$. The matrix representing an element from $\Or(1,N)$ whose first column is the unit vector is a diagonal matrix with $\pm 1$ at the diagonal. As we have seen already in the proof of  Theorem \ref{crem} this implies that $\tilde{w} = w$.

\end{proof}


One can apply Theorem \ref{charr} to list the types $(m_0,m_1,\ldots,m_N)$ of all homaloidal linear systems with $N$ indeterminacy points. They correspond to the orbit of the vector $e_0$ with respect to the group $W(E_N)$. In particular, the number of types is finite only for $N\le 8$.

\subsection{Cremona action of $W(p,q,r)$}

Consider the natural diagonal action of the group $G = \PGL(n+1,\bbC)$ on  $(\bbP^n)^{N},$ where
$m= N-n-2\ge 0$. A general orbit contains a unique point set $(p_1,\ldots,p_N)$ with the first $n+2$ points equal to the set of \emph{reference points} $(1,0,\ldots,0),\ldots,(0,\ldots,0,1),(1,\ldots,1)$. This easily implies that the field of $G$-invariant rational functions on $(\bbP^n)^N$ is isomorphic to the field of rational functions on $(\bbP^n)^{m}$ and hence is isomorphic to the field of rational functions  $\bbC(z_1,\ldots,z_{nm})$.
The symmetric group $\Sigma_{N}$ acts naturally on this field via its action on $(\bbP^n)^N$ by permuting the factors. Assume $n\ge 2$ and consider $\Sigma_N$ as a subgroup $W(1,n+2,m+1)$ of the Coxeter group of type $W(2,n+1,m+1)$ corresponding to the subdiagram of type $A_{N-1}$ of the Coxeter diagram of $W(2,n+1,m+1)$. In 1917 A. Coble  extended the action of $\Sigma_N$ on the field $\bbC(z_1,\ldots,z_{nm})$ to the action of the whole group $W(2,n+1,m+1)$. This construction is explained in modern terms in \cite{DO}. In Coble's action the remaining generator of the Coxeter group acts as a standard quadratic transformation $\bbP^n-\to \bbP^n$ defined  by 
$$T:(x_0,\ldots,x_n)\mapsto (x_0^{-1},\ldots,x_n^{-1}).$$
 One takes a point set $(p_1,\ldots,p_{N})$, where the first $n+2$ points are the reference points. Then  applies $T$ to the remaining points to get a new set 
$$(p_1,\ldots,p_{n+1}, T(p_{n+2}),\ldots,T(p_N)).$$ The \emph{Cremona action} is the corresponding homomorphism of groups
$$W(2,n+1,m+1)\to \Aut_\bbC(\bbC(z_1,\ldots,z_{nm})).$$
One can show  that for $N \ge 9$, this homomorphism does not arise from a regular action of the Coxeter group on any Zariski open subset of $(\bbP^n)^{m}$. 

The following result of  S. Mukai \cite{Mu} extends the Cremona action to any  group $W(p,q,r)$.

\begin{theorem}  Let $X = (\bbP^{q-1})^{p-1}$. Consider the natural diagonal action of the group $\PGL(q,\bbC)^{p-1}$ on $X$ and extend it to the  diagonal action  on 
$X_{p,q,r}: = X^{q+r}$. Let $K(p,q,r)$ be the field of invariant rational functions  isomorphic to $\bbC(t_1,\ldots,t_d), d = (p-1)(q-1)(r-1)$. Then there is a natural homomorphism 
$$\text{cr}_{p,q,r}:W(p,q,r) \to \textup{Aut}_\bbC(K(p,q,r)).$$
It coincides with the Coble action when $p = 2$. 
\end{theorem}

It seems that  the homomorphism $\text{cr}_{p,q,r}$ is always injective when the group is infinite. At least it is true in the case when one of the numbers $p,q,r$ is equal to 2.  The geometric meaning of the kernels in the  case of  finite groups $W(2,q,r)$ are discussed in \cite{DO} and \cite{DV2}.

\smallskip

The reflections corresponding to the vertices on the branches of the $T_{p,q,r}$-diagram with $q$ and $r$ vertices act by permuting the factors of $X$. The reflections corresponding to $p-2$ last vertices of the $p$-branch  permute the factors of $X$. The second vertex on the $p$-branch acts via the Cremona transformation in $(\bbP^{q-1})^{p-1}$  
{\small$$\bigl((x_0^{(1)},\ldots,x_{q-1}^{(1)}),\ldots, (x_0^{(p-1)},\ldots,x_{q-1}^{(p-1)})\bigr)\to 
\bigl((\frac{1}{x_0^{(1)}},\ldots,\frac{1}{x_{q-1}^{(1)}}),\ldots, (\frac{x_0^{(p-1)}}{x_{0}^{(1)}},\ldots,\frac{x_{q-1}^{(p-1)}}{x_{q-1}^{(1)}})\bigr).$$}

\smallskip
Let $Y(p,q,r)$ be a birational model of the field $K(p,q,r)$ on which $W(p,q,r)$ acts birationally via $c_{pqr}$. For any $g\in W(p,q,r)$ let $\textup{dom}(g)$ be the domain of definition of $g$. Let $Z$ be a closed irreducible subset of $Y(p,q,r)$ with generic point $\eta_Z$  Let 
$$G_Z = \{g\in W(p,q,r):\eta_Z\in \textup{dom}(g)\cap \textup{dom}(g^{-1}), g(\eta_Z) = \eta_Z\}$$
be the \emph{decomposition subgroup} of $Z$ in $W(p,q,r)$ and $G_Z^i$ be the \emph{inertia subgroup} of $Z$, the kernel of the natural map 
$G_Z \to \Aut(R(Z))$, where $R(Z)$ is the field of rational functions of $Z$. These groups were introduced  for  any group of birational transformations by M. Gizatullin \cite{Gi2}.  Define a $\PGL(q,\bbC)^{p-1}$-invariant closed irreducible subset $S$ of $X_{p,q,r}$ to be \emph{special} if it defines a closed subset $Z$ on some birational model $Y(p,q,r)$  such that  $G_Z = W(p,q,r)$ and $G_Z^i$ is a subgroup of finite index of $W(p,q,r)$. 

Consider a general point $s\in S$ as a set of $q+r$ distinct points in $(\bbP^{q-1})^{p-1}$ and let $V(s)\to (\bbP^{q-1})^{p-1}$ be the blow-up of this set.  One can show that  the group $G_Z'$ is realized as a group of  pseudo-automorphisms of $V(s)$.\footnote{A pseudo-automorphism is a birational transformation which is an isomorphism outside a closed subset of  codimension $>2$.} 

I know only a few examples of special subsets when $W(p,q,r)$ is infinite.  Here are some examples. 
\begin{itemize}
\item $(p,q,r) = (2,3,6)$, $S$ parametrizes 
 ordered sets of base points of a  pencil of plane cubic curves  in $\bbP^2$\cite{Coble}, \cite{DO};
\item $(p,q,r) = (2,4, 4)$, $S$ parametrizes ordered sets of base points  of a  net of quadrics in $\bbP^3$ \cite{Coble}, \cite{DO};
\item $(p,q,r) = (2,3,7)$,  $S$ parametrizes ordered sets of double points of a rational plane sextic \cite{Coble}, \cite{DO};
\item $(p,q,r) = (2,4,6)$,  $S$ parametrizes ordered sets of double points of a quartic symmetroid surface \cite{Coble},\cite{CD}.

\end{itemize}

The inertia subgroups  of finite index of $W(p,q,r)$ defined in these examples have the quotient groups  isomorphic to simple groups $\Or(8,\bbF_2)^+, \Sp(6,\bbF_2), \Or(10,\bbF_2)^+$ and  $\Sp(8,\bbF_2)$, respectively. 
 
\begin{remark} It is popular in group theory to represent a sporadic simple group or a related group as a finite quotient of a Coxeter group $W(p,q,r)$. For example, the Monster group $F_1$ is a quotient of $W(4,5,5)$. The Bimonster group $F_1\wr 2$ is a quotient of $W(6,6,6)$ by a single relation \cite{Iv},\cite{No}. Is there a geometric interpretation of these presentations in terms of the Cremona action of $W(p,q,r)$ on some special subset of points in $X_{p,q,r}$?  Mukai's construction should relate the Monster group with some special configurations of 10 points in $(\bbP^3)^4$ or $9$ points in $(\bbP^4)^4$.  The Bimonster group could  be related to  special configurations of  $12$ points in $(\bbP^5)^5$ (see  related speculations in \cite{All}).
\end{remark}

\section{Invariants of finite complex reflection groups}

Let $\Gamma \subset \GL(n+1,\bbC)$ be a finite linear complex reflection group in $\bbC^{n+1}$ and let $\bar{\Gamma}$ be its image in $\PGL(n,\bbC)$. The reflection hyperplanes of $G$ define a set of hyperplanes in $\bbP^n$ and the zeroes of $G$-invariant polynomials define hypersurfaces in $\bbP^n$. The geometry, algebra, combinatorics and topology of arrangements of reflecting hyperplanes of  finite complex reflection groups is a popular area in the theory of hyperplane arrangements (see \cite{OS}, \cite{OS2}). On the other hand, classical algebraic geometry is full of interesting examples of projective hypersurfaces whose symmetries are described in terms of a complex reflection group. We discuss only a few examples. 

We begin with the group $J_3(4)$ of order 336 (No  24 in the list). It has fundamental invariants of degrees $4,6$ and $14$. Its center is of order 2 and the group $\bar{\Gamma}$ is a simple group of order 168 isomorphic to $\PSL(2,\bbF_7)$. The invariant curve of degree 4 is of course the famous Klein quartic which is projectively equivalent to the curve
$$F_4= T_0^3T_1+T_1^3T_2+T_2^3T_0 = 0.$$
There are 21 reflection hyperplanes in $\bbP^2$.  They intersect the curve at 84 points forming an orbit with stabilizer subgroups  of order 2. The invariant $F_6$ of degree 6 defines a nonsingular curve of degree 6, the Hessian curve of the Klein quartic. Its equation is given by the Hesse determinant of second partial derivatives of $F_4.$

The double cover of $\bbP^2$ branched along the curve $F_6=0$ is a K3 surface $X$. The automorphism group of $X$ is an infinite group which contains a subgroup isomorphic to $\bar{\Gamma}$.

Next we consider the group $L_3$ of order 648 (No 25). The group $\bar{\Gamma}$ is of order 216 and is known  as the Hesse group \footnote{Another Hesse group is related to 28 bitangents of a plane quartic.}. It is isomorphic to the group of projective transformations leaving invariant the  \emph{Hesse pencil} of plane cubic curves 
\begin{equation}
\lambda (t_0^3+t_1^3+t_2^3)+\mu t_0t_1t_2 = 0.
\end{equation}
It is known that any nonsingular plane cubic curve is projectively isomorphic to one of the curves in the pencil. The base points of the pencil (i.e. points common to all curves from the pencil) are inflection points of each nonsingular member from the pencil. The singular members of the pencil correspond to the values of the parameters 
$(\lambda,\mu) = (0,1),\ (1,-3),\  (1,-3e^{2\pi i/3}),\  (1,-3e^{-2\pi i/3}).$
The corresponding cubic curves are the unions of 3 lines, altogether we get 12 lines which form the 12 reflection hyperplanes.

The smallest degree invariant of $L_3$ in $\bbC^4$ is a polynomial of degree 6
$$F_6 = T_0^6+T_1^6+T_2^6-10(T_0^3T_1^3+T_0^3T_2^3+T_1^3T_2^3).$$
The double cover of $\bbP^2$ branched along the curve $F_6 = 0$ is a K3 surface. Its group of automorphisms is an infinite group containing a subgroup isomorphic to  the Hesse group. 

\medskip
Next we turn our attention to  complex reflection groups of types of types $K_5, L_4$ and $E_6$.
Their orders are all divisible by $6!\cdot 36 = 25,920$ equal to the order of the simple group $\text{PSp}(4,\bbF_3)$. The group of type $E_6$ (the Weyl group of the lattice $E_6$) contains this group as a subgroup of index 2 which consists of words in fundamental reflections of even length. The  group $K_5$ (No 35) is the direct product of $\text{PSp}(4,\bbF_3)$ and a group of order 2.  The group $L_4$ (No 32) is the direct product of  a group of order 3 and $\text{Sp}(4,\bbF_3)$.

Let $\Gamma$ be of type $K_5$. It acts in $\bbP^4$ with  45 reflecting hyperplanes. The hypersurface defined by its invariant of degree 4 is isomorphic to the \emph{Burkhardt quartic} in $\bbP^4$. Its equation can be given in more symmetric form in $\bbP^5$
$$\sum_{i=0}^5T_i = \sum_{i=0}^5T_i^4 = 0.$$
These equations exhibit the action of the symmetric group $S_6$ contained in $\bar{\Gamma}$. It is easy to see that the hypersurface has 45 ordinary double points. This is a record for a hypersurfaces of degree 4 in $\bbP^4$ and this property characterizes Burkhardt quartics. The dual representation of $\Gamma$ is a linear reflection representation too, it is obtained from the original one by composing it with an exterior automorphism of the group. The double points correspond to reflection hyperplanes in the dual space. The reflection hyperplanes in the original space cut out the Burkhardt quartic in special quartic surfaces (classically known as \emph{desmic quartic surfaces}). They are birationally isomorphic to the Kummer surface of the product of an elliptic curve with itself. 

Finally note that the Burkhardt quartic is a compactification of the moduli space of principally polarized abelian surfaces with level 3 structure. All of this and much much more can be found in \cite{Hunt}.

Let $\Gamma$ be of type $L_4$. The number of reflection hyperplanes is 40. The stabilizer subgroup of each hyperplane is the Hesse group of order 648. The smallest invariant is of degree 12. The corresponding hypersurface cuts out in each reflecting hyperplane  the 12 reflecting lines of the Hesse group.  Again for more of this beautiful geometry we refer to Hunt's book.

The geometry of the Weyl group $W(E_6)$ in $\bbP^5$ is fully discussed in Hunt's book. He calls the invariant hypersurface of degree 5  the \emph{gem of the universe}.

\section{Monodromy  groups}

\subsection{Picard-Lefschetz transformations}  Let  $f:X\to S$ be a  holomorphic map of complex manifolds which is a locally trivial $C^\infty$-fibration. One can construct a  complex local coefficient system  whose fibres are the cohomology with compact support  $H_c^n(X_s,\Lambda)$ with some coefficient group $\Lambda$ of the fibres $X_s = f^{-1}(s)$. The  local coefficent system  defines the \emph{monodromy map}
\begin{equation}\label{mon}
\rho_{s_0}:\pi_1(S,s_0) \to \Aut(H_c^n(X_{s_0},\Lambda)).\end{equation}
We will be interested in the cases when $\Lambda = \bbZ, \bbR$, or $\bbC$ that lead to  integral, real or complex monodromy representations. 

The image of the monodromy representation is called the \emph{monodromy group} of the map $f$ (integral, real, complex).

 We refer to \cite{Ar}, \cite{GZ}, \cite{Lo3} for some of  the material which follows.

Let $f:\bbC^{n+1}\to \bbC$ be a holomorphic function with an isolated critical point at $x_0$. We will be interested only in germs  $(f,x_0)$ of $f$ at $x_0$. Without loss of generality we may assume that $x_0$ is the origin and $f(x_0) = 0$. The  level set  $V = f^{-1}(0)$ is an analytic subspace of $\bbC^{n+1}$ with isolated singularity at $0$. The germ of  $(V,0)$ is an \emph{$n$-dimensional isolated hypersurface singularity}. In general the isomorphism type of the germ of $(V,0)$ does not determine the isomorphism type of the germ $(f,0)$. However, it does in one important case when  $f$ is  a \emph{weighted homogeneous polynomial}.\footnote{This means that $f(z_1,\ldots,z_{n+1})$ is a linear combination of monomials of the same degree,  where each variable is weighted with some positive number.} In this case we say that the germ $(V,0)$ is a weighted homogeneous isolated hypersurface singularity.

It was shown by J. Milnor  \cite{Milnor} that for sufficiently small $\epsilon$ and $\delta$, the restriction of $f$ to 
$$X = \{z\in \bbC^{n+1}:||z|| < \epsilon, 0 < |f(z)|  <\delta\}$$
is a locally trivial $C^\infty$  fibration whose fibre is an open $n$-dimensional complex manifold. Moreover, each fibre has the homotopy type of a bouquet of $n$-spheres.  The number $\mu$ of the spheres is equal to the multiplicity of $f$ at $0$ computed as  the dimension of the jacobian algebra
\begin{equation}\label{jac}
J_f =
 \dim_\bbC \bbC[[z_1,\ldots,z_{n+1}]]/(\frac{\partial f}{\partial z_1},\ldots,\frac{\partial f}{\partial z_{n+1}}).\end{equation}

Let $D_\epsilon^* = \{t\in \bbC:0 < |t| < \epsilon\}$ and $f:X\to D_\delta^*$ be the above fibration, a \emph{MIlnor fibration} of $(V,0). $  Let 
\begin{equation}
M_t = H_c^n(F_t,\bbZ) \cong H_n(F_t,\bbZ) .\end{equation}

 The bilinear pairing  
$$H_c^n(F_t,\bbZ)\times H_c^n(F_t,\bbZ)\to H_c^{2n}(F_t,\bbZ)\cong \bbZ$$
is symmetric if $n$ is even, and skew-symmetric otherwise. To make the last isomorphism unique we fix an orientation on $F_t$ defined by the complex structure on the open ball $B_\epsilon = \{z\in \bbC^{n+1}:||z|| < \epsilon\}$.  
Thus, if $n$ is even, that we will assume from now on, the bilinear pairing 
equips $M_t$ with a structure of a lattice, called the \emph{Milnor lattice} of isolated critical point of $f$ at $0$. Its isometry class does not depend on the choice of  the point $t$ in $\pi_1(D_\delta^*;t_0)$. Fixing a point $t_0\in D_\epsilon^*$ we obtain the \emph{classical monodromy map}
$$\pi_1(D_\epsilon^*;t_0) \cong \bbZ \to \Or(M_{t_0}).$$
 Choosing a generator of $\pi_1(D_\epsilon^*;t_0)$ it defines an isometry of the Milnor lattice which is called a \emph{classical monodromy operator}.

\begin{example}\label{quadric} Let 
$$q(z_1,\ldots,z_{n+1})= z_1^2+\ldots+z_{n+1}^2.$$
It has a unique critical point at the origin with critical value $0$. An isolated critical point $(f,x_0)$ (resp. isolated hypersurface singularity $(V,x_0)$) locally analytically isomorphic to  $(q,0)$ (resp. $(q^{-1}(0),0)$)  is called  \emph{nondegenerate critical point} (resp. \emph{ordinary double point} or \emph{ordinary node}). 

Let $\epsilon$ be any positive number and $\delta < \sqrt{\epsilon}$. For any $t \in D_\delta^*$, the intersection  
$F_t = q^{-1}(t)\cap B_\epsilon$ is non-empty and is given by the equations
$$\sum_{i=1}^{n+1}z_i^2 = t, \quad ||z|| < \epsilon.$$
Writing  $z_i = x_i+\sqrt{-1}y_i$ we find the real equations
$$||x||^2 -||y||^2 = t, \quad x\cdot y = 0, \quad ||x||^2+||y||^2  < \epsilon^2.$$
After some smooth  coordinate change, we get the equations
$$||x||^2 = 1, \quad x\cdot y = 0, \quad ||y||^2 < 1.$$
It is easy to see that these are the equations of the open unit ball subbundle of the tangent bundle of the  $n$-sphere $S^{n}$. Thus we may take $F_t$ to be a Milnor fibre of $(q,0)$.
 The sphere $S^n$ is contained in $F_t$ as the zero  section of the tangent bundle, and $F_t$ can be obviously retracted to $S^n$. Let $\alpha$ denote the fundamental class  $[S^n]$ in $H_c^n(F_t,\bbZ)$. We have 
$$H_c^n(F_t,\bbZ) = \bbZ  \alpha.$$
Our orientation on $F_t$ is equal to the orientation of the tangent bundle of the sphere taken with the sign $(-1)^{n(n-1)/2}$. This gives
\[
(\alpha,\alpha) = (-1)^{n(n-1)/2}\chi(S^n) = (-1)^{n(n-1)/2}(1+(-1)^n).
\]
and hence in the case when $n$ is even,
\begin{equation}\label{self}
(\alpha,\alpha) = \begin{cases}
      -2& \text{if} \ n\equiv 2 \mod 4\\
      2& \text{if}\ n\equiv 0 \mod 4.
\end{cases}.
\end{equation}

By collapsing the zero section $S^n$ to the point we get a map $F_t$ to $B_\epsilon \cap q^{-1}(0)$ which is a diffeomorphism outside $S^n$. For this reason the homology class $\delta$ is called the \emph{vanishing cycle}.

The classical \emph{Picard-Lefschetz formula} shows that the monodromy operator is a reflection transformation
\begin{equation}\label{pl}
T(x) = x-(-1)^{n(n-1)/2}(x,\delta)\delta.
\end{equation}

\end{example}

\medskip
Let $(V,0)$ be an isolated $n$-dimensional hypersurface singularity. A \emph{deformation} of $(V,0)$ is  a holomorphic map-germ 
$\phi:(\bbC^{n+k},0)\to (\bbC^k,0)$ with fibre $\phi^{-1}(0)$ isomorphic to $(V,0)$. By definition $(V,0)$ admits a deformation map $f:(\bbC^{n+1},0)\to (\bbC,0)$. Let $J_f$ be the jacobian algebra of $f$ \eqref{jac} and 
\begin{equation}\label{jacc}
J_{f=0} = J_f/(\bar{f}),
\end{equation}
where $\bar{f}$ is the coset of $f$ in $J_f$. Its dimension $\tau$ (the \emph{Tjurina number}) is less or equal than  $\mu$. The equality occurs if and only if $(V,0)$ is isomorphic to a weighted homogeneous singularity.
 
Choose a basis $z^{a_1},\ldots,z^{a_\mu}$ of the  algebra \eqref{jacc} represented by monomials with exponent vectors $a_1,\ldots,a_\tau$ with $a_\tau = 0$ (this is always possible). Consider a deformation 
\begin{equation}\label{sud}
\Phi:(\bbC^{n+\tau},0)\to (\bbC^\tau,0), (z,u)\to (f(z)+\sum_{i=1}^{\tau-1}u_iz^{a_i}, u_1,\ldots,u_{\tau-1}).
\end{equation}  
This deformation represents a \emph{semi-universal} (\emph{miniversal}) deformation of $(V,0)$.   Roughly speaking this means that any deformation $\phi:(\bbC^{n+k},0)\to (\bbC^k,0)$ is obtained from \eqref{sud} by mapping $(\bbC^k,0)$ to $(\bbC^\tau,0)$ and taking the pull-back of the map $\phi$ (this explains \emph{versal} part). The map-germ $(\bbC^k,0)\to (\bbC^\tau,0)$ is not unique, but its derivative at $0$ is unique (whence the semi in semi-universal).

Let 
$$\Delta = \{(t,u)\in \bbC^\tau: \Phi^{-1}(t,u)\  \text{has a singular point  at some $(z,u)$}\}.$$
Let $(\Delta,0)$ be the germ of $\Delta$ at $0$. It is called the \emph{bifurcation diagram}  or \emph{discriminant} of $f$.
Choose a representative of $\Phi$ defined by restricting the map  to some open ball $B$ in 
$\bbC^{n+\tau}$ with center at $0$ and let $U$ be some open neighborhood  of $0$ in $\bbC^\tau$ such that $\pi^{-1}(U) \subset B$. When $B$ is small enough one shows that the restriction map
$$\pi:\pi^{-1}(U\setminus U\cap \Delta)\to U\setminus \Delta$$
is a locally trivial $C^\infty$-fibration with fibre diffeomorphic to a Milnor fibre $F$ of $(f,0)$. Fixing a point
$u_0\in U\setminus \Delta$, we get the \emph{global monodromy map} of $(X,0)$.
\begin{equation}
\pi_1(U\setminus \Delta;u_0)\to \Or(H_c^n(\pi^{-1}(u_0),\bbZ)) \cong \Or(M),
\end{equation}
where $M$ is a fixed lattice in the isomorphism class of Milnor lattices of $(f,0)$. 

One can show that the conjugacy class of the image of the global monodromy map (the \emph{global monodromy group} of $(X,0)$)  does not depend on the choices of $B, U, u_0$.

\begin{theorem} The global monodromy group $\Gamma$ of an isolated hypersurface singularity of even dimension 
$n = 2k$ is a subgroup of  $\Ref_{-2}(M)$ if $k$ is odd and $\Ref_2(M)$ otherwise.
\end{theorem}

\begin{proof} This follows from the Picard-Lefschetz theory. Choose a point $(t_0,u^{(0)})\in \bbC^\tau\setminus \Delta$ close enough to the origin and pass  a general line 
$u_i = c_i(t-t_0)+u_i^{(0)}$
through this point. The pre-image of this line in $\bbC^{n+\tau}$ is the 1-dimensional deformation
$\bbC^{n+1}\to \bbC$ given  the function $t = \tilde{f}(z)$ implicitly defined by the equation
$$ f(z)+t(-1+\sum_{i=1}^{\tau-1} c_iz_i^{a_i})+\sum_{i=1}^{\tau-1}(u_i^{(0)}-c_it_0)z_i^{a_i} = 0.$$
By choosing the line general enough, we may assume that all critical points of this function are non-degenerate and critical values are all distinct. Its fibre over the point $t=t_0$ is equal to the fibre of the function 
$f(z)+\sum_{i=1}^{\tau-1}u_i^{(0)}z_i^{a_i}$ over $t_0$. Since the semi-universal deformation is a locally trivial fibration over the complement of $\Delta$ all nonsingular fibres are diffeomorphic. It follows from the Morse theory that the number of critical points of $\tilde{f}$ is equal to the Milnor number $\mu$. Let $t_1,\ldots,t_\mu$ be the critical values.  Let $S(t_i)$ be a small circle around $t_i$ and $t_i'\in S_i$. The function $\tilde{f}$ defines a locally trivial $C^\infty$-fibration over $\bbC\setminus \{t_1,\ldots,t_\mu\}$ so that we can choose a diffeomorphism of the fibres
$$\phi_i:\tilde{f}^{-1}(t_0) \to \pi^{-1}(u_0).$$
The Milnor fibre $F_i$ of $\tilde{f}$ at $t_i$ is an open subset of $\tilde{f}^{-1}(t_i')$ and hence defines an inclusion of the lattices
$$H_n(F_i,\bbZ)\hookrightarrow H_n(\tilde{f}^{-1}(t_i'),\bbZ) \cong H_n(\tilde{f}^{-1}(t_0),\bbZ)\cong H_n(\pi^{-1}(u_0),\bbZ)\cong M.$$
The image $\delta_i$ of the vanishing cycle $\delta_i'$ generating $H_n(F_i,\bbZ)$ in $M$ is called a \emph{vanishing cycle} of the Milnor lattice $M$. The vanishing cycles $(\delta_1,\ldots,\delta_\mu)$ form a basis in $M$. 

For each $t_i$ choose a path $\gamma_i$ in $\bbC\setminus \{t_1,\ldots,t_\mu\}$ which connects $t_0$ with $t_i'$ and when continues along the circle $S(t_i)$ in counterclockwise way until returning to $t_i'$ and going back to $t_0$ along the same path but in the opposite direction. The homotopy classes $[\gamma_i]$ of the paths $\gamma_i$ generate $\pi_1(\bbC\setminus \{t_1,\ldots,t_\mu\};t_0)$, the map 
$$s_*:\pi_1(\bbC\setminus \{t_1,\ldots,t_\mu\};t_0)\to \pi_1(U\setminus \Delta;u_0)$$
is surjective and the images $g_i$ of $[\gamma_i]$ under the composition of the monodromy map and $s_*$ generate the monodromy group $\Gamma$. Applying the Picard-Lefschetz formula \eqref{pl} we obtain, for any $x\in M$;
$$g_i(x) = x-(-1)^{n(n-1)/2}(x,\delta_i)\delta_i, \quad i = 1,\ldots,\mu.$$
This shows that $\Gamma$ is generated by $\mu$ reflections in elements $\delta_i$ satisfying \eqref{self}. This proves the claim.
\end{proof}

\begin{remark}\label{whcase}
It is known that for any isolated hypersurface singularity given by a weighted homogeneous polynomial $P$ with isolated critical point at $0$ the Milnor lattice is isomorphic to
$H_c^n(P^{-1}(\epsilon),\bbZ)$. In this case we can define the  monodromy group as the image of the monodromy map
$$\pi_1(\bbC^\mu\setminus \Delta;u_0) \to \Or(H_c^n(\pi^{-1}(u_0),\bbZ)),$$
where $u_0\in \bbC^\mu\setminus \Delta$. 
\end{remark}

One can relate the global monodromy group with the classical monodromy operator.

\begin{theorem}\label{cox} Let $T$ be a classical monodromy operator and $s_1,\ldots,s_\mu$ be the reflections in $\Or(M)$ corresponding to a choice of paths $\gamma_i$ as in above.  Then $\pm T$ is conjugate to the product $s_1\cdots s_\mu$.
\end{theorem}

The product $s_1\cdots s_\mu$ is an analog of a  Coxeter element in  the reflection  group of $M$ (see \ref{ex3}). However, the Weyl group of $M$ is not a Coxeter group in general.

\subsection{Surface singularities}   Assume now that $n = 2$.   We will  represent a surface singularity by an isolated singular point $0$ of an affine surface $X\subset \bbC^3$. 

 Let $\pi:X'\to X$ be a \emph{resolution of the singularity} $(X,0)$. This means that $X'$ is a nonsingular algebraic surface, and $\pi$ is a proper holomorphic map which is an isomorphism over $U = X\setminus\{0\}$.  We  may assume to be minimal, i.e. does not contain $(-1)$-curves in its fibre over $0$. This assumption makes it unique up to isomorphism.

One defines the following invariants of $(X,0)$. The first invariant is the \emph{genus} $\delta p_a$ of $(X,0)$. This is the dimension of the first cohomology space of $X'$ with coefficient in the structure sheaf $\calO_{X'}$ of regular (or holomorphic) functions on $X'$. In spite of $X'$ not being compact, this space is finite-dimensional. 

Our second invariant is the \emph{canonical class square} $\delta K^2$. To define it we use that there is a rational differential 2-form on $X$ which has no poles or zeros on $X\setminus \{0\}$. Its extension to $X'$ has divisor $D$ supported at the exceptional fibre. It represents an element $[D]$ in $ H_c^2(X',\bbZ)$. Using the cup-product pairing we get the number 
$[D]^2$ which we take for $\delta K^2$.  

Finally we consider the fibre $E= \pi^{-1}(0)$ of a minimal resolution. This is a (usually reducible) holomorphic curve on $X$. We denote its Betti numbers  by $b_i$. For example, $b_2$ is the number of irreducible components of $E$. 

Since $X'$ is defined uniquely up to isomorphism, the numbers $\delta p_a,\delta K^2, b_i$   are well-defined. 

\begin{remark} The notations $\delta p_a$ and $\delta K^2$ are explained as follows. Suppose $Y$ is a projective surface  of degree $d$ in $\bbP^{3}$ with isolated singularities $y_1,\ldots,y_k$. Let $Y'$ be its minimal resolution.  Define the arithmetic genus of a nonsingular projective surface $V$ as
$p_a = -q+p_g$, where $q = \dim H^1(V,\calO_V) = \dim \Omega^1(X)$ is the dimension of the space of holomorphic 1-forms on $V$ and $p_g = \dim H^2(V,\calO_V) = \dim \Omega^2(X)$ is the dimension of the space of holomorphic 2-forms on $V$. Let $F_d$ be a nonsingular hypersurface of degree $d$ in $\bbP^3$. Then
$$
p_a(Y') = p_a(F_d)-\sum_{i=1}^k \delta p_a(x_i), \quad 
K_{Y'}^2 = K_{F_d}^2- \sum_{i=1}^k \delta K^2(x_i).$$
The numbers  $ p_a(F_d)$ and   $K_{F_d}^2$ are easy to compute. We have
$$
p_a(F_d) = p_g(F_d) = (d-1)(d-2)(d-3)/6, \quad 
K_{F_d}^2= d(d-4)^2. 
$$
\end{remark}

\begin{theorem}(J. Steenbrink \cite{Steenbrink}) Let $\mu$ be the Milnor number of  a surface singularity $(V,0)$. Then the signature  $(\mu_+,\mu_-,\mu_0)$ of the Milnor lattice is given as follows
$$\mu_0= b_1,\quad  \mu_++\mu_- +\mu_0 = \mu,\quad  \mu_+-\mu_- = -\delta K^2-b_2-8\delta p_a.$$
\end{theorem}

\begin{example}  The following properties are equivalent:
\begin{itemize} 
\item $\delta p_a = 0$;
\item  $\mu_+ = \mu_0 = 0$;
\item $\delta K^2 = 0$;
\item $b_2 = \mu$;
\item the exceptional curve of a minimal resolution is the union of  nodal curves.
\item $V$ can be given by equation $P(z_1,z_2,z_3) = 0$, where $P$ is a weighted homogeneous polynomial of degree $d$ with respect to positive weights $q_1,q_2,q_3$ such that $d-q_1-q_2-q_3 < 0$.
\item $(V,0)$ is isomorphic to an affine surface with ring of regular functions isomorphic to the ring of invariant polynomials of a   finite subgroup $G\subset \SL(2,\bbC)$.
\end{itemize}

The following table gives the list of isomorphism classes of surface singularities characterized by the previous properties. They go under many different names  \emph{simple singularities, ADE singularities, Du Val singularities, double rational points, Gorenstein quotient singularities, Klein singularities}.

\begin{table}[h]
 \caption{Simple surface singularities }\label{simple }
\begin{center} 
\begin{tabular}{|| l | r | r | r | r |r|}\hline
 Type&Polynomial&Weights&Degree&G\\ \hline
   $A_{2k}, k\ge 1$&$Z_1Z_2+Z_3^{2k+1}$&$(2k+1,2k+1,2)$ &$4k+2$&C$_{4k}$ \\ \hline
    $A_{2k+1}, k\ge 0$&$Z_1Z_2+Z_3^{2k+2}$&$(k+1,k+1,1)$ &$2k+2$&C$_{2k+1}$ \\\hline
    $D_n, n\ge 4$&$Z_1^2+Z_2Z_3^2+Z_2^{n-1}$&$(n-1,2,n-2)$ &$2n-2$&$\bar{D}_{4n-4} $\\ \hline
    $E_{6}$&$Z_1^2+Z_2^3+Z_3^4$&$(6,4,3)$ &$12$&$\bar{T}_{24}$ \\ \hline
    $E_7$&$Z_1^2+Z_2^3+Z_2Z_3^3$&$(9,6,4)$ &$18$&$\bar{O}_{48}$ \\ \hline
    $E_8$&$Z_1^2+Z_2^3+Z_3^{5}$&$(15,10,6)$ &$30$&$\bar{I}_{120}$ \\ \hline
  
\end{tabular}
\vspace{10pt}
\end{center}
\end{table}

The following result follows immediately from  Brieskorn's and Tjurina's construction of simultaneous resolution of simple surface singularities \cite{Br2}, \cite{Pinkham}, \cite{Tj}.

\begin{theorem} The Milnor lattice of a simple surface singularity is a finite root lattice of the type indicated in the first column of Table with quadratic form multiplied by $-1$. The monodromy group is the corresponding finite real reflection group.
\end{theorem}

It is known that the incidence graph of the irreducible components of a minimal resolution of singularities is the Coxeter graph of the corresponding type. This was first observed by P. Du Val \cite{DV}.

Consider the monodromy group of a simple surface singularity $\Gamma$ of $(X,0)$ as the image of the monodromy map \eqref{whcase}. For each Coxeter system $(W,S)$ with Coxeter matrix $(m_{ss'})$ one defines the associated \emph{Artin-Brieskorn braid group} $B_W$ by presentation
$$B_W = \{g_s, s\in S:\underbrace{(g_sg_{s'}g_s)\cdots (g_sg_{s'}g_s)}_{m(s,s')-2} = \underbrace{(g_{s'}g_{s}g_{s'})\cdots (g_{s'}g_{s}g_{s'})}_{m(s,s')}\}$$
If we impose additional relations $g_s^2 = 1, s\in S$ we get the definition of $(W,S)$. This defines an extension of groups
\begin{equation}
1\to \tilde{B}_W\to B_W \to W \to 1,
\end{equation}
where $\tilde{B}_W$ is the normal subgroup of $B_W$ generated by conjugates of $g_s^2$.

Brieskorn also proves the following.

\begin{theorem}\label{bri} Let $\Delta \subset \bbC^\mu$ be the discriminant of a simple surface singularity. Then $\pi_1(\bbC^\mu \setminus \Delta,u_0)$ is isomorphic to the braid group $B_\Gamma$ of the monodromy group $\Gamma$. The regular covering $U\to \Delta \subset \bbC^\mu$ corresponding to the normal subgroup $\tilde{B}_\Gamma$ can  be $\Gamma$-equivariantly extended to the covering $V\to V/\Gamma \cong \bbC^\mu$, where $V$ is a complex vector space of dimension $\mu$ on which $\Gamma$ acts as a reflection group. The pre-image of $\Delta$ in $V$ is the union of reflection hyperplanes.
\end{theorem}


\end{example}

\begin{example} The following properties are equivalent.
\begin{itemize}
\item $\mu_+ = 0, \mu_0 > 0$;
\item $\mu_+ = 0, \mu_0 = 2$;
\item the exceptional curve of a minimal resolution is a nonsingular elliptic curve;
\item $(V,0)$ can be represented by the zero level of one of the polynomials given in Table \ref{se}.
\end{itemize}

These singularities are called \emph{simple elliptic singularities}.

\begin{table}[h]
 \caption{Simple elliptic surface singularities }\label{se}
\begin{center} 
\begin{tabular}{|| l | r | r | r | r |}\hline
 Type&Polynomial&Weights&Degree\\ \hline
     $P_8$&\small{$z_1^3+z_2^3+z_3^4+\lambda z_1z_2z_3$}&$(1,1,1)$ &$3$ \\ \hline
    $X_9$&\small{$z_1^2+z_2^4+z_3^4+\lambda z_1z_2z_3$}&$(2,1,1)$ &$4$ \\ \hline
    $J_{10}$&\small{$z_1^2+z_2^3+z_3^6+\lambda z_1z_2z_3$}&$(3,2,1)$ &$6$ \\ \hline
  
\end{tabular}
\vspace{10pt}
\end{center}
\end{table}

Here the subscript is equal to $\mu$. 

\begin{theorem} (A. Gabrielov \cite{Ga1}) Let $M$ be the Milnor lattice of a simple elliptic singularity. Then $M^\perp$ is of rank 2 and $M/M^\perp$ is  isomorphic to the root lattice of type $E_{\mu-2}$.  The image $\bar{G}$ of the monodromy group $\Gamma$ in $(M/M^\perp)$ is the finite  reflection group of type $E_{\mu-2}$. The monodromy group is isomorphic to the semi-direct product $(M^\perp\otimes M/M^\perp)\rtimes W(E_{\mu-2})$ and can be naturally identified with  an affine complex crystallographic reflection group with linear part $W(E_{\mu-2})$. 
\end{theorem}

There is a generalization of Theorem \ref{bri} to the case of simple elliptic singularities due to E. Looijenga \cite{Lo1} and \cite{Pi1}. It involves affine crystallographic reflection groups and uses Theorem \ref{BS}.  
\end{example}

\begin{example} The following properties are equivalent:
\begin{itemize} 
\item $\mu_+ = 1$;
\item  $\mu_+ =  1, \mu_0 = 1$;
\item $V$ can be given by equation $P(z_1,z_2,z_3) = 0$, where 
$$P = z_1^a+z_2^b+z_3^c+\lambda z_1z_2z_3,  \quad \frac{1}{a}+\frac{1}{b}+\frac{1}{c} < 1, \ \lambda\ne 0.$$
\end{itemize}
These singularities are called \emph{hyperbolic unimodal singularities}

\begin{theorem}(A. Gabrielov \cite{Ga1}) The Milnor lattice $M$ of a hyperbolic singularity is isomorphic to the lattice
$$E_{p,q,r}(-1)\perp \langle 0\rangle.$$
The monodromy group is the semi-direct product $\bbZ^{\mu}\rtimes W(p,q,r).$ Its image in $\Or(M/M^\perp)$ is the reflection group $W(p,q,r)$.
\end{theorem}
\end{example} 

The previous classes of isolated surface singularities are characterized by the condition $\mu_+\le 1$. If  $\mu_+\ge 2$,  the monodromy group is always of finite index  in $\Or(M)$ (see \cite{Ebeling1}, \cite{Ebeling2}). Together with the previous theorems this implies that the monodromy group is always of finite index in $\Or(M)$ except in the case of hyperbolic unimodal singularities with $(p,q,r) =  (2,3,7), (2,4,5), (3,3,4)$
(see Example \ref{ex4}).

There is a generalization of Brieskorn's theorem \ref{bri} to the case of  hyperbolic singularities  due to E. Looijenga \cite{Lo1},\cite{Lo2}.

\section{Symmetries of singularities}
\subsection{Eigen-monodromy groups}
Suppose we have a holomorphic map $f:X\to S$ as in the section defining the monodromy map \eqref{mon}.  Suppose also that a finite group $G$ acts on all fibres of the map in a compatible way. This means that there is an action of $G$ on $X$ which leaves fibres invariant. Then the cohomology groups $H_c^n(X_s,\bbC)$ become representation spaces for $G$ and we can decompose them into irreducible components
$$H_c^n(X_s,\bbC) = \bigoplus_{\chi\in \textup{Irr}(G)}H_c^n(X_s,\bbC)_\chi.$$
One checks that the monodromy map decomposes too and defines the $\chi$-monodromy map
$$\rho_{s_0}^\chi:\pi_1(S;s_0)\to \GL(H_c^n(X_s,\bbC)_\chi).$$
Let  $E$ be a real vector space equipped with a bilinear form  $(v,w)$, symmetric or skew-symmetric. Let $E_\bbC$ be its complexification with the conjugacy map $v\mapsto \bar{v}$.  We extend the bilinear form on $E$ to $E_\bbC$ by linearity. It is easy to see that it satisfies 
$(\bar{x},\bar{y}) = \overline{(x,y)}$. Next we equip $E_\bbC$ with a hermitian form defined by
$$\langle x,y \rangle =\begin{cases}
   (x,\bar{y}), & \text{if $(x,y)$ is symmetric}, \\
     i(x,\bar{y}) & \text{otherwise}.
\end{cases}
$$
We apply this to $E_\bbC = H_c^n(X_s,\bbC)$ with the bilinear map defined by the cup-product.  The $\chi$-monodromy map leaves the corresponding hermitian form invariant and defines a homomorphism
$$\rho_{s_0}:\pi_1(S;s_0)\to \U(H_c^n(X_s,\bbC)_\chi).$$
We are interested in examples when the image of this homomorphism is a complex reflection group.

\subsection{Symmetries of  singularities}
Assume that the germ of an isolated hypersurface singularity $(X,0)$ can be represented by a polynomial $f$ which is invariant with respect to some finite subgroup $G$ of $\GL(n+1,\bbC)$.  One can define the notion of a $G$-equivariant  deformation of $(X,0)$ and show that a semi-universal $G$-equivariant deformation of $(X,0)$ can be given by the germ of the map
\begin{equation}
\Phi_G:(\bbC^{n+\tau'},0)\to (\bbC^{\tau'},0), (z,u)\to (f(z)+\sum_{i=1}^{\tau'-1}u_ig_i, u_1,\ldots,u_{\tau'-1}),\end{equation}
where $(g_1,\ldots,g_{\tau'})$ is a basis of the  subspace $J_{f=0}^{G,\chi}$ of relative invariants of the algebra \eqref{jacc}. In the case when $G$ cyclic, we can choose $g_i$'s to be monomials. The equivariant discriminant $\Delta_G$ is defined the same as in the case of the trivial action. Now we define the $G$-equivariant monodromy group $\Gamma_G$ of $(X,x_0)$ following the definition in the case $G = \{1\}$. The group $G$ acts obviously on the Milnor fibre of $f$ and hence on the Milnor lattice  $M$ giving it a structure of a $\bbZ[G]$-module. Since 
$$\bbC^{\tau'}\setminus \Delta_G \subset \bbC^{\tau}\setminus \Delta$$
we can choose a point $s_0\in \bbC^{\tau'}\setminus \Delta_G$ to define a homomorphism
$$i_{s_0}:\pi_1(\bbC^{\tau'}\setminus \Delta_G;s_0) \to \pi_1(\bbC^{\tau}\setminus \Delta;s_0).$$
This homomorphism induces  a natural injective homomorphism of the monodromy groups $i_*:\Gamma_G \to \Gamma$. This allows us to identify $\Gamma_G$ with a subgroup of $\Gamma$

\begin{proposition}(P. Slodowy \cite{Sl1}) \label{pro1} Let $\Gamma$ be the monodromy group of $(X,x_0)$ and $\Gamma_G$ be the $G$-equivariant monodromy group. Then 
$$\Gamma_G = \Gamma \cap \Aut_{\bbZ[G]}(M).$$
\end{proposition}

Consider the natural action of  $G$ on the jacobian algebra \eqref{jac} via its action on the partial derivatives of the function $f$. 

\begin{theorem}(C.T.C. Wall \cite{Wall}) \label{wall} Assume $(X,x_0)$ is an  isolated hypersurface singularity defined by a holomorphic function $f:\bbC^{n+1}\to \bbC$. Let $J_f$ be its jacobian algebra. There is an isomorphism of $G$-modules
$$M_\bbC \cong J_f\otimes \textup{det}_G,$$
where $\det_G$ is the one-dimensional representation of $G$ given by the determinant.
\end{theorem}

\begin{example} Let $f(z_1,z_2,z_3) = z_1z_2+z_3^{2k}$ be a simple surface singularity of type $A_{2k-1}$. Consider the action of the group $G = \bbZ/2\bbZ$ by $(z_1,z_2,z_3)\mapsto (z_1,z_2,-z_3)$. We take $1,z_3,\ldots,z_3^{2k-2}$ to be a basis of the jacobian algebra. Thus we have $k$ invariant monomials $1,z_3^2,\ldots,z_3^{2k-2}$ and $k-1$ anti-invariant monomials $z_3,\ldots,z_3^{2k-3}$.
It follows from Theorem \ref{wall} that $M_\bbC$ is the direct sum of $k-1$-dimensional invariant part $M_+$ and $k$-dimensional anti-invariant part $M_-$.  We have
$$\Or(M) = W(A_{2k-1})\rtimes (\tau),$$
where $\tau$ is the non-trivial symmetry of the Coxeter diagram of type $A_{2k-1}$. In fact it is easy to see that the semi-direct product is the direct product. Let $\alpha_1,\ldots,\alpha_{2k-1}$ be the fundamental root vectors. We have $\tau(\alpha_i) = \alpha_{2k-i}$, hence $\dim M^\tau = k$. This shows that the image $\sigma$ of the generator of $G$ in $\Or(M)$ is not equal to $\tau$. In fact, it must belong to  $W(A_{2k-1})$. To see this we use that all involutions in $W(M)\cong \Sigma_{2k}$ are the products of at most $k$ transpositions, hence their fixed subspaces are of dimension $\le k-1$. The group $\Or(M)$ is the product of $W(M)$ and $\pm 1$, thus all involution in $\Or(M)\setminus W(M)$ have fixed subspaces of dimension $\le k-1$. The conjugacy class of $\sigma\in W(A_{2k-1})\cong \Sigma_{2k}$ is determined by the number $r$ of disjoint transpositions in which it decomposes. We have $\dim M^\sigma = 2k-r-1$. Thus $\sigma$ is conjugate to the product of $k$ disjoint transpositions. It follows from the model of the lattice $A_{2k-1}$ given in \eqref{space} that $\sigma$ is conjugate to the transformation $\alpha_i\mapsto -\alpha_{2k-i}$. Hence  sublattice $M_-$ is generated by $\beta_i = \alpha_i+\alpha_{2k-i}, i = 1,\ldots,k-1,$ and $\beta_k = \alpha_k$. We have 
$$(\beta_i,\beta_j) = \begin{cases}
      -4& \text{if}\  i= j\ne k, \\
      -2& \text{if}\  i= j = k,\\
      2&\text{if}\  |i-j| = 1,\\
      0&\text{if}\ |i-j| > 1.
\end{cases}
$$
Comparing this with Example \ref{rootlat} we find that $M = N(2)$, where $N$ is the lattice defining an integral structure for the reflection group of type $B_k$. In other words, the reflections $r_{\beta_i}$ generate the Weyl group of type $B_k$. 

Similar construction for a symmetry of order 2 of singular points of types $D_{k+1}, \\n\ge 4,$ (resp. $E_6$, resp. $D_4$) lead to the  reflection groups of type  $B_k$ (resp. $F_4$, resp. $G_2= I_2(6)$). The Milnor lattices obtained as the invariant parts of the Milnor lattice of type $A_{2k-1}$ and $D_{k+1}$ define the same reflection groups but their lattices are similar but not isomorphic.

 \begin{figure}[h]
\xy
(-10,0)*{};(0,10)*{};
@={(0,3),(5,3),(10,3),(20,3),(25,3),(30,3),(35,0),(0,-3),(5,-3),(10,-3),(20,-3),(25,-3),(30,-3)}@@{*{\bullet}};
(0,3)*{};(12,3)*{}**\dir{-};(0,-3)*{};(12,-3)*{}**\dir{-};
(18,3)*{};(30,3)*{}**\dir{-};(18,-3)*{};(30,-3)*{}**\dir{-};
(30,3)*{};(35,0)*{}**\dir{-};(30,-3)*{};(35,0)*{}**\dir{-};
(15,3)*{\ldots};(15,-3)*{\ldots};
@={(0,-10),(5,-10),(10,-10),(20,-10),(25,-10),(30,-10),(35,-10)}@@{*{\bullet}};
(0,-10)*{};(12,-10)*{}**\dir{-};(18,-10)*{};(30,-10)*{}**\dir{-};(15,-10)*{\ldots};
(30,-10)*{};(35,-10)*{}**\dir2{-};
(43,0)*{A_{2k-1}}; (43,-10)*{B_{k}};
@={(60,-3),(60,3),(65,0),(70,0),(80,0),(85,0),(90,0),(95,0)}@@{*{\bullet}};
(60,3)*{};(65,0)*{}**\dir{-};(60,-3)*{};(65,0)*{}**\dir{-};
(65,0)*{};(72,0)*{}**\dir{-};
(78,0)*{};(95,0)*{}**\dir{-};
(75,0)*{\ldots};
@={(60,-10),(65,-10),(70,-10),(80,-10),(85,-10),(90,-10),(95,-10)}@@{*{\bullet}};
(65,-10)*{};(72,-10)*{}**\dir{-};(78,-10)*{};(95,-10)*{}**\dir{-};
(60,-10)*{};(65,-10)*{}**\dir2{-};(75,-10)*{\ldots};
(102,0)*{D_{k+1}}; (102,-10)*{B_{k}};
\POS (-3,-3)*{}="a"
, (-3,3)*{}="b"
\POS "a" \ar@{<->} "b"
\POS (57,-3)*{}="a"
, (57,3)*{}="b"
\POS "a" \ar@{<->} "b"
\endxy

\xy
(-10,0)*{};(0,-10)*{};
@={(0,-22),(5,-22),(0,-28),(5,-28),(10,-25),(15,-25)}@@{*{\bullet}};
(0,-22)*{};(5,-22)*{}**\dir{-};(0,-28)*{};(5,-28)*{}**\dir{-};
(5,-22)*{};(10,-25)*{}**\dir{-}; (5,-28)*{};(10,-25)*{}**\dir{-};
(10,-25)*{};(15,-25)*{}**\dir{-};
(22,-25)*{E_6};
@={(0,-35),(5,-35),(10,-35),(15,-35)}@@{*{\bullet}};
(0,-35)*{};(5,-35)*{}**\dir{-};(10,-35)*{};(15,-35)*{}**\dir{-};
(5,-35)*{};(10,-35)*{}**\dir2{-};
(22,-35)*{F_4};
@={(60,-25),(60,-29),(60,-22),(66,-25)}@@{*{\bullet}};
(66,-25)*{};(60,-25)*{}**\dir{-};(66,-25)*{};(60,-29)*{}**\dir{-};(66,-25)*{};(60,-21)*{}**\dir{-};
(60,-35)*{\bullet};(66,-35)*{\bullet}**\dir3{-};
(72,-25)*{D_4}; (72,-35)*{G_2};
\POS (-3,-28)*{}="a"
, (-3,-22)*{}="b"
\POS "a" \ar@{<->} "b"
\POS (55,-20)*{}="a"
, (55,-30)*{}="b"
\POS "a" \ar@{<->} "b"
\POS (57,-20)*{}="a"
, (57,-24)*{}="b"
\POS "a" \ar@{<->} "b"
\POS (57,-30)*{}="a"
, (57,-26)*{}="b"
\POS "a" \ar@{<->} "b"

\endxy

\caption{}

\end{figure}

\end{example}

 \begin{remark} The appearance of the Dynkin diagrams of type $B_n,F_4,G_2$ in the theory of simple singularities was first noticed by P. Slodowy \cite{Sl3}, and, from different, but equivalent,  perspective in the work of Arnol'd on critical points on manifolds with boundary \cite{Arnold},\cite{Ar}. In the theory of simple surface singularities over non-algebraically closed field they appear in  \cite{Li}. 
 \end{remark}
  
\begin{example} Let $(X,0)$ be a simple surface singularity of type $E_6$ given by equation
$z_1^2+z_2^3+z_3^4 = 0$. Consider the group $G$ generated by the symmetry $g$ of order 3 given by 
$(z_1,z_2,z_3)\mapsto (z_1,\eta_3z_2,z_3)$, where $\eta_3 = e^{2\pi i/3}$. A monomial basis of the jacobian algebra is $(1,z_2,z_3,z_3^2, z_2z_3, z_2z_3^2)$. By Theorem \ref{wall}, we have 
$$M_\bbC = (M_\bbC)_\chi\oplus (M_\bbC)_{\bar{\chi}},$$
where $\chi(g)(x) = \eta_3x$. 
The characteristic polynomial of $g$ in $M_\bbC$ is equal to $(t^2+t+1)^3$. We have
$$\Or(M) = W(E_6)\rtimes (\bbZ/2\bbZ),$$
where the extra automorphism is defined by the symmetry of the Coxeter diagram. Since $g$ is of order 3, its  image $w$  in $\Or(E_6)$ belongs to $W(E_6)$. The classification of elements of order 3 in the Weyl group $W(E_6)$ shows  the conjugacy class of $w$  corresponds to the primitive embedding of lattice $A_2^3\hookrightarrow E_6$ so that $w$ acts as the product $c_1c_2c_3$ of Coxeter element in each copy of $A_2$. It is known that the centralizer of $w$ is a maximal subgroup of $W(E_6)$ of order 648 \footnote{Not to be confused with another maximal subgroup of $W(E_6)$ of the same order which is realized as the stabilizer subgroup of the sublattice $A_2^3$}.  This group is isomorphic to the unitary complex reflection group $L_3$ (No 25 in the list) and $(M_\bbC)_\chi$ is its three-dimensional reflection representation. 

This and other examples of appearance of finite unitary complex reflection groups as the  $G$-equivariant monodromy groups were first constructed by  V. Goryunov \cite{Go1}, \cite{Go2}.
\end{example}

Affine complex crystallographic reflection groups   can be also realized as $G$-equivariant monodromy groups. We give only one example, referring for more to \cite{Go3}, \cite{Go4}.

\begin{example} Consider a simple elliptic singularity of type $J_{10}$ with parameter $\lambda = 0$. Its equation is given in Table \ref{se}. Consider the symmetry of order $3$ defined by an automorphism $g:(z_1,z_2,z_3)\mapsto (z_1,\zeta_3z_2,z_3)$.  A monomial basis of the jacobian algebra is 
$(1,z_2,z_3,z_3^2,z_3^4,z_2z_3,z_2z_3^2,z_2z_3^4).$
We have $5$ invariant monomials $1,z_3,z_3^2,z_3^3,z_3^4$. Applying Theorem \ref{wall} we obtain that $M_\bbC = (M_\bbC)_\chi\oplus (M_\bbC)_{\bar{\chi}}$, both summands of dimension 5. The characteristic polynomial of $g$ is equal to $(1+t+t^2)^5$. Obviously, $g$ leaves $M^\perp$ invariant and the image $w =\bar{g}$ of $g$ in $\Or(M/M^) = \Or(E_8) = W(E_8)$ has characteristic polynomial $(1+t+t^2)^4$. 
It follows from the classification of conjugacy classes in $W(E_8)$ that $w$ is the product of the 
Coxeter elements in the sublattice $A_2^4$ of $E_8$. Its centralizer is a subgroup of index 2 in the wreath product 
$\Sigma_3^4\wr \Sigma_4$ of order 31104. This group is a finite complex reflection group $L_4$ (No 32 in the list). The centralizer of $g$ is the unique complex crystallographic group in affine space of dimension 5 with linear part $L_4$.
\end{example}

\section{Complex ball quotients}
\subsection{Hypergeometric integrals} 
Let 
$S$ be an ordered set of $n+3$ distinct points  $z_1,\ldots,z_{n+3}$  in $\bbP^1(\bbC)$. We assume that $(z_{n+1},z_{n+2},z_{n+3}) = (0,1,\infty)$.  Let $U = \bbP^1\setminus S$ and $\gamma_1,\ldots,\gamma_{n+3}$ be its standard generators of $\pi_1(U;u_0)$ satisfying the relation $\gamma_1\cdots \gamma_{n+3} = 1$.  We have a canonical surjection of the fundamental group of $U$ to the group $A = (\bbZ/d\bbZ)^{n+3}/\Delta(\bbZ/d\bbZ)$ which defines an \'etale covering $V\to U$ with the Galois group $A$. The open Riemann surface $V$ extends $A$-equivariantly to  a compact Riemann surface $X(z)$ with quotient $X(z)/A$ isomorphic to $\bbP^1(\bbC)$. 

Let $\boldsymbol{\mu} = (m_1/d,\ldots,m_{n+3}/d)$ be a collection of rational numbers in the interval $(0,1)$ satisfying
\begin{equation}\label{k}
\frac{1}{d}\sum_{i=1}^{n+3} m_i= k\in \bbZ.\end{equation}
They define a surjective homomorphism
$$\chi: A\to  \bbC^*, \quad \bar{\gamma}_i\mapsto e^{2\pi \sqrt{-1}m_i/d},$$
where $\bar{\gamma}_i$ is the image of $\gamma_i$ in $A$.

The following  computation can be found in \cite{DM} (see also \cite{Do}). 

\begin{lemma} Let $H^1(X(z),\bbC)_\chi$ be the $\chi$-eigensubspace of the natural representation of the Galois group $A$ on $H^1(X(z),\bbC)$. Then
$$\dim H^1(X(z),\bbC)_\chi = n+1.$$
Let $\Omega(X(z))$ be the space of holomorphic 1-forms on $X(z)$ and $\Omega(X(z))_\chi$ be the $\chi$-eigensubspace of $A$ in its natural action on the space $\Omega(X(z))$. Then 
$$\dim \Omega(X(z))_\chi = k-1,$$
where $k$ is defined in \eqref{k}.
\end{lemma}

Recall that $H^1(X(z),\bbC) = \Omega(X(z))\oplus \bar{\Omega}(X(z))$, so we get
$$H^1(X(z),\bbC)_\chi = \Omega(X(z))_\chi\oplus \bar{\Omega}(X(z))_\chi$$
and we can consider the hermitian form in $H^1(X(z),\bbC)_\chi$ induced by the skew-symmetric 
cup-product
$$H^1(X(z),\bbC) \times H^1(X(z),\bbC) \to H^2(X(z),\bbC) \cong \bbC.$$

The signature of the hermitiain form is equal to  $(k-1,n-k+2)$.

Now let us start vary the points $z_1,\ldots,z_n$ in $\bbP^1(\bbC)$ but keeping them distinct and not equal to $0,1$ or $\infty$. Let $\calU \subset \bbP^1(\bbC)^n$ be the corresponding set of parameters. 
Its complement in $(\bbP^1(\bbC)^n)$ consists of $N = \binom{n}{2}+3n$ hyperplanes $H_{ij}:z_i -z_j = 0$ and 
$H_i(0):z_i =0, H_i(1):z_i = 1, H_i(\infty): z_i = \infty$. Fix a point $z^{(0)}\in \calU$. For each of these hyperplanes $H$ consider a path  which starts at $z^{(0)}$ goes to a point on a small circle normal bundle of an open subset of $H$, goes along the circle, and then returns to the starting point. The homotopy classes $s_1,\ldots,s_N$ of these paths generate $\pi_1(\calU;z^{(0})$.

It is not difficult to construct a fibration over $\calU$ whose fibres are the curves $X(z)$. This defines a local coefficient system $\calH(\chi)$ over $\calU$  whose fibres are the spaces $H^1(X(z),\bbC)_\chi$. and the monodromy map
$$\pi_1(\calU;z^{(0})\to \U(H^1(X(z^{(0)}),\bbC)_\chi)$$
Denote the monodromy group by $\Gamma(\boldsymbol{\mu})$.

The important case for us is when $|\boldsymbol{\mu}| = 2$. In this case the signature is $(1,n)$ and we can consider the image of the monodromy group $\Gamma(\boldsymbol{\mu})$ in $\bbP U(1,n))$
which acts in the complex hyperbolic space $H_\bbC^n$. 

Here is the main theorem from \cite{DM}, \cite{Mo}.

\begin{theorem} The image of each generator $s_i$ of $\pi_1(\calU;z^{(0})$ in $\Gamma(\boldsymbol{\mu})$ acts as a complex reflection  in the hyperbolic space $H_\bbC^n$. The group $\Gamma(\boldsymbol{\mu})$ is a crystallographic reflection group  in $H_\bbC^n$ if and only if one of the following conditions is satisfied
\begin{itemize}
\item $(1-\frac{m_i}{d}-\frac{m_j}{d})^{-1}\in \bbZ,\  i\ne j, m_i+m_j < 1$;
\item $2(1-\frac{m_i}{d}-\frac{m_j}{d})^{-1}\in \bbZ,$ if  $m_i=m_j, i\ne j.$
\end{itemize}
\end{theorem}

All possible $\boldsymbol{\mu}$ satisfying the conditions  from the theorem can be enumerated. We have 59 cases if $n = 2$,  20 cases if $n = 3$, 10 cases if $n = 4$, 6 cases when $n = 5$, 3 cases if $n = 6$, 2 cases when $n = 7$, and 1 case if $n = 8$ or $n = 9$.

There are several cases when the monodromy group is cocompact. It does not happen in dimension $n > 7$. 

The orbit spaces $H_\bbC^n/\Gamma(\boldsymbol{\mu})$ of finite volume  have a moduli theoretical interpretation. It is isomorphic to the geometric invariant theory quotient 
$(\bbP^1)^{n+3}/\!/\SL(2)$ with respect to an appropriate choice of linearization of the action.

We refer to Mostow's survey paper \cite{Mo2}, where he explains a relation between the monodromy groups $\Gamma(\boldsymbol{\mu})$ to the monodromy groups of hypergeometric integrals.

\subsection{Moduli space of Del Pezzo surfaces as complex ball quotients}
In the last section we will discuss some recent work on complex ball uniformization of some moduli spaces in algebraic geometry. 

It is well-known that  a nonsingular  cubic curve in the projective plane
is isomorphic as a complex manifold to a complex torus $\bbC/\bbZ+\bbZ\tau$, where $\tau$ belongs to the upper-half plane $\calH  = \{a+bi\in \bbC:b > 0\}$. Two such tori  are isomorphic if and only if the corresponding $\tau$'s belong to the same orbit of the group $\Gamma = \SL(2,\bbZ)$ which acts on the $\calH$ by Moebius transformations $z\mapsto (az+b)/(cz+d)$. This result  implies that the moduli space of plane cubic curves is isomorphic to the orbit space $\calH/\Gamma$. Of course, the upper-half  plane is a model of the one-dimensional complex hyperbolic space $H_\bbC^1$ and the group $\Gamma$ acts as a crystallographic reflection group. 

In a beautiful paper of D. Allcock, J. Carlson and D. Toledo \cite{ACT}, the complex ball uniformization of the moduli space of plane cubics is generalized to the case of the moduli space of cubic surfaces in $\bbP^3(\bbC)$. It is known since the last century that the linear space $V$ of homogeneous forms of degree 3 in 4 variables admits a natural action of the group $\SL(4)$ such that the algebra of invariant polynomial functions on $V$  of  degree divisible by $8$ is freely generated by invariants $I_8, I_{16}, I_{24}, I_{32}, I_{40}$ of degrees indicated by the subscript. This can be interpreted as saying that the moduli space of nonsingular cubic surfaces admits a compactification isomorphic to the weighted projective space $ \bbP(1,2,3,4,5)$. Let 
$F(T_0,T_1,T_2,T_3) = 0$ be an equation of a nonsingular cubic surface $S$. Adding the cube of a new variable $T_4$ we obtain an equation 
$$F(T_0,T_1,T_2,T_3)+T_4^3 = 0$$
of a nonsingular cubic hypersurface $X$ in $\bbP^4(\bbC)$. There is a construction of an abelian variety  of dimension 10 attached to $X$ (the \emph{intermediate jacobian}) $\Jac(X)$. The variety $X$ admits an obvious automorphism of order 3 defined by  multiplying the last coordinate by a third root of unity. This  makes $\Jac(X)$ a principally polarized  abelian variety of dimension 10 with complex multiplication of certain type\footnote{This beautiful idea of assigning to a cubic surface a certain abelian variety was independently suggested by B. van Geemen and B. Hunt.}. The moduli space of such varieties is known to be isomorphic to a quotient of a 4-dimensional complex ball by certain discrete subgroup $\Gamma$.  It is proven in \cite{ACT} that the group $\Gamma$ is a hyperbolic complex crystallographic reflection group  and the quotient  $H_\bbC^4/\Gamma$ is isomorphic to the moduli space of  cubic surfaces with at most  ordinary double points as  singularities. By adding one point one obtains a compactification of the moduli space isomorphic to the weighted projective space $\bbP(1,2,3,4,5)$.

 The geometric interpretation of reflecting hyperplanes is also very nice, they form one orbit representing singular surfaces. The group $\Gamma$ contains a normal subgroup $\Gamma'$ with quotient isomorphic to the Weyl group $W(E_6)$. The quotient subgroup $H_\bbC^4/\Gamma'$ is the moduli space of marked nodal cubic surfaces. For a nonsingular  surface a marking is a fixing of order on the set of 27 lines on the surface. 

We mentioned  before that some complex ball quotients appear as the  moduli space of K3 surfaces which admit an action of a  cyclic group $G$ with fixed structure of the sublattice $(S_X)^G$.
This idea was used by S. Kond\=o to construct an action of a crystallographic reflection group $\Gamma$ in a complex ball $H_\bbC^6$ (resp. $H_\bbC^9$) with orbit space containing the moduli space of nonsingular plane quartic curves of genus 3\footnote{It is isomorphic  to the moduli space of  Del Pezzo surfaces of degree 2} (resp. moduli space of canonical curves of genus 4). In the first case he assigns to a plane quartic $F(T_0,T_1,T_2) = 0$ the quartic K3-surface 
$$F(T_0,T_1,T_2) +T_3^4 = 0$$
with automorphism of order 4 and $S_X^G \cong U(2)\perp A_1(-1)^6$. This leads to a new example of a crystallographic reflection group in $H_\bbC^6$. It is known that a canonical curve $C$ of genus 4 is isomorphic to a complete intersection of a quadric and cubic in $\bbP^3(\bbC)$. To each such curve Kondo assigns the K3 surface isomorphic to the triple cover of the quadric branched long the curve $C$. It has an action of a cyclic group of order 3 with $(S_X)^G \cong U\perp A_2(-2)$.  The stabilizer of one of a reflection hyperplane gives a complex reflection group in $H_\bbC^8$ with quotient isomorphic to a partial compactification of the moduli space of Del Pezzo surfaces of degree 1. Independently such a construction was found in \cite{HL}. 

Finally, one can also reprove the result of Allcock-Carlson-Toledo by using K3 surfaces instead of intermediate jacobians (see \cite{DGK}).

\begin{remark} All reflection groups arising in these complex ball uniformization constructions are not contained in the Deligne-Mostow list (corrected in \cite{Th}). However, some of them  are commensurable\footnote{this means that the two groups share a common subgroup of finite index.} with some groups from the list. For example, the group associated to cubic surfaces is commensurable to the group $\Gamma(\boldsymbol{\mu})$, where $d = 6, m_1 = m_2 = 1, m_3 = \ldots = m_7 =2$. We refer to  \cite{CHL} for a construction of complex reflection subgroups of finite volume which are not  commensurable to the groups from the Deligne-Mostow list. No algebraic-geometrical interpretation of these groups is known so far. \end{remark}

\bibliographystyle{amsplain}

\begin{thebibliography}{110}

\bibitem {Al}  D. Allcock, \textit{The Leech lattice and complex hyperbolic reflections}, Invent.  Math. J. {\bf 140} (2000), 283--31. MR1756997 (2002b:11091)


\bibitem {ACT}  D. Allcock, J.  Carlson,  and D. Toledo, \textit{The complex hyperbolic geometry of the moduli space of cubic surfaces}, J. Alg. Geom. \textbf{11} (2002), 659--724. MR1910264 (2003m:32011)

\bibitem {All} D. Allcock,  \textit{A monstrous proposal},  9 pages, math.GR/0606043.

\bibitem {ACT2}  D. Allcock, J.  Carlson,  and D. Toledo, \textit{The moduli space of cubic threefolds as a ball quotient}, math.AG/0608287.

\bibitem {An} E. Andreev, \textit{
Convex polyhedra of finite volume in Lobachevskii space}, (Russian) 
Mat. Sb. (N.S.) 83 (125) 1970 256--260. MR0273510 (42 \#8388)

\bibitem {Arnold} V. Arnol'd, \textit{Critical points of functions on a manifold with boundary, the simple Lie groups $B_k, C_k, F_4$ and sequences of evolutes}, Uspkhi Mat. Nauk {\bf 33} (1978), 91--105. MR0511883.

\bibitem {Ar} V. Arnol'd, S.  Gusein-Zade and A. Varchenko, \textit{
Singularities of differentiable maps. Vol. I.}
Translated from the Russian. Monographs in Mathematics, 82. 
BirkhŠuser Boston, Inc., Boston, MA, 1985. MR0777682 (86f:58018) 

\bibitem {BK} E. Bedford, K. Kim, \textit{Dynamics of rational surface automorphisms: linear fractional recurrences}, math.DS/0611297.

\bibitem {BP}  W. Barth and C.  Peters,  \textit{Automorphisms of Enriques surfaces}, Invent. Math. {\bf 73} (1983), 383--411. MR0718937 (85g:14052)

\bibitem {BPV}  W. Barth, K. Hulek, C. Peters, and A. Van de Ven, \textit{Compact complex surfaces} Second edition. Ergeb. der Mathematik und ihrer Grenzgebiete,  Springer-Verlag, Berlin, 2004. MR2030225 (2004m:14070)

\bibitem {BS} I. 
Bernstein and O. Shvarzman, \textit{ 
Chevalley's theorem for complex crystallographic Coxeter groups} (Russian) 
Funktsional. Anal. i Prilozhen. {\bf 12} (1978), no. 4, 79--80. MR0515632 (80d:32007) 

\bibitem {Bo1} R. Borcherds, \textit{Automorphism groups of Lorentzian lattices}, J. Algebra {\bf 111} (1987), 133--153. MR0913200 (89b:20018) 



\bibitem {Bo3}  R. Borcherds, \textit{Coxeter groups, Lorentzian lattices, and $K3$ surfaces},  Internat. Math. Res. Notices 1998, {\bf 19}, 1011--1031. MR1654763 (2000a:20088)

\bibitem  {Bou} N. Bourbaki, \textit{Lie groups and Lie algebras, Chapters 4--6}, Translated from the 1968 French original, Elements of Mathematics (Berlin). Springer-Verlag, Berlin, 2002.  MR1890629 (2003a:17001) 


\bibitem  {Br2} E. Brieskorn, \textit{Die Aufl\"osung der rationalen SingularitŠten holomorpher Abbildungen}, Math. Ann. {\bf 178} (1968), 255--270. MR0233819 (38 \#2140)

 \bibitem  {Br} E.  Brieskorn, \textit{Singular elements of semi-simple algebraic groups}, Actes du Congrs International des MathŽmaticiens (Nice, 1970), Tome 2, pp. 279--284. Gauthier-Villars, Paris, 1971. MR0437798 (55 \#10720)


\bibitem  {Bries} E. Brieskorn, \textit{Die Milnorgitter der exzeptionellen unimodularen SingularitŠten},  Bonner Mathematische Schriften [Bonn Mathematical Publications], 150. UniversitŠt Bonn, Mathematisches Institut, Bonn, 1983. MR0733785 (85k:32014)

 \bibitem {Bu} V. Bugaenko, \textit{Arithmetic crystallographic groups generated by reflections, and reflective hyperbolic lattices},  in `` Lie groups, their discrete subgroups, and invariant theory'', pp. 33--55, Adv. Soviet Math., 8, Amer. Math. Soc., Providence, RI, 1992. MR1155663 (93g:20094)




\bibitem {Ch} C. Chevalley, \textit{Invariants of finite groups generated by reflections},
Amer. J. Math. {\bf 77} (1955), 778--782. MR0072877 (17,345d) 

\bibitem {Coble} A. Coble, \textit{The ten nodes of the rational sextic and of the Cayley symmetroid},
Amer J. Math. {\bf 41} (1919), 243--265.

\bibitem  {Co}  A. Coble, \textit{Algebraic geometry and theta functions} (reprint of the 1929 edition), A. M. S. Coll. Publ., v. 10. A. M. S., Providence, R.I., 1982. MR0733252 (84m:14001)

\bibitem {Con} J. 
Conway, \textit{The automorphism group of the $26$-dimensional even unimodular Lorentzian lattice}, 
J. Algebra, {\bf 80} (1983),  159--163. MR0690711 (85k:11030) 
 
 \bibitem{CS} J. H. Conway and N.J.A. Sloane, \text{Sphere packings, lattices, and groups},  Grundlehren der Mathematischen Wissenschaften, 290. Springer-Verlag, New York, 1999.  MR1662447 (2000b:11077)

\bibitem{CD} F.  Cossec and I. Dolgachev, \emph{On automorphisms of nodal Enriques surfaces}, Bull. Amer. Math. Soc. (N.S.) {\bf 12} (1985), 247--249. MR0776478 (86f:14028)

 \bibitem {CHL} W. Couwenberg, G.  Heckman, E.  Looijenga, \textit{Geometric structures on the complement of a projective arrangement},  Publ. Math. Inst. Hautes ƒtudes Sci. No. 101 (2005), 69--161. MR2217047.

\bibitem{Coxeter2} H.S.M. Coxeter, \textit{The pure archimedean polytopes in six and seven dimensions}, Proc. Cambridge Phil. Soc., {\bf 24} (1928), 7--9.

\bibitem {Cox} H.S.M. Coxeter, \textit{Discrete groups generated by reflections}, Annals of Math. {\bf 35} (1934), 588--621.

\bibitem {DM} P. Deligne and G. Mostow, \textit{Monodromy of hypergeometric functions and nonlattice integral monodromy}, Inst. Hautes \'Etudes Sci. Publ. Math. \textbf{63} (1986), 5--89. MR0849651 (88a:22023a)

 
\bibitem {Dos} I. Dolgachev, 
\textit{Integral quadratic forms: applications to algebraic geometry (after V. Nikulin)}, 
Bourbaki seminar, Vol. 1982/83, 251--278, 
AstŽrisque, 105-106, 
Soc. Math. France, Paris, 1983. MR0728992 (85f:14036)

\bibitem {Do} I. Dolgachev, \textit{
On automorphisms of Enriques surfaces}, 
Invent. Math. {\bf 76} (1984), 163--177. MR0739632 (85j:14076)

\bibitem {Do2} I. Dolgachev, \textit{Infinite Coxeter groups and automorphisms of algebraic surfaces},  in ``The Lefschetz centennial conference, Part I (Mexico City, 1984)'', 91--106, Contemp. Math., {\bf 58}, Amer. Math. Soc., Providence, RI, 1986. MR0860406 (87j:14068)

\bibitem {DGK} I. Dolgachev,  B. van Geemen,  and S. Kond\=o, \textit{
A complex ball uniformization of the moduli space of cubic surfaces via periods of $K3$ surfaces}, 
J. Reine Angew. Math. \textbf{588} (2005), 99--148. MR2196731 (2006h:14051) 

\bibitem {DI} I. Dolgachev and V. Iskovskikh, \textit{Finite subgroups of the plane Cremona group}, math.AG/06510595. 

\bibitem {DoKe} I. Dolgachev and J. Keum, \textit{Birational automorphisms of quartic Hessian surfaces}, Trans. Amer. Math. Soc. {\bf 354} (2002), 3031--3057. MR1897389 (2003c:14045) 

\bibitem {DK} I. Dolgachev and  S. Kond\=o, \textit{A supersingular $K3$ surface in characteristic 2 and the Leech lattice},  Int. Math. Res. Not. {\bf 2003},  (2003) 1--23. MR1935564 (2003i:14051)

\bibitem {DO} I. Dolgachev and D. Ortland, \textit{Point sets in projective spaces and theta functions}, 
Ast\'erisque No. 165 (1988), 210 pp. (1989). MR1007155 (90i:14009)

\bibitem{DV} P. Du Val, \textit{On singularities which do not affect the conditions of adjunction}, Proc. Cambridge Phil. Society, {\bf 30} (1934), 434--465.

\bibitem{DuVal2}  P. Du Val, \emph{On the Kantor group of a set of points in a plane}, Proc. London Math. Soc. {\bf 42} (1936), 18--51.

\bibitem{DV2}  P. Du Val, \textit{Crystallography and Cremona transformations}, in ``The geometric vein'', 
pp. 191--201, Springer, New York-Berlin, 1981. MR0661778 (84h:52013)


\bibitem{Ebeling1} W. Ebeling, 
\textit{On the monodromy groups of singularities},  
Singularities, Part 1 (Arcata, Calif., 1981), 327Ð336, Proc. Sympos. Pure Math., 40, Amer. Math. 
Soc., Providence, RI, 1983. MR713071 (85d:14003

\bibitem{Ebeling2} W. Ebeling, \textit{The monodromy groups of isolated singularities of complete intersections} 
Lect.  Notes in Math., vol. 1293, Springer-Verlag, Berlin, 1987. MR0923114 (89d:32051)

\bibitem {Es} F. Esselmann, \textit{\"Uber die maximale Dimension von Lorentz-Gittern mit 
coendlicher Spiegelungsgruppe},  J. Num. Theory {\bf  61} (1996), 103--144. MR1418323 (97g:11073)

\bibitem {Fano} G. Fano, \textit{Superficie algebriche di genere zero e bigenere uno, e loro casi particolari},
Palermo Rend. {\bf 29} (1910), 98-118. 


\bibitem {FM} R. Friedman, J. Morgan and E. Witten, \textit{
Principal $G$-bundles over elliptic curves}, 
Math. Res. Lett. {\bf 5} (1998),  97--118. MR1618343 (99j:14037) 

\bibitem {Ga1} A.  Gabrielov, \text{Dynkin diagrams of unimodal singularities} (Russian) Funkcional. Anal. i Prilo\v zen. {\bf 8} (1974), no. 3, 1--6. MR0367274 (51 \#3516)

\bibitem {Gi}  M. Gizatullin, \textit{Rational $G$-surfaces}, (Russian) Izv. Akad. Nauk SSSR Ser. Mat. {\bf 44} (1980),  110--144, 239. MR0563788 (81d:14020)


\bibitem {Gi2} M. Gizatullin, 
\textit{The decomposition, inertia and ramification groups in birational geometry}, in `` Algebraic geometry and its applications'', eds. A. Tikhomirov and A. Tyurin, Aspects Math., E25, Vieweg, Braunschweig, 1994, pp. 39--45. MR1282018 (95d:14018) 

\bibitem  {Go1} V. Goryunov, \textit{Unitary reflection groups and automorphisms of simple hypersurface singularities}, in  ``New developments in singularity theory (Cambridge, 2000)'', 305--328, NATO Sci. Ser. II Math. Phys. Chem., 21, Kluwer Acad. Publ., Dordrecht, 2001. MR1849314 (2002m:32046)

\bibitem {Go2}  V. Goryunov, \textit{Unitary reflection groups associated with singularities of functions with cyclic symmetry} (Russian) Uspekhi Mat. Nauk {\bf 54} (1999), 3--24; [English transl.: Russian Math. Surveys {\bf 54} (1999),  873--893]. MR1741660 (2001b:32055)

\bibitem {Go3} V. Goryunov and S. Man, 
\textit{The complex crystallographic groups and symmetries of $J_{10}$}, 2004, preprint.

\bibitem {Go4} V. Goryunov, 
\textit{Symmetric $X_9$ singularities and the complex affine reflection groups}, 2005, preprint.

\bibitem{GN} V. Gritsenko and V. Nikulin, \textit{On classification of Lorentzian Kac-Moody algebras}, Russian Math. Surveys, {\bf 57:5} (2002), 921-979. MR1992083 (2004f:17034)

\bibitem {GZ}  S. Gusein-Zade, \textit{Monodromy groups of isolated singularities of hypersurfaces} (Russian), Uspehi Mat. Nauk {\bf 32} (1977), (194), 23--65, 263. MR0476738 (57 \#16295)

\bibitem {Harb} 
B. Harbourne, \textit{
Rational surfaces with infinite automorphism group and no antipluricanonical curve}, 
Proc. Amer. Math. Soc. {\bf 99} (1987),  409--414. MR0875372 (88a:14043) 

\bibitem {Ha}  R. Hartshorne, \textit{Algebraic geometry}, Graduate Texts in Mathematics, No. 52. Springer-Verlag, New York-Heidelberg, 1977. MR0463157 (57 \#3116)


\bibitem {Hi} A. Hirschowitz, \textit{
Sym\'etries des surfaces rationnelles g\'en\'eriques}, 
Math. Ann. {\bf 281} (1988), no. 2, 255--261. MR0949832 (89k:14064) 

\bibitem {HL} G. Heckman and E.  Looijenga, \textit{
The moduli space of rational elliptic surfaces},
Algebraic geometry 2000, Azumino (Hotaka), 185--248, 
Adv. Stud. Pure Math., 36, 
Math. Soc. Japan, Tokyo, 2002. MR1971517 (2004c:14068) 

\bibitem  {Hu} J. Humphreys, \textit{Reflection groups and Coxeter groups},  Cambridge Studies in Advanced Mathematics, 29. Cambridge University Press, Cambridge, 1990 MR1066460 (92h:20002)

\bibitem {Hunt}  B. Hunt, \textit{The geometry of some special arithmetic quotients},  Lecture Notes in Mathematics, 1637. Springer-Verlag, Berlin, 1996. MR1438547 (98c:14033)
 

\bibitem {Igusa} J.  Igusa, 
\textit{On the structure of a certain class of Kaehler varieties} 
Amer. J. Math. 76, (1954), 669--678. MR0063740 (16,172c)

\bibitem {Iv}  A.  Ivanov, \textit{A geometric characterization of the Monster}, in `` Groups, combinatorics \& geometry (Durham, 1990)'', 46--62, London Math. Soc. Lecture Note Ser., 165, Cambridge Univ. Press, Cambridge, 1992. MR1200249 (94c:20033)

\bibitem{Kane}  R. Kane, \emph{Reflection groups and invariant theory}, CMS Books in Mathematics/Ouvrages de MathŽmatiques de la SMC, 5. Springer-Verlag, New York, 20. MR1838580 (2002c:20061)

\bibitem  {Ka} S.  Kantor, \textit{Theorie der endlichen Gruppen von eindeutigen Transformationen in der Ebene},
Berlin. Mayer \& MŸller. 111 S. gr. $8^\circ$. 1895.

\bibitem {KK}  J. Keum and S. Kond\=o, \textit{The automorphism groups of Kummer surfaces associated with the product of two elliptic curves}, Trans. Amer. Math. Soc. {\bf 353} (2001), 1469--1487.  MR1806732 (2001k:14075)

\bibitem {Kh1} A. Khovanskii, \textit{
Hyperplane sections of polyhedra, toric varieties and discrete groups in Lobachevskii space} (Russian) 
Funktsional. Anal. i Prilozhen.  {\bf 20} (1986), no. 1, 50--61, 96. MR0831049 (87k:22015) 

\bibitem 
{Kh2} A. Khovanskii, \textit{
Combinatorics of sections of polytopes and Coxeter groups in Lobachevsky spaces},
The Coxeter legacy, 129--157, 
Amer. Math. Soc., Providence, RI, 2006. MR2209026 


\bibitem {Ko} M. Koitabashi, \textit{
Automorphism groups of generic rational surfaces}, 
J. Algebra {\bf 116} (1988), no. 1, 130--142. MR0944150 (89f:14045) 

\bibitem {Kond}  S. Kond\=o, \textit{Enriques surfaces with finite automorphism groups},  Japan. J. Math. (N.S.) {\bf 12} (1986), 191--282. MR0914299 (89c:14058)

\bibitem{KondoNew} S. Kond\=o, \textit{Automorphisms of algebraic K3 surfaces which act trivially on Picard groups}, J. Math. Soc. Japan {\bf 44} (1992), 75--98. MR1139659 (93e:14046)

 \bibitem {Kondo} S. Kond\=o, \textit{The automorphism group of a generic Jacobian Kummer surface},  J. Algebraic Geom. {\bf 7} (1998),  589--609. MR1618132 (99i:14043)

\bibitem  {Ko1} S. Kond\=o, \textit{The moduli space of curves of genus 4 and Deligne-Mostow's complex reflection groups},  Algebraic geometry 2000, Azumino (Hotaka), 383--400, Adv. Stud. Pure Math., 36, Math. Soc. Japan, Tokyo, 2002. MR1971521 (2004h:14033)

\bibitem {Ko2}  S. Kond\=o, \textit{A complex hyperbolic structure for the moduli space of curves of genus three}, J. Reine Angew. Math. \textbf{525} (2000), 219--232. MR1780433 (2001j:14039)

\bibitem{Li} J.  Lipman, \textit{Rational singularities, with applications to algebraic surfaces and unique factorization}, Inst. Hautes ƒtudes Sci. Publ. Math. No. 36 (1969) 195--279. MR0276239 (43 \#1986)

 \bibitem {Lo1}  E. Looijenga, \textit{Homogeneous spaces associated to certain semi-universal deformations}, Proceedings of the International Congress of Mathematicians (Helsinki, 1978), pp. 529--536, Acad. Sci. Fennica, Helsinki, 1980.  MR0562651 (81j:14006)

\bibitem {Lo2} E. 
Looijenga, \textit{Invariant theory for generalized root systems}, 
Invent. Math. {\bf 61} (1980), no. 1, 1--32. MR0587331 (82f:17011) 

\bibitem {Lo3}  E. Looijenga, \textit{Isolated singular points on complete intersections}, London Mathematical Society Lecture Note Series, 77. Cambridge University Press, Cambridge, 1984. MR0747303 (86a:32021)


\bibitem{Lo4} E. Looijenga and R. Swiersa, \textit{The period map for cubic threefolds}, math.AG/0608279.

\bibitem {Mc} C. McMullen, \textit{Dynamics of blowups of the projective plane}, 2006, preprint.

\bibitem {Milnor}  J. Milnor, \textit{Singular points of complex hypersurfaces},  Annals of Mathematics Studies, No. 61 Princeton University Press, Princeton, N.J.; University of Tokyo Press, Tokyo 1968. MR0239612 (39 \#969)

\bibitem  {Mo} G. Mostow,\textit{Generalized Picard lattices arising from half-integral conditions}, 
Inst. Hautes ƒtudes Sci. Publ. Math. \textbf{63} (1986), 91--106. MR0849652 (88a:22023b) 


\bibitem {Mo2} G. Mostow, 
\textit{Braids, hypergeometric functions, and lattices}, 
Bull. Amer. Math. Soc. (N.S.) {\bf 16} (1987), no. 2, 225--246. MR0876959 (88e:22017) 

\bibitem {Mu} S. Mukai, \textit{Geometric realization of $T$-shaped root systems and counterexamples to Hilbert's fourteenth problem}, in `` Algebraic transformation groups and algebraic varieties'', 123--129, Encyclopaedia Math. Sci., 132, Springer, Berlin, 2004. MR2090672 (2005h:13008)

\bibitem {Na}  Y. Namikawa, \textit{Periods of Enriques surfaces}, Math. Ann. {\bf 270} (1985),  201--222. MR0771979 (86j:14035)

 
 \bibitem {Ni2}   V. Nikulin, \textit{Quotient-groups of groups of automorphisms of hyperbolic forms by subgroups generated by $2$-reflections. Algebro-geometric applications. Current problems in mathematics} Vol. 18, pp. 3--114, Akad. Nauk SSSR, Vsesoyuz. Inst. Nauchn. i Tekhn. Informatsii, Moscow, 1981. MR0633160 (83c:10030) [English translation: J. Soviet Math., {\bf 22} (1983), 1401--1475].
 
 \bibitem {Ni3} V. 
Nikulin, \textit{Description of automorphism groups of Enriques surfaces} (Russian) 
Dokl. Akad. Nauk SSSR {\bf 277} (1984),  1324--1327. MR0760514 (86c:14033) 


\bibitem{Ni4} V. Nikulin, \emph{
$K3$ surfaces with a finite group of automorphisms and a Picard group of rank three} (Russian) 
in ``Algebraic geometry and its applications'',  
Trudy Mat. Inst. Steklov. {\bf 165} (1984), 119--142. MR0752938 (86e:14018) 
 
 
 \bibitem {Ni1} V. Nikulin, \textit{Discrete reflection groups in Lobachevsky spaces and algebraic surfaces}, Proceedings of the International Congress of Mathematicians, Vol. 1, 2 (Berkeley, Calif., 1986), 654--671, Amer. Math. Soc., Providence, RI, 1987. MR0934268 (89d:11032)


\bibitem {No} S. Norton, \text{
Constructing the Monster}, in `` 
Groups, combinatorics \& geometry (Durham, 1990)'', 63--76, 
London Math. Soc. Lecture Note Ser., 165, 
Cambridge Univ. Press, Cambridge, 1992.  MR1200250 (94c:20034) 

\bibitem {OS} P. Orlik and L. Solomon, \textit{Arrangements defined by unitary reflection groups}, Math. Ann. {\bf 261} (1982), no. 3, 339--357. MR0679795 (84h:14006) 

\bibitem {OS2} P. Orlik, H. Terado, \textit{Arrangements of hyperplanes}, Springer-Verlag, 1992, MR1217488.

\bibitem {Pi1}  H. Pinkham, \textit{Simple elliptic singularities, Del Pezzo surfaces and Cremona transformations}, in ``Several complex variables (Proc. Sympos. Pure Math., Vol. XXX, Part 1, Williams Coll., Williamstown, Mass., 1975)'', pp. 69--71. Amer. Math. Soc., Providence, R. I., 1977 MR0441969 (56 \#358)

\bibitem {Pinkham} H. Pinkham, \textit{R\'esolution simultan\'ee de points doubles rationnels} in 
``S\'eminaire sur les Singularit\'es des Surfaces'',  Ed.  by M. Demazure, H. Pinkham and B. Teissier. Lecture Notes in Mathematics, 777. 
Springer, Berlin, 1980. MR0579026 (82d:14021) 


\bibitem {Po} V. Popov, \textit{Discrete complex reflection groups}, Communications of the Mathematical Institute, Rijksuniversiteit Utrecht, 15. Rijksuniversiteit Utrecht, Mathematical Institute, Utrecht, 1982. 89 pp. MR0645542 (83g:20049)









\bibitem {PS} I. Pjatecki-Shapiro and I.  Shafarevich, \textit{Torelli's theorem for algebraic surfaces of type ${\rm K}3$},  Izv. Akad. Nauk SSSR Ser. Mat. 35 1971 530--572. MR0284440 (44 \#1666) [English translation: Math. USSR Izv.,
{\bf 5} (1971), 547--587].

\bibitem {Pr} M.  Prokhorov, \textit{Absence of discrete groups of reflections with a noncompact fundamental polyhedron of finite volume in a Lobachevskii space of high dimension} (Russian),  Izv. Akad. Nauk SSSR Ser. Mat.  {\bf 50} (1986), no. 2, 413--424. MR0842588 (87k:22016)


\bibitem {RS}  A. Rudakov and I. Shafarevich, \textit{Surfaces of type $K3$ over fields of finite characteristic} (Russian) Current problems in mathematics, Vol. 18, pp. 115--207, Akad. Nauk SSSR, Vsesoyuz. Inst. Nauchn. i Tekhn. Informatsii, Moscow, 1981. MR0633161 (83c:14027) [English translation:  I. Shafarevich,  \textit{Collected mathematical papers}, Springer-Verlag, Berlin, 1989. MR0977275 (89m:01142)]

\bibitem {SW} 
R. Scharlau and C.  Walhorn, \textit{
Integral lattices and hyperbolic reflection groups},   
Journ\'ees Arithm\'etiques, 1991 (Geneva). 
Ast\'erisque No. 209, (1992),  279--291. MR1211022 (94j:11057) 


\bibitem{Schoute} P. Schoute, \textit{Over het verband tusschen de hoekpunten van een bepaald zesdimsionaal polytoop en de rechten van een kubisch oppervlak}, Proc. Konig. Akad. Wis. Amsterdam, {\bf 19} (1910), 375--383.

 \bibitem {Se} J.-P. Serre, \textit{A course in arithmetic}, Translated from the French. Graduate Texts in Mathematics, No. 7. Springer-Verlag, New York-Heidelberg, 1973. MR0344216 (49 \#8956)

\bibitem {Severi} F. Severi, \textit{Complementi alla teoria della base per la totalitˆ delle curve di una superficie algebrica}.
Palermo Rend. {\bf 30} (1910), 265--288.

\bibitem {Sha} I. Shafarevich, \textit{Le th\'eor\`eme de Torelli pour les surfaces alg\'ebriques de type K3}. Actes du Congrs International des MathŽmaticiens (Nice, 1970), Tome 1, pp. 413--417. Gauthier-Villars, Paris, 1971. MR0419459 (54 \#7480)


\bibitem {ST} G. Shephard and J. Todd,  
\textit{Finite unitary reflection groups}
Canadian J. Math. {\bf 6}, (1954). 274--304. MR0059914 (15,600b)

 \bibitem {Sh} O. Shvarzman, \textit{
Reflectivity of three-dimensional hyperbolic lattices}, (Russian) 
Problems in group theory and homological algebra (Russian), 135--141, 
Matematika, 
Yaroslavl Gos. Univ., Yaroslavl, 1990. MR1169974 (93d:51031) 

\bibitem {Si}  C. Siegel, \textit{
\"Uber die analytische Theorie der quadratischen Formen},
Ann. of Math. (2) {\bf 36} (1935), no. 3, 527--606.


\bibitem {Sl2} P.  Slodowy, \textit{Simple singularities and simple algebraic groups}, Lecture Notes in Mathematics, 815. Springer, Berlin, 1980. MR0584445 (82g:14037)

\bibitem {Sl3}   P. Slodowy, \textit{Four lectures on simple groups and singularities}, Communications of the Mathematical Institute, Rijksuniversiteit Utrecht, 11. Rijksuniversiteit Utrecht, Mathematical Institute, Utrecht, 1980. MR0563725 (82b:14002)

\bibitem {Sl1}  P. Slodowy, \textit{Simple singularities and complex reflections}, in  ``New developments in singularity theory (Cambridge, 2000)'', 329--348, NATO Sci. Ser. II Math. Phys. Chem., 21, Kluwer Acad. Publ., Dordrecht, 2001. MR1849315 (2002h:14061)

\bibitem {Steenbrink}  
J. Steenbrink, \textit{Mixed Hodge structures associated with isolated singularities}, 
in ``Singularities, Part 2 (Arcata, Calif., 1981)'', 513--536, 
Proc. Sympos. Pure Math., v. 40, 
Amer. Math. Soc., Providence, RI, 1983. MR0713277 (85d:32044) 




\bibitem {Te} T. Terada, \textit{ 
Probl\`eme de Riemann et fonctions automorphes provenant des fonctions hyperg\'eom\'etriques de plusieurs variables}, J. Math. Kyoto Univ. \textbf{13} (1973), 557--578. MR0481156 (58 \#1299) 

\bibitem  {Th} W. Thurston, \textit{Shapes of polyhedra and triangulations of the sphere},
The Epstein birthday schrift, 511--549 (electronic), 
Geom. Topol. Monogr., \textbf{1}, 
Geom. Topol. Publ., Coventry, 1998.  MR1668340 (2000b:57026) 
 
 \bibitem {Tj} G. Tjurina, \textit{Resolution of singularities of flat deformations of double rational points}, Funkcional. Anal. i Prilo\v zen. {\bf 4} 1970 no. 1, 77--83.  MR0267129 (42 \#2031) 


\bibitem {Vi1} E. Vinberg, \textit{
Discrete groups generated by reflections in Lobachevskii spaces}, (Russian)
Mat. Sb. (N.S.) {\bf 72} (114) 1967 471--488; correction, ibid. {\bf 73} (115) (1967), 303 MR0207853 (34 \#7667) 

\bibitem {Vi2} E. Vinberg, \textit{The groups of units of certain quadratic forms}, (Russian) 
Mat. Sb. (N.S.) {\bf 87}(129) (1972), 18--36. MR0295193 (45 \#4261) 


\bibitem {Vi3} E. Vinberg, \textit{Some arithmetical discrete groups in Lobachevsky spaces}. in ``Discrete subgroups of Lie groups and applications to moduli'', Internat. Colloq., Bombay, 1973, pp. 323--348. Oxford Univ. Press, Bombay, 1975. MR0422505 (54 \#10492)

 \bibitem {VK} E. Vinberg and I. Kaplinskaja, \textit{The groups $O_{18,1}(Z)$ and $O_{19,1}(Z)$},  (Russian) Dokl. Akad. Nauk SSSR 238 (1978), no. 6, 1273--1275. MR0476640 (57 \#16199)

\bibitem {Vi4} E. Vinberg, \textit{The two most algebraic $K3$ surfaces}, Math. Ann. {\bf 265} (1983), no. 1, 1--21. MR0719348 (85k:14020)

 \bibitem {Vi5} 
E. Vinberg, \textit{
Discrete reflection groups in Lobachevsky spaces}, 
Proceedings of the International Congress of Mathematicians, Vol. 1, 2 (Warsaw, 1983), 593--601, 
PWN, Warsaw, 1984. MR0804716 (87h:22016) 


\bibitem {Vi6} 
E. Vinberg, \textit{
Absence of crystallographic groups of reflections in Lobachevski spaces of large dimension}, (Russian) 
Trudy Moskov. Mat. Obshch. 47 (1984), 68--102, 246. MR0774946 (86i:22020) 

\bibitem {V7} 
E. Vinberg, \textit{
Hyperbolic groups of reflections} (Russian) 
Uspekhi Mat. Nauk {\bf 40} (1985), no. 1(241), 29--66, 255. MR0783604 (86m:53059) [English transl.: Russian Math. Surveys {\bf 40} (1985), no. 1, 31--75)
 
 \bibitem {VS} E. Vinberg and O. Shvartsman, \textit{Discrete groups of motions of spaces of constant curvature}, Geometry, II, 139--248, Encyclopaedia Math. Sci., 29, Springer, Berlin, 1993.  MR1254933 (95b:53043)

\bibitem {Wall} 
C. T. C.  Wall, \textit{
A note on symmetry of singularities},  
Bull. London Math. Soc. {\bf 12} (1980), 169--175. MR0572095 (81f:32009) 

\end{thebibliography}

\end{document}